\documentclass[12pt]{article}
\usepackage{amssymb,amsfonts,amsmath,amsthm,mathrsfs,enumerate}
\usepackage{enumitem}

\usepackage{graphics,xcolor, colortbl,multirow, tcolorbox}

\usepackage{graphicx}
\usepackage{color}
\usepackage[colorlinks = true,
linkcolor = blue,
urlcolor  = blue,
citecolor = blue,
anchorcolor = blue]{hyperref}

\newtheorem{theorem}{Theorem}[section]
\newtheorem{lemma}[theorem]{Lemma}
\newtheorem{Corollary}[theorem]{Corollary}
\newtheorem{proposition}[theorem]{Proposition}
\newtheorem{remark}{Remark}[section]

\newtheorem{Assumption}{Assumption}[section]

\numberwithin{equation}{section}

\title{Dynamics of a class of extensible beams with degenerate and non-degenerate nonlocal damping\thanks{The authors wish to dedicate this paper to the memory of Professor Igor Chueshov.}}

\author{
	\small {\bf  E. H. Gomes Tavares\thanks{
	Postdoctoral fellow. Email: 
	\href{mailto:eduardogomes7107@gmail.com}{eduardogomes7107@gmail.com}.}}\\
	\small Department of Mathematics, State University of Londrina, \\
	\small 86057-970 Londrina, PR, Brazil. \\ 
 	\small {\bf M. A. Jorge Silva\thanks{Supported by the CNPq, grant  301116/2019-9. Email: 
	\href{mailto:marcioajs@uel.br}{marcioajs@uel.br}.  (Corresponding Author)}}\\
	\small Department of Mathematics, State University of Londrina, \\
	\small 86057-970 Londrina, PR, Brazil. \\
	\small {\bf V. Narciso}  \footnote{Supported by the FUNDECT, grant 219/2016. Email: 
	\href{mailto:vnarciso@uems.br}{vnarciso@uems.br}.}\\
	\small Nucleus of Exact and Technological Sciences,
	\small State University of Mato Grosso do Sul \\
	\small 79804-970 Dourados, MS, Brazil. \\
	\small {\bf A. Vicente}  \footnote{Email:
		\href{mailto:andre.vicente@unioeste.br}{andre.vicente@unioeste.br}.}\\
	\small Center of Exact and Technological Sciences,
	\small State University of Paraná West\\
	\small 85819-110 Cascavel, PR, Brazil. }
\date{}

\begin{document}

\maketitle

\begin{abstract}
This work is concerned with new results on long-time dynamics of a class of hyperbolic evolution equations related to extensible beams with three distinguished nonlocal  nonlinear  damping terms.  
In the first possibly degenerate case, the results feature the existence of a family of compact global attractors and a thickness estimate for their Kolmogorov's
$\varepsilon$-entropy. Then, in the non-degenerate context, the structure of the helpful nonlocal damping leads to the existence of finite-dimensional compact global and exponential attractors. Lastly, in a degenerate and critical framework, it is proved  the existence of a bounded closed global attractor but not compact. To the proofs, we provide several new technical results by means of refined estimates that open up perspectives for a new branch of nonlinearly damped problems.
\end{abstract}


\smallskip

\noindent{{\bf 2020 MSC:} 35B40; 35B41; 35L75; 37L25; 37L30.

{\footnotesize		
\begin{quote}
	\tableofcontents
\end{quote}	
}

\section{Introduction}\label{sec-intro}

In the present article, we address the following   evolution problem of  hyperbolic type  with nonlocal  nonlinear damping term
\begin{equation} \label{P-intro}
\left\{\begin{array}{l}\displaystyle{
	u_{tt}+\kappa A u +
	A_1 u+f(u)+k\big(\mathcal{E}_{\alpha}(u,u_t)\big)u_t= h_\lambda ,\quad t> 0,} \medskip \\
\displaystyle{u(0)=u_0, \ \ u_t(0)=u_1,}
\end{array}\right.
\end{equation} 
where $\kappa$ is non-negative  parameter, $A$ and $A_1$ are linear self-adjoint positive definite  operators related to Laplacian and bi-harmonic  differential operators, 
respectively, $f(u)$ corresponds to a nonlinear source of lower order and growth exponent $p$, 
$h_\lambda=\lambda h$ represents an external force  with $\lambda\in[0,1]$ and $h$ lying in a Hilbert space $(H,\|\cdot\|)$, and 
\begin{equation}\label{Ealpha}
\mathcal{E}_{\alpha}(u,u_t)=\|A^{\alpha} u\|^2+\|u_t\|^2, \ \  \alpha\in[0,1], 
\end{equation}
with fractional powers $A^\alpha$ to be well defined later. Our study encompasses the abstract model  \eqref{P-intro} including three possibilities concerning the scalar function $k(\cdot)$ as follows:
\begin{enumerate}[label=(\textbf{k.\arabic*})]
	
	\item\label{k1} $k(\cdot)$ is 
	a monomial-like function on $[0,\infty)$ with	exponent $ q\geq 1/2$, namely, 
$$
k(s)= \gamma s^{q}, \, s\geq0, \ \ \mbox{ with } \ \ \gamma>0;
$$

	\item\label{k2} $k(\cdot)$ is any  $C^1$-function on $[0,\infty)$ such that $k(s)>0, \, s\geq0$;
	
	\item\label{k3} $k(\cdot)$ is a bounded Lipschitz function on  $[0,\infty)$ such that $k\equiv0$ on $[0,1]$ and $k(s)$ is strictly increasing for $s>1.$
\end{enumerate}

Through the coming statements,  we try to be more transparent as possible in clarifying that problem \eqref{P-intro}-\eqref{Ealpha} with damping coefficient obeying the three possibilities \ref{k1}-\ref{k3} has not been studied in the literature, by reaching in each specific case a level of results concerning its long-time behavior. In what follows, we state our main achievements.

\subsection{The possibly degenerate case {\rm \ref{k1}}}

In such a case the ``polynomial'' structure of the function $k(\cdot)$ says that the nonlocal damping coefficient  can  degenerate
whenever the argument $\mathcal{E}_{\alpha}(u,u_t)$ vanishes, and the rate of degeneracy (of such unknown degenerate points)  is determined by  the  exponent $q$. Additionally, 
the long-time dynamics of the system 
 \eqref{P-intro}  is characterized by a ``polynomial behavior'' of type $1/q,$ which can very slow when $q$ is very large. 
Within this scenario, our main achievements   are stated in Sections \ref{sec-tech-results} and \ref{sec-dynamics} whose highlights are presented as follows: 
\begin{enumerate}
	\item We prove a new key inequality (see Proposition \ref{theo_inequality}) that 
	gives a suitable estimate for the difference of two
	trajectory solutions of \eqref{P-intro}. Moreover, it will be very useful in the proof of attractors (smooth properties) for critical aspects in terms of the exponent $p$. To its proof, we launch a generalized Nakao's inequality (see Proposition \ref{Nak-lemma}) whose proof is only given in Appendix  \ref{sec-proof-nak-lemma}. We point out that such a generalized Nakao's lemma will be also applicable to a wide class of autonomous and 
	nonautonomous dynamical systems in the near future.
	
\item  Theorem \ref{theo_main} brings out  the existence of a family of attractors for the dynamical system associated with problem \eqref{P-intro} as well as their geometric and continuity properties. All statements in this result are reached for 
critical parameters $p$ and subcritical power $\alpha\in[0,1)$ due to lack of compactness when $\alpha=1$ (which is the subject of the case \ref{k3}). To the proof of asymptotic smoothness/compactness properties, we have combined our new generalized Nakao's inequality with an extended version of Khanmamedov's limit (see Proposition \ref{Khan-lemma}). For
the continuity properties with respect to (w.r.t. for short)	parameter $\lambda\in[0,1]$, we still prove a new Lipschitz continuous result 
(Proposition \ref{prop-esp-conti})

\item In the subcritical framework w.r.t. both parameters $\alpha$ and $p$, we can compute an estimate of Kolmogorov's
$\varepsilon$-entropy of the family of attractors, see Theorem \ref{theo_mainIII}, which roughly speaking measures the ``thickness'' of the compact attractors. For such an achievement we have regarded  the stabilizability estimate provided by  Corollary $\ref{cor-key-subcritic}$ instead of the above mentioned key inequality and this reveals why we must go down to the subcritical case concerning the growth exponent $p$.

\end{enumerate}

As far as we know, the model \eqref{P-intro}-\eqref{Ealpha} in case 
\ref{k1}, as well as the above mentioned results,  
have not been addressed in the literature in the sense of dynamical systems. In what concerns the stability to equilibrium, closer  models we have found  are  \cite{cavalcante2-marcio-vando,jorge-narciso-vicente}. In  \cite{cavalcante2-marcio-vando} an interesting local stability result is presented with degenerate coefficient $\gamma\left[\|Au\|^2\right]^q u_t$, but only when initial data is taken regular and bounded, cf. \cite[Theorem 3.1]{cavalcante2-marcio-vando}, which is not viable in the study of dynamical systems in the pattern weak phase space. In \cite{jorge-narciso-vicente} the authors address the problem with 
nonlocal structural (strong) damping 
$\gamma\left[\|Au\|^2+\|u_t\|^2\right]^q Au_t$ stead of the nonlocal weak one
$\gamma\left[\|A^{\alpha} u\|^2+\|u_t\|^2\right]^qu_t$.
However, due to technical difficulties, 
in \cite[Theorems 2.1 and 3.1]{jorge-narciso-vicente} the authors only deal with existence and stability of regular solution when $q\geq1$. Here, we have surpassed such difficulties and all results are provided for   $q \geq \frac{1}{2}$. Moreover, it is worth pointing out  that  our particular  
Corollary \ref{P3}  (see also Remark \ref{rem-sque-polin}) gives the precise answer on  polynomial stability to the related homogeneous  problem, which finally clarifies the prediction estimate stated in \cite[Theorem 4.1]{jorge-narciso-vicente}.
We still stress that, in this case, all results encompass 
the particular one with nonlocal averaged damping 
$\gamma \|u_t\|^{2q} u_t$, that is, by neglecting the potential energy $\|A^{\alpha} u\|^2$ in \eqref{Ealpha}.

Finally, it is worth notifying  that we also put some strength in the attempts of proving finite 
dimensionality for the family of compact global attractors obtained in Theorem \ref{theo_main}. Indeed, throughout the whole Subsection \ref{sec-final-comments}, we  clarified that it seems a delicate task due to the 
 nature of such a damping term in case \ref{k1}. Our  conclusion is that we could not reach the assumptions of the current abstract results in dynamical systems   to prove the finiteness of dimension   in this case of nonlocal 
  nonlinear  possibly degenerate damping and a keener theory must be done for this purpose. All technical details concerning this issue are elucidated in Subsection \ref{sec-final-comments}.

 \subsection{The non-degenerate case {\rm \ref{k2}}} 
 
 It is for sure a more touchable case 
where we do can thrive up to the existence of exponential attractors.
As a matter of fact, such an assumption in \ref{k2} is already employed by the authors in \cite{jorge-narciso-DIE,jorge-narciso-DCDS,jorge-narciso-EECT} for related extensible beam models with nonlocal averaged damping coefficient 
$
k\big(\|A^{1/4} u\|^2 \big) 
$
instead of the energy damping coefficient $k\big(\|A^{\alpha} u\|^2+\|u_t\|^2\big)$.
 Nevertheless,  our results here improve and generalize those provided in these references.
Indeed, this case is treated in 
Section \ref{sec-dynamics-k2} and the highlights are:
\begin{enumerate}
	\item In Theorem \ref{theo_main-k2} we catch up the same results as in Theorem \ref{theo_main} by aggregating 
	finite fractal dimension and regularity of the existing global attractors, and also the existence of  generalized fractal exponential attractors, that is, exponential attractors whose fractal dimensional is finite only in an extended space. Moreover, thanks to a new stabilizability estimate in this case (see Proposition \ref{prop-stabi-ineq-k2}), 	
	all these properties in Theorem \ref{theo_main-k2} are achieved for the critical source exponent $p$ but still in the subcritical case w.r.t. fractional power $\alpha$.

	\item We also draw attention to the fact that positive constant functions are also incorporated in this case, and for such a very specific situation we can go further and prove the existence of time-dependent exponential attractors  (with finite fractal dimensional in the standard phase space). This is the subject of Theorem \ref{theo-exponential-attractor}
whose proof is mainly achieved by means of a new 	smoothing property (see Proposition \ref{prop-compac-embedd}) in the subcritical aspect w.r.t. $p$ and assuming the commutative case $A=A_1^{1/2}$.
\end{enumerate}

Therefore, the right above facts furnish a considerable extension of the previous results achieved by authors in \cite{jorge-narciso-DIE,jorge-narciso-DCDS,jorge-narciso-EECT} concerning the criticality of the exponent $p$ and the existence of exponential attractors in the linear damping case. Additionally, we note that standard examples of functions $k(\cdot)$ are:
$$
k(s)=\gamma e^{\pm s}, \ \ k(s)=\frac{\gamma}{1+s}, \  \ k(s)=\gamma, \ \ s\geq0, \ \gamma>0.
$$

\subsection{The degenerate case {\rm \ref{k3}}} 

This is exactly the case where we address the critical parameter $\alpha=1$ and it justifies why in the previous ones we only get compact global attractors for subcritical powers $\alpha\in[0,1).$ In fact, due to lack of compactness of the damping coefficient in the standard phase space and since $k(\cdot)$ vanishes on $[0,1]$, 
then  we shall prove in Section \ref{sec-dynamics-k3}  the existence of a noncompact global attractor even by assuming $f=h_\lambda=0$ in the model \eqref{P-intro}. 
Here, our main result is Theorem \ref{theo-noncompact-attract} that proves the existence of a closed bounded (but not compact) global attractor for $\kappa>0$ small enough, and due to the uniqueness of a global attractor (when it there exists), this  precludes the existence of a compact global attractor under the assumptions considered in the present third case.
We also point out that our main result is an extension of the 
one stated in \cite[Proposition 5.3.9]{chueshov-solo} where the particular case $\kappa=0$ is considered. To our best knowledge, the case $\kappa>0$ has never been approached and, although similar, it requires different computations.

Examples of functions $k(\cdot)$ in this case are given as follows.
\begin{equation*}
k(s)=\left\{\begin{array}{l}
0, \ \ 0 \leq s \leq 1, \smallskip \\
\gamma (1-s^{-1}), \ \  s>1,
\end{array}\right.
\ \  
k(s)=\left\{\begin{array}{l}
0, \ \  0 \leq s \leq 1, \smallskip  \\
\gamma (1-e^{-s}), \ \  s>1,
\end{array}\right.
\quad  \gamma>0.
\end{equation*}

\noindent
We finish the introduction with the organization of the  remaining paper:

\begin{itemize}
	\item[-]  in Section \ref{sec-physic-abstract} we provide a  physical motivation for studying the problem \eqref{P-intro} and  its abstract formulation through concrete problems in mathematical physics;
	
	\item[-] in Section \ref{sec-syst-din} we state the well-posedness of the model \eqref{P-intro}-\eqref{Ealpha} under the  three conditions \ref{k1}-\ref{k3} and set  the corresponding dynamical system;
	
	\item[-] From Section \ref{sec-tech-results} to Section  \ref{sec-dynamics-k3} we state and prove our main results concerning the long-time dynamics of problem \eqref{P-intro} as well as it is exhibited additional remarks on the novelties introduced in the present article  comparing to previous literature;
	
	\item[-]  In  Appendix \ref{sec-appendix} we present some supplementary technical proofs;  
	
	\item[-] Last, aiming  reader's convenience and also for the sake of  self-containment, in Appendix \ref{sec-appendixB} we remind several concepts and results coming from the meaningful  literature  in dynamical systems.

\end{itemize}

\section{Physical and mathematical patterns}\label{sec-physic-abstract}

In this section, our main goal is to clarify that the abstract problem \eqref{P-intro}-\eqref{Ealpha} is motivated by concrete problems in mathematical physics.  More precisely, we are going to set it up as an abstract version of the following generalized  $n$-dimensional extensible beam equation with nonlinear source and nonlocal damping terms 
\begin{equation}\label{P1}
u_{tt}-\kappa\Delta u +
\Delta^2 u + f(u) + k\big(\|(-\Delta)^{\alpha} u\|^2+\|u_t\|^2\big)u_t = h_\lambda  \ \mbox{ in } \ \Omega\times (0,\infty),
\end{equation}
where $\Omega\subset \mathbb{R}^n$ is a bounded domain with smooth boundary $\partial \Omega$,  $\kappa$ is a non-negative constant,  $f(u)$ is a nonlinear source whose assumption will be given in Section \ref{sec-syst-din} (see Assumption \ref{assumption}), $\alpha\in[0,1]$, $h_\lambda=\lambda h$ with $\lambda\in[0,1]$ and $h\in L^2(\Omega)$,  the notation $\|\cdot \|$ stands for the usual norm in $L^2(\Omega)$, and the   scalar function $k(\cdot)$ is given in some class of functions encompassed by  \ref{k1}-\ref{k3}.
Additionally, problem \eqref{P1} is  considered with initial data
\begin{equation}\label{data}
\begin{array}{l}
u(x,0)=u_0(x),\quad u_t(x,0)=u_1(x),\quad  x\in\Omega,
\end{array}
\end{equation}
and either physical boundary conditions: clamped 
\begin{equation}\label{boundary_clamped} 
u=\frac{\partial u}{\partial \nu}=0\quad \mbox{on}\quad \partial\Omega, 
\end{equation}
where $\nu$ is the outward normal to $\partial\Omega$, or hinged (simply supported)
\begin{equation}\label{boundary_hinged}
u=\Delta u =0\quad \mbox{on}\quad \partial\Omega.
\end{equation}

\subsection{Physical motivation}
In Balakrishnan \cite{Balakrishnan}  it is presented the damping phenomena in flight structures with a free response.  Accordingly, the one-dimensional model is described in terms of the basic second-order dynamics for the displacement variable $y(t)$  as follows
\begin{equation}\label{problem_orig}\ddot{y}(t)+\omega^2y(t)+\gamma D(y(t),\dot{y}(t))=0,
\end{equation}
where $\omega$ is the mode frequency and $\gamma>0$ corresponds to a (small) damping coefficient. Here,  $\dot{y}=\frac{dy}{dt}$ stands for the time derivative. Since  $D(y(t),\dot{y}(t))$  is responsible by damping effects, it has been firstly considered as function 
depending on $\dot{y}(t)$ only, namely, 
$$\mbox{sign}\dot{y}(t),\quad |\dot{y}(t)|\dot{y}(t),\quad |\dot{y}(t)|^{\alpha}\dot{y}(t),\quad \alpha\in(0,1).$$

After that, employing  Krylov-Bogoliubov's approximation, Balakrishnan and Taylor \cite{Balakrishnan-Taylor} 
 suggested a new class of damping models, called by {\it energy damping,} that are based on the instantaneous total energy of the system. More specifically,
 denoting the instantaneous energy associated with the system \eqref{problem_orig} by
$$
\mathcal{E}(t)=\frac{ w^2 }{2}[y(t)]^2+\frac{1}{2}[\dot{y}(t)]^2,
$$ 
then $D(y(t),\dot{y}(t))$ is represented as
\begin{equation}\label{damping-D}
D(y(t),\dot{y}(t))=[\mathcal{E}(t)]^q
\dot{y}(t), \quad q>0.
\end{equation}
Thus, a straightforward computation with \eqref{problem_orig}-\eqref{damping-D} reveals 
$$
\frac{d}{dt}\mathcal{E}(t) \, =\,  -\gamma [\mathcal{E}(t)]^q  [\dot{y}(t)]^2,
$$
from where one sees that the stability of the system is driven by a possibly degenerate nonlocal nonlinear dissipative term involving the energy as a damping coefficient. 

The same conclusion can be done with  Krasovskii's system presented in \cite{krasovskii}, namely,
\begin{equation}\label{krassov}
\ddot{y}(t)+y(t)+ 
k\left([y(t)]^{2}+[\dot{y(t)}]^{2}\right) \dot{y}=0,
\end{equation}
where the damping coefficient function  $k(\cdot)$ is assumed to satisfy suitable properties for the generation of a dissipative dynamical system in $\mathbb{R}^{2}.$ See also \cite[Subsection 5.3.3]{chueshov-solo} for more details on the  infinite-dimensional version of \eqref{krassov}.

By following the same spirit of \cite[Section 4]{Balakrishnan-Taylor}, where  the authors derive some prototypes of 
models for  uniform Bernoulli beams, we consider the following one-dimensional beam bending for a beam of length $2L$ and nonlocal damping coefficient in terms of the energy
\begin{eqnarray}\label{problem_orig_2}
u_{tt}-2\zeta\sqrt{\lambda} u_{xx}+ \lambda u_{xxxx} -\gamma\left[\int_{-L}^{L}\big(\lambda |u_{xx}|^2 + |u_t|^2\big)dx\right]^{q}   u_{xxt}= 0,\end{eqnarray}
where $u=u(x,t)$ represents the transversal deflection of the beam, $\gamma>0$ is a  damping coefficient, $\zeta$ is the constant coming from the Krylov-Bogoliubov approximation  and $\lambda=\frac{2\zeta w}{\sigma^2}$ with $w$ being the mode frequency   and $\sigma^2$ the spectral density of a Gaussian external force. Moreover, for materials whose viscosity can be  essentially seen as friction between moving solids,  one may  interpret beam modes like \eqref{problem_orig_2}  with nonlocal frictional damping term instead
of the viscous one. In this way,
and also in accordance with the energy damping \eqref{damping-D}, the following equation emerges  
\begin{eqnarray}\label{problem_orig_3}
u_{tt}-\kappa u_{xx}+ u_{xxxx} +\gamma\left[\int_{-L}^{L}\big(|u_{xx}|^2 + |u_t|^2\big)dx\right]^{q}   u_{t}= 0,\end{eqnarray}
where we have normalized the equation with respect to an appropriate structural constant $(\lambda=1)$ and also denoted $\kappa=2\zeta$ for the sake of notation. 

Furthermore, in order to see \eqref{problem_orig_3} in the $n$-dimensional scenario, we consider the representative mathematical model
\begin{equation}\label{ncontext}
u_{tt}-\kappa \Delta u + \Delta^2 u+ \gamma \left[\int_{\Omega}\left(|\Delta u|^2+|u_t|^2\right)dx\right]^q u_t=0,
\end{equation}
where  $\Omega$ is a bounded domain of the Euclidian space $\mathbb{R}^n, \, n\geq1$. 

Now, we note that  \eqref{ncontext} represents \eqref{P1} with 
\begin{equation*}
\begin{array}{l}
f=h_\lambda=0,  \ \ k(s)=\gamma s^q, \ \ \alpha=1,  \smallskip    \\ \displaystyle\mathcal{E}_{1}(u,u_t)=\int_{\Omega}\left(|\Delta u|^2+|u_t|^2\right)dx:= \|\Delta u\|^2+\|u_t\|^2,
\end{array}
\end{equation*}
and, roughly speaking,  problem \eqref{P1} is  a particular case of \eqref{P-intro}-\eqref{Ealpha} with 
\begin{equation}\label{partic} 
\begin{array}{l}
A=-\Delta \quad \mbox{and} \quad A_1=\Delta^2.   
\end{array}
\end{equation} 
We still remark  \eqref{P-intro}
can be seen as an
infinite-dimensional generalization of Balakrishnan-Taylor's model \eqref{problem_orig}-\eqref{damping-D} and  Krasovskii's system \eqref{krassov}.

Finally, since the long-time dynamics of problem \eqref{P1}-\eqref{boundary_hinged} seems to be not studied in the literature so far,  we feel motivated to investigate this issue by means of the abstract version \eqref{P-intro} whose precise details on its abstract configuration are presented thereupon.

\subsection{Abstract formulation of problem}\label{sec-wp}

Below we provide the precise details to set problem 
\eqref{P1}-\eqref{boundary_hinged} up  in an abstract formulation as given in \eqref{P-intro}, not only formally taking operators as in \eqref{partic}.
 We advance 
that the theory on functional analysis used below can be found e.g. in  \cite{Brezis,chueshov-solo,chueshov-yellow,lions,lions-magenes,temam1,yosida}.

Throughout this work, the notation $(\cdot , \cdot)$
stands for the $L^2$-inner product,  and $\|\cdot\|_p$  denotes the $L^p$-norm, $p\geq1$.
For commodity, when $p=2$, we design $\|\cdot\|_2=\|\cdot\|$.  As usual, we denote by $H^s(\Omega)$ the $L^2$-based Sobolev space of the order $s>0$ and by $H^s_0(\Omega)$ the closure of $C^{\infty}_0(\Omega)$ in $H^s(\Omega)$. Let us also set $H=L^2(\Omega).$

 We first introduce  the second-order (Laplacian) unbounded  linear operator $A$ by the formula
$$Au=-\Delta u,\quad u\in \mathcal{D}(A)=H^2(\Omega)\cap H^1_0(\Omega).
$$
It is well-known that $A$ is a  positively self-adjoint operator in $H$ and, consequently, we can  define the fractional powers $A^s$ of $A$, $s\in \mathbb{R}$, with domains  $\mathcal{D}(A^{s})$. The spaces $\mathcal{D}(A^{s})$ are Hilbert spaces with inner product  and  norm 
$$(u,v)_{\mathcal{D}(A^{s})}=(A^{s}u,A^{s}v),\quad \|u\|_{\mathcal{D}(A^{s})}=\|A^{s}u\|.$$

Now,  we consider the fourth-order (Biharmonic) unbounded  linear operator  $A_1$ as follows
$$
A_1u=\Delta^2 u,\quad u\in \mathcal{D}(A_1)=\left\{\begin{array}{ll}
H^4(\Omega)\cap H^2_0(\Omega) & \mbox{for } \ \eqref{boundary_clamped} , \smallskip \\
\{u\in H^4(\Omega); \ u=\Delta u =0 \;\;\mbox{on}\;\;\partial\Omega \}& \mbox{for } \ \eqref{boundary_hinged}. 
\end{array}\right.
$$
From the spectral theory, we know that there exists a complete orthonormal basis   $\{w_j\}_{j\in \mathbb{N}}$ in $H$ consisting of eigenvectors of $A_1$, say  with eigenvalues  $\{ \sigma_j\}_{j\in \mathbb{N}}$, such that 
$$
\begin{array}{l}
w_j\in \mathcal{D}(A_1), \quad  A_1w_j=\sigma_jw_j,\ \ j\geq1, \smallskip \\
0<\sigma_1\le \sigma_2\le \cdots \quad \mbox{with} \quad \sigma_j\rightarrow+\infty \ \ \mbox{as} \ \ j\rightarrow +\infty. 
\end{array}
$$
Likewise, we can also define the fractional powers $A^s_1$ of $A_1$,  $s\in \mathbb{R}$, with domains $\mathcal{D}(A^s_1)$  which are  
Hilbert spaces with inner product and  norm
$$(u,v)_{\mathcal{D}(A_1^{s})}=\left(A^s_1 u,A^s_1v\right),\quad \|u\|_{\mathcal{D}(A^{s}_1)}=\|A^{s}_1u\|.$$

In both cases, we have densely  inclusions
$$
\mathcal{D}(A^{s_1})\subset \mathcal{D}(A^{s_2}) \quad \mbox{ and } \quad \mathcal{D}(A^{s_1}_1)\subset \mathcal{D}(A^{s_2}_1), \quad    s_1,s_2\in\mathbb{R}, \, s_1\ge s_2, 
$$
with   continuous embedding when $s_1\ge s_2$ and compact embedding in case $s_1> s_2$.

Let us still stress some particular properties of the above operators and spaces.  In fact, we first note that
$\mathcal{D}(A_1^0)=\mathcal{D}(A^0)=H$, and
$$
\mathcal{D}(A^{1/2}_1)=\left\{\begin{array}{ll}
H^2_0(\Omega) & \mbox{for } \ \eqref{boundary_clamped}, \medskip \\
H^2(\Omega)\cap H^1_0(\Omega) & \mbox{for } \ \eqref{boundary_hinged}.
\end{array}\right.
$$
Hence, in any case, concerning the boundary conditions we have
\begin{equation}\label{equal-norm}
\mathcal{D}(A^{1/2}_1)\subseteq \mathcal{D}(A) \quad \text{ and } \quad \|A^{1/2}_1u\|= \|Au\|, \quad  \forall \,  u\in \mathcal{D}(A^{1/2}_1).
\end{equation}
More particularly, in the specific  case of hinged boundary condition (\ref{boundary_hinged}), one knows that 
\begin{equation}\label{commutative-bc}
A^s u=A^{s/2}_1u, \quad  \forall \,  u\in \mathcal{D}(A^s)=\mathcal{D}(A^{s/2}_1),
\end{equation}
providing a good symmetry  (say commutativity) to the extensible beam model \eqref{P1}, and consequently to problem \eqref{P-intro}. Moreover, the next particular embedding inequalities will be useful throughout this text  
\begin{equation*}
\sigma_1 \|u\|^2 \leq \|A^{1/2}_1 u\|^2, \quad \sigma_1^{1/2}\|A^{1/4}_1 u\|^2 \leq \|A^{1/2}_1 u\|^2, \quad \forall \ u\in \mathcal{D}(A_1^{1/2}),
\end{equation*}
which can be taken for both boundary conditions. 

Under the above construction, we are in the position to rewrite the concrete problem \eqref{P1}-\eqref{boundary_hinged}   in the following 
abstract  class of second-order evolution problems 
\begin{equation*} \label{P-intro-formulation}
\left\{\begin{array}{l}\displaystyle{
	u_{tt}+\kappa A u +
	A_1 u+f(u)+k\big(\|A^{\alpha} u\|^2+\|u_t\|^2\big)u_t= h_\lambda ,\quad t> 0,} \medskip \\
\displaystyle{u(0)=u_0, \quad u_t(0)=u_1,}
\end{array}\right.
\end{equation*} 
which fairly corresponds  to the beginning problem \eqref{P-intro}-\eqref{Ealpha}, where the damping coefficient is assumed to satisfy \ref{k1}-\ref{k3}.

\section{Well-posedness and  associated dynamical system}\label{sec-syst-din}

The asymptotic behavior of the solutions of  problem  \eqref{P-intro}  will be considered on the   Hilbert phase space
\begin{equation*}
\mathcal{H}=\mathcal{D}(A^{1/2}_1)\times H, \quad  ||(u,v)||^{2}_{\mathcal{H}}=\|A^{1/2}_1 u\|^2+\|v\|^{2}, \quad (u,v)\in
\mathcal{H}.
\end{equation*}
As we shall see later, $\mathcal{H}$ is natural finite energy space for \eqref{P-intro}, and if we define  the  Sobolev phase space 
$$
\mathcal{H}_{\alpha}=\mathcal{D}(A^{\alpha})\times H, \quad ||(u,v)||^2_{\mathcal{H}_{\alpha}}=\|A^{\alpha}u\|^2+\|v\|^2, \ \ \alpha \in [0,1],
$$
then
$$
\mathcal{H}=\mathcal{D}(A^{1/2}_1)\times H \subseteq \mathcal{D}(A)\times H= \mathcal{H}^1, \quad ||(\cdot,\cdot)||_{\mathcal{H}}=||(\cdot,\cdot)||_{\mathcal{H}_1}.
$$

The  energy functional $E(u(t),u_t(t)):=E(t),\, t\geq0,$ corresponding to problem \eqref{P-intro} is expressed by
\begin{equation}\label{energy}
E(t)=\frac{1}{2}\left[ \,\|u_t(t)\|^2 +  \|A^{1/2}_1 u(t)\|^2 +
\kappa\|A^{1/2} u(t)\|^2\,\right]
+ \left(\widehat{f}(u(t)),1\right)-\left(h_\lambda,u\right),
\end{equation}
where we set hereafter $\widehat{f}(u)=\int_0^{u}f(\tau)d\tau$ as the primitive of the function $f$, whose assumptions  are given below in order to address the Hadamard well-posedness of problem \eqref{P-intro} as well as its long-time dynamics results.

\begin{Assumption}\label{assumption}
	Let  $f:\mathbb{R}\to\mathbb{R}$ be a  $C^1$-function with $f(0)=0$  
	and satisfying
	\begin{eqnarray}\label{assumption_f'}
	|f'(u)| \,   \le \,   C_{f'}(1+|u|^{p}),\quad   u\in \mathbb{R},
	\end{eqnarray}
	\begin{equation} \label{assumption_f}
	-C_f-\frac{c_f}{2} |u|^2 \leq  \widehat{f}(u)\leq  f(u)u +\frac{c_f}{2} |u|^2, \quad   u \in \mathbb{R},
	\end{equation}
	for some  constants $C_f,C_{f'}>0$,    $c_f\in[0,\sigma_1)$  and growth exponent $p
	\leq\frac{4}{n-4}$ for $n\geq5$.
\end{Assumption}

\begin{remark}\label{rem-criti-growth}
	We observe that  $p^*:=2(p+1)\leq\frac{2n}{n-4}$
	can be the critical (Sobolev) exponent for the continuous embedding $\mathcal{D}(A^{1/2}_1)\hookrightarrow L^{p^*}(\Omega),$ which is not compact for the critical case
	$p=\frac{4}{n-4} $. From Section 
	$\ref{sec-tech-results}$ to  Section $\ref{sec-dynamics-k2}$ we present several results on long-time behavior in critical and subcritical frameworks.
	We do not work in lower dimensions 
	$1\leq n\leq 4$ since they do not require compactness issues w.r.t. $p$, by holding the same results with unchanged computations. 
	
\end{remark}

\subsection{Well-posedness}

\begin{theorem}[Hadamard Well-Posedness]\label{theo-w-p-general}
Let us assume that $\alpha,\lambda\in[0,1],$   $h_\lambda:=\lambda h\in H,$  and $f$ satisfies {Assumption $\ref{assumption}$}. Then, problem \eqref{P-intro} is Hadamard well-posed for $k(\cdot)$ given in any case  {\rm \ref{k1}}-{\rm \ref{k3}}.
\end{theorem}

The proof of Theorem \ref{theo-w-p-general} is  given below separately concerning  the cases \ref{k1}-\ref{k3}.
We start with the possibly degenerate case  \ref{k1} in what concerns the function $k(\cdot)$. The other cases \ref{k2}-\ref{k3} shall be treated at the end of this section.

\smallskip 

\noindent {\bf Well-posedness: case \ref{k1}.} In this case,    problem \eqref{P-intro} can be  written explicitly as
\begin{equation} \label{P}
\left\{\begin{array}{l}\displaystyle{
	u_{tt}+\kappa A u +
	A_1 u+\gamma\left[\,\|A^{\alpha} u\|^2+\|u_t\|^2\,\right]^qu_t+f(u)= h_\lambda ,\quad t> 0,} \medskip \\
\displaystyle{u(0)=u_0, \quad u_t(0)=u_1,}
\end{array}\right.
\end{equation}

\begin{theorem}[Existence and Uniqueness] \label{theo-existence} Let    $\gamma>0$, $ q \geq \frac{1}{2}$ and  $\kappa\geq0$   be given constants. Additionally, let us take  on $\alpha,\lambda\in[0,1],$  $h_\lambda=\lambda h\in H,$  and {Assumption $\ref{assumption}$}. 
\begin{itemize}
\item[$(i)$]  If $(u_0,u_1)\in \mathcal{H}$, then there exists
		$T_{\max}>0$  such that problem  \eqref{P} has a unique mild (weak) 
		solution  $(u^\lambda,u_t^\lambda):=(u,u_t)$ in the class
\begin{equation*}
(u,u_t)\in C([0,T_{\max}),\mathcal{H}).
\end{equation*}		
		
\item[$(ii)$]  If  $(u_0,u_1)\in \mathcal{D}(A_1)\times\mathcal{D}(A^{1/2}_1)$, then the   solution $U$ is more regular (strong). 		
\end{itemize}
In both cases, we have    $T_{\max}=+\infty$.
\end{theorem}
\begin{proof}

 We first define   vector-valued function $
 U(t):=(u(t),v(t)), \, t\geq0,$ with 
 $v=u_{t}.$   Then we can rewrite  system  $(\ref{P})$  as the following first order abstract problem
 \begin{eqnarray}\label{abstrac-cauchy}
 \left\{\begin{array}{l}
 U_t =\mathcal{A}U +\mathcal{M}(U), \quad t>0, \medskip \\
 U(0)=(u_0,u_1):=U_0,
 \end{array}\right.
 \end{eqnarray}
where   $\mathcal{A}:\mathcal{D}(\mathcal{A})\subset\mathcal{H}\to\mathcal{H}$ is the linear  operator given by
\begin{equation}\label{def_A}
\mathcal{A}U = ( v , -A_1 u),
\quad    U \in \mathcal{D}(\mathcal{A}):=\mathcal{D}(A_1)\times \mathcal{D}(A_1^{1/2}),
\end{equation}
and $\mathcal{M}:\mathcal{H}\to\mathcal{H}$ is the nonlinear operator
\begin{equation}\label{def_B}
\mathcal{M}(U)= (0, -\kappa A u+\gamma||U||^{2q}_{\mathcal{H}^{\alpha}}v  -f(u)+h_\lambda), \quad U \in\mathcal{H}.
\end{equation}
Therefore, the existence and uniqueness of solution to the system $(\ref{P})$ rely on the study of 
problem (\ref{abstrac-cauchy}). To this purpose, and according to   Pazy   \cite[Chapter 6]{Pazy}, it is enough to prove that $\mathcal{A}$  given in  (\ref{def_A})  is the infinitesimal generator of a $C_{0}$-semigroup of contractions $e^{\mathcal{A}t}$ (which is very standard) and $\mathcal{M}$ set in  (\ref{def_B}) is locally Lipschitz   on $\mathcal{H}$ (which will be done in Appendix 
\ref{subsec-A1} for the sake of completeness -- here is the precise moment where, for a  technical reason, we consider $q\geq1/2$ as clarified later). From this, one concludes the proof of items $(i)$ and $(ii)$.

It remains to   check that $T_{\max}=+\infty,$ that is, both mild and regular solutions are globally defined in time.    Indeed,
taking the scalar product in $H$ of (\ref{P}) with $u_t$, we obtain
\begin{eqnarray}\label{energy_relation_AA}
\frac{d}{dt}E(t)+\gamma ||(u(t),u_t(t))||_{\mathcal{H}^{\alpha}}^{2q}\,\|
u_t(t)\|^{2}=0, \quad t>0.
\end{eqnarray}
Integrating (\ref{energy_relation_AA}) over $(0,t), \, t>0$, we get
\begin{eqnarray}\label{energy_relation}
E(t)
+\gamma\int_0^t||(u(\tau),u_t(\tau))||_{\mathcal{H}^{\alpha}}^{2q}\,\|
u_t(\tau)\|^{2} \,d\tau = E(0), \quad  t>0,
\end{eqnarray}

Now, using  (\ref{assumption_f}) and Young's inequality with    $\sigma_1\omega:=\sigma_1-{c_f}>0$, we have
\begin{equation*}
\left[\,\left(\widehat{f}(u(t)),1\right)-\left(h_\lambda,u(t)\right)\,\right]\ge -\frac{c_f}{2\sigma_1}\|A^{1/2}_1 u(t)\|^2-C_f|\Omega|-\frac{1}{\omega\sigma_1}\|h_\lambda\|^2-\frac{\omega}{4}\|A^{1/2}_1 u(t)\|^2.
\end{equation*}
From the energy defined (\ref{energy}), we infer
\begin{align*}
E(t)  & \geq \frac{\omega}{4}\|A^{1/2}_1 u(t)\|^2+ \frac{1}{2}\|u_t(t)\|^2 -C_f|\Omega|-\frac{1}{\omega\sigma_1}\|h\|^2 \\ & \ge \frac{\omega}{4}||(u(t),u_t(t))||^2_{\mathcal{H}}-K_\lambda, 
\end{align*}
where we denote $K_\lambda=\left[C_f|\Omega|+\frac{1}{\omega\sigma_1}\|h_\lambda\|^2\right]>0$. From this  and \eqref{energy_relation}, we arrive at
\begin{align}\label{global-in-time}
\frac{\omega}{4}\|(u(t),u_t(t))\|_{\mathcal{H}}^2 \leq   {E}(t)+K_\lambda \le {E}(0)+K_\lambda , \quad \forall \ t\in[0,T_{\max}).
\end{align}
If $T_{\max}<+\infty,$ we have that $\|(u(t),u_t(t))\|_{\mathcal{H}}$ blows up in finite time (cf. \cite[Theorem $1.4$]{Pazy}), which is contraction with  (\ref{global-in-time}). Therefore, $T_{\max}=+\infty.$
\end{proof}

Hereafter, for the sake of notation, we still omit the parameter $\lambda\in[0,1]$ indexed to the solution of \eqref{P}, by  simply writing down $u^\lambda:=u$.

\smallskip

To the continuous dependence result, we invoke the following technical lemma on the power function $s^r$ for $r\geq1$ (cf. \cite{Aloiui-Ben-Haraux}).

\begin{lemma}[{\cite[Lemma 2.2]{Aloiui-Ben-Haraux}}]
	\label{Lemma-Haraux}
Let $X$ be a normed space with norm $\|\cdot\|_X$. Then, for any $r\geq 1$, we have
\begin{equation}\label{Haraux-variant}
\left|\|u\|_X^r-\|v\|_X^r\right| \leq r\max\{\|u\|_X,\|v\|_X\}^{r-1}\|u-v\|_X, \quad \forall \,\, u,v \in X.
\end{equation} 
\end{lemma}

\begin{theorem}[Continuous Dependence]\label{theo-dependence}
Under the assumptions of Theorem  $\ref{theo-existence}, $ let 
$U^1=(u^1,u^1_{t})$ and $ U^2=(u^2,u^2_{t})$
be two (strong or weak) solutions of problem $(\ref{P})$
corresponding to initial data $U^1_0=(u^1_{0},u^1_{1})$ and $U^2_0=(u^2_{0},u^2_{1})$, respectively, and let $T>0$ be any positive time. Then,  there exists
a positive non-decreasing function
$\mathcal{Q}(t)=\mathcal{Q}\big(||U^1_0||_{\mathcal{H}},||U^2_0||_{\mathcal{H}},t\big)$ such that
\begin{equation}\label{dependence}
||U^1(t)-U^2(t)||_{\mathcal{H}}  \, \le  \, \mathcal{Q}(t)
||U^1_0-U^2_0||_{\mathcal{H}},\quad t\in[0,T].
\end{equation}
\end{theorem}
\begin{proof}
Setting $w=u^1-u^2$ and $F(w)=f(u^1)-f(u^2)$,
the difference $U^1-U^2=(w,w_{t})$ is a solution (in
the strong and weak sense) of the following problem 
\begin{equation}\label{dependenceI}
\left\{\begin{array}{l}
\displaystyle{w_{tt}+\kappa A w+A_1 w+\frac{\gamma}{2}\Pi_1 w_t+\frac{\gamma}{2}\Pi_2\left[\,u^1_t+u^2_t\,\right]+F(w)=0,} \smallskip \\
\displaystyle{w(0)= u^1_0-u^2_0:=w_0,\quad w_t(0)=u^1_1-u^2_1:=w_1,}
\end{array}\right.
\end{equation}
where
$$\Pi_i(t)=||U^1(t)||^{2q}_{\mathcal{H}^{\alpha}}+(-1)^{1-i}||U^2(t)||^{2q}_{\mathcal{H}^{\alpha}},\quad i=1,2.$$
The next estimates are done for strong solutions, and the same result  can be achieved for weak solution through density arguments.

We still denote
\begin{equation}\label{def_V}
\mathcal{E}_w(t)=\|w_t(t)\|^2+\kappa\|A^{1/2} w(t)\|^2+\|A^{1/2}_1 w(t)\|^2, \quad t\geq0.
\end{equation}
Taking the inner product in $H$ of (\ref{dependenceI}) with $w_t$, we obtain
\begin{eqnarray}\label{dependenceII}
\frac{1}{2}\frac{d}{dt}\mathcal{E}_w(t)+\frac{\gamma}{2}\Pi_1(t)\|w_t(t)\|^2\,=\,\mathcal{J}_1+\mathcal{J}_2,
\end{eqnarray}
where
\begin{eqnarray*}
\mathcal{J}_1&=&-\,\left(\,F(w),w_t\,\right),\\
\mathcal{J}_2&=&-\,\frac{\gamma}{2}\Pi_2(t)\left(\,(u^1_t+
u^2_t),w_t\,\right).
\end{eqnarray*}
Using that
$\mathcal{D}(A^{1/2}_1)\hookrightarrow \mathcal{D}(A)\hookrightarrow \mathcal{D}(A^{1/2})$, it is easy
to see that
\begin{equation}\label{dependenceIII}
||U^1(t)-U^2(t)||_{\mathcal{H}}^2\le \mathcal{E}_w(t)\le (1+\kappa
\mu_0)||U^1(t)-U^2(t)||_{\mathcal{H}}^2.
\end{equation}
where $\mu_0>$ is a constant independent of  initial data. In what follows, we will denote by $C$ several constants that depend on the initial data.
From assumption (\ref{assumption_f'}), H\"older's inequality with
$\frac{p}{p^*} + \frac{1}{p^*} + \frac{1}{2} = 1$, the embedding
$\mathcal{D}(A^{1/2}_1)\hookrightarrow L^{p^*}(\Omega)$ and
(\ref{dependenceIII}) we can estimate the term $\mathcal{J}_1$ as
follows
\begin{eqnarray}\label{J1}
\left| \mathcal{J}_1\right|&\le&C_{f'}\left[\,|\Omega|+\|u^1(t)\|^{p^*}_{p^*}+\|u^2(t)\|^{p^*}_{p^*}\,\right]^{\frac{p}{p^*}}\|w(t)\|_{p^*}\|w_t(t)\| \nonumber\\
&\le& C\|w(t)\|_{p^*}\|w_t(t)\|\le
C\|A^{1/2}_1w(t)\|\|w_t(t)\|\le C \mathcal{E}_w(t),
\end{eqnarray}
for some constant $C>0$ depending on initial data. This is the precise moment where we should use $q\geq1/2$ to the handling of $\Pi_2(t)$. Indeed, since  $2q\geq 1$, we can apply Lemma \ref{Lemma-Haraux} to obtain the following estimate 
$$\left|\Pi_2(t)\right| \leq 2q \max\{\|U^1(t)\|_{\mathcal{H}^{\alpha}},\|U^2(t)\|_{\mathcal{H}^{\alpha}}\}^{2q-1}\left[\mathcal{E}_w(t)\,\right]^{1/2},$$
and then,
\begin{eqnarray}\label{J2}
\left|\mathcal{J}_2\right|
 &\le& \frac{\gamma}{2}\left|\Pi_2(t)\right|\left[\,\|
u^1_t(t)\|+ \|u^2_t(t)\|\,\right]\|w_t(t)\|\\
&\le& C\mathcal{E}_w(t),\nonumber
\end{eqnarray}
for some constant $C>0$ depending on initial data.
Replacing the above estimates \eqref{J1}-\eqref{J2} in
(\ref{dependenceII}), we obtain
\begin{eqnarray}\label{dependenceV}
\frac{1}{2}\frac{d}{dt}\mathcal{E}_w(t)+\frac{\gamma}{2}\Pi_1(t)||w_t(t)||^2\le
C\mathcal{E}_w(t),
\end{eqnarray}
for all $t\in [0,T]$ and  some constant $C>0$ depending on initial data.  Hence,  integrating (\ref{dependenceV}) on
$[0,t]$, applying Gronwall's inequality and using (\ref{dependenceIII}), we arrive at
\begin{eqnarray*}
||U^1(t)-U^2(t)||^2_{\mathcal{H}}  \leq C_0e^{C_1 t  } ||U^1_0-U^2_0||^2_{\mathcal{H}},
\end{eqnarray*}
for some  constants $C_i=C_i(||(u^i_0,u^i_1)||_{\mathcal{H}})>0, \, i=0,1,$ which proves (\ref{dependence}) with non-decreasing function 
$\mathcal{Q}(t):=C_0e^{C_1 t}$ as desired. 
\end{proof}

\noindent {\bf Well-posedness: cases \ref{k2}-\ref{k3}.} Concerning the well-posedness of problem \eqref{P-intro} with function $k(\cdot)$ given in cases  {\rm \ref{k2}}  or  {\rm\ref{k3}}, we can proceed verbatim as in the proofs of Theorems  $ \ref{theo-existence} $ and $\ref{theo-dependence} $ by noting that such cases do not interfere in the locally Lipschitz property of operator  
	 $\mathcal{M}$ set by  \eqref{def_B} nor in the computations for the continuous dependence  \eqref{dependence}. Thus, we shall omit the details here.\qed
	  
 For similar approaches in these cases, we refer to 	\cite[Sect. 2]{jorge-narciso-DCDS} and \cite[Subsect. 5.3.3]{chueshov-solo}.

\subsection{Generation of the dynamical system}\label{subsec-generation-syst-din}

For every $h_\lambda =\lambda h \in H, \, \lambda\in[0,1],$ Theorem 
\ref{theo-w-p-general}
 ensures the Hadamard well-posedness of problem \eqref{P-intro} and, consequently,  the  definition  of a family of   nonlinear $C_0$-semigroups  
$S_\lambda(t):\mathcal{H}\rightarrow \mathcal{H}$ given by
\begin{equation}\label{sist-din}
S_\lambda(t)(u_{0},u_{1})=(u^\lambda(t),u^\lambda_{t}(t)):=(u(t),u_{t}(t)), \quad t\geq0,
\end{equation}
where  $(u,u_{t})$ is the unique  solution of the abstract system
$(\ref{P-intro})$.  Moreover, through the condition (\ref{dependence}) one sees that $S_\lambda(t)$ is  locally
Lipschitz continuous on the phase space $\mathcal{H}$. 

Accordingly, we also note the change of the time variable 
$t \mapsto-t$ in problem \eqref{P-intro} turns itself into the same problem unless the damping term, which is replaced  by $-k\big(\mathcal{E}_{\alpha}(u,u_t)\big)u_t$ instead of $k\big(\mathcal{E}_{\alpha}(u,u_t)\big)u_t$ in \eqref{P-intro}. For such a
``inverse-time'' problem, we can use the same arguments as in  
Theorems
\ref{theo-existence} and \ref{theo-dependence} for all cases concerning $k(\cdot)$
to prove the well-posedness.
Therefore, it allows us to set $S_\lambda(t)$ as an evolution $C_0$-group. This remark will be important for future discussions on full trajectories and invariance properties of bounded sets in 
$\mathcal{H}$.

Therefore, in what follows the
dynamics of problem \eqref{P-intro} shall be studied through its corresponding 
dynamical system $(\mathcal{H},S_\lambda(t))$ originated by \eqref{sist-din}. We start with the case \ref{k1} and then we treat the other cases \ref{k2} and \ref{k3}  as well.

To deal with computations in the next sections, 
it is worth keeping in mind the
energy   $E(t)$   given in \eqref{energy} and, additionally, we  set the notation  $\omega:=1-\frac{c_f}{\sigma_1}>0$, as well as the following parameters and modified energy 
\begin{equation}\label{modified-energy}
K_\lambda=\left[C_f|\Omega|+\frac{1}{\sigma_1\omega}\|h_\lambda\|^2\right]>0  \quad \mbox{and} \quad  \widetilde{E}(t)=E(t)+K_\lambda, \ \ t\geq0.
\end{equation}

\section{Technical results: case \ref{k1}}\label{sec-tech-results}

   In what follows, except in Subsection \ref{subsec-lip-cont}, we will continue using the simplified notation  $u^\lambda:=u$, by neglecting (whenever there is no confusion) the index $\lambda\in[0,1]$.

\subsection{Dissipativity and gradient property}\label{sec-useful-results}

The first proposition gives us  lower and upper estimates for  $\widetilde{E}(t)$. It reads as follows: 
\begin{proposition}\label{P2} 
	Under the assumptions of Theorem  $\ref{theo-existence} $,   
	there exist constants $C_{\alpha,q,\gamma}>0$ and  $C_{\widetilde{E}(0)}>0$ such that the modified energy $\widetilde{E}(t)$ satisfies
	\begin{eqnarray}\label{inequality_pricipal}
	\left[\frac{q}{C_{\alpha,q,\gamma}}t+\frac{1}{\big[\widetilde{E}(0)\big]^{q}}\,\right]^{-1/q}\le \widetilde{E}(t)\le\left[\frac{q}{C_{\widetilde{E}(0)}}(t-1)^{+}+\frac{1}{\big[\widetilde{E}(0)\big]^{q}}\right]^{-1/q}+8K_\lambda,
	\end{eqnarray}
	for all $t>0$, where we use the standard notation $s^+:=(s+|s|)/2$.
\end{proposition}
\begin{proof}
	
	We initially prove the first inequality of \eqref{inequality_pricipal}.	 Taking the scalar product in $H$ of (\ref{P}) with $u_t$ and using that $\frac{d}{dt}\widetilde{E}(t)=\frac{d}{dt}E(t)$, we obtain
	\begin{eqnarray}\label{absorbing_CC}
	\frac{d}{dt}\widetilde{E}(t)=-\,\gamma||(u(t),u_t(t))||^{2q}_{\mathcal{H}^{\alpha}}\|u_t(t)\|^2, \quad t>0.
	\end{eqnarray}
	Now, using the embedding $\mathcal{D}(A^{1/2}_1)\hookrightarrow \mathcal{D}(A) \hookrightarrow \mathcal{D}(A^{\alpha})$, with constant $C_\alpha>0$ to the second one,  the expression for $\widetilde{E}(t)$ in \eqref{modified-energy}, and also (\ref{global-in-time}),  we get
	\begin{eqnarray*}
		\gamma|| (u(t),u_t(t))||^{2q}_{\mathcal{H}^{\alpha}}\|u_t(t)\|^2&\le&\gamma C_{\alpha}^{2q}|| (u(t),u_t(t))||^{2q}_{\mathcal{H}}\|u_t(t)\|^2 \\
		&=& \gamma \left[\frac{4C_{\alpha}}{\omega}\right]^q\left[\,\widetilde{E}(t)\,\right]^q\left[2\widetilde{E}(t)\right]\\
		&=&\frac{2^{2q+1}C_{\alpha}^{q}\gamma}{\omega^{q}}\left[\,\widetilde{E}(t)\,\right]^{q+1}.
	\end{eqnarray*}
	Returning to (\ref{absorbing_CC}), one sees  that   
	\begin{eqnarray*}
		\frac{d}{dt}\widetilde{E}(t)\left[\,\widetilde{E}(t)\,\right]^{-(q+1)}\ge -\,\frac{2^{2q+1}C_{\alpha}^{q}\gamma}{\omega^{q}},
	\end{eqnarray*} 
	or else, 
	\begin{eqnarray}\label{LLL}
	-\,\frac{1}{q}\frac{d}{dt}\left[\,\widetilde{E}(t)\,\right]^{-q}\ge -\,\frac{2^{2q+1}C_{\alpha}^{q}\gamma}{\omega^{q}}.
	\end{eqnarray}
	Simply solving this differential inequality, we arrive at 
	\begin{equation*} 
	\widetilde{E}(t)\ge  \left[\,\frac{2^{2q+1}C_{\alpha}^{q}\gamma q}{\omega^{q}}t+\frac{1}{\big[\widetilde{E}(0)\big]^{q}}\,\right]^{-1/q},
	\end{equation*}
	which proves the first inequality in (\ref{inequality_pricipal}) with $C_{\alpha,q,\gamma}= \frac{\omega^{q}} {2^{2q+1}C_{\alpha}^{q}\gamma}>0.$

	Now, we are going to prove the second inequality of \eqref{inequality_pricipal}. To do so, we provide some proper estimates and then apply Nakao's method (cf. \cite{Nakao}).
	
	Again from   (\ref{global-in-time}), we note that
	\begin{equation}\label{energy_norm} 
	\frac{\omega}{4}||(u(t),u_t(t)||^2_{\mathcal{H}}\leq\widetilde{E}(t) \leq   {E}(t)+K_1, \quad t\geq0, \ \lambda\in[0,1].
	\end{equation}
	We also observe that
	\begin{eqnarray}\label{absorbing_BB}
	\gamma||(u(t),u_t(t)||^{2q}_{\mathcal{H}^{\alpha}}\|u_t(t)\|^2\ge \gamma\|u_t(t)\|^{2q}\|u_t(t)\|^2= \gamma\|u_t(t)\|^{2(q+1)},
	\end{eqnarray}
	and replacing  (\ref{absorbing_BB}) in (\ref{absorbing_CC}), we get
	\begin{eqnarray}\label{absorbing_C'}
	\frac{d}{dt}\widetilde{E}(t)+\gamma\|u_t(t)\|^{2(q+1)}\le 0, \quad t>0,
	\end{eqnarray}
	which implies that $\widetilde{E}(t)$ is non-increasing. Also, integrating (\ref{absorbing_C'}) from $t$ to $t+1$, we obtain
	\begin{eqnarray}\label{C''}
	\gamma\int_t^{t+1}\|u_t(s)\|^{2(q+1)}\,ds\le \widetilde{E}(t)-\widetilde{E}(t+1):= [\,D(t)\,]^2.
	\end{eqnarray}
	Using H\"older's inequality with $\frac{q}{q+1}+\frac{1}{q+1}=1$ and (\ref{C''}), we infer
	\begin{eqnarray}
	\int_t^{t+1}\|u_t(s)\|^2ds\le\left[\int_t^{t+1}1^{\frac{q+1}{q}}ds\right]^{\frac{q}{q+1}} \left[\int_t^{t+1}\|u_t(s)\|^{2(q+1)}ds\right]^{\frac{1}{q+1}}\le \frac{1}{\gamma^{\frac{1}{q+1}}}[\,D(t)\,]^{\frac{2}{q+1}}.\label{G}
	\end{eqnarray}
	From the Mean Value Theorem for integrals, there exist $t_1\in [t,t+\frac{1}{4}]$ and $t_2\in [t+\frac{3}{4},t+1]$ such that
	\begin{eqnarray}\label{eq37}
	\|u_t(t_i)\|^2\le 4\int_t^{t+1}\|u_t(s)\|^2ds\le \frac{4}{\gamma^{\frac{1}{q+1}}}[\,D(t)\,]^{\frac{2}{q+1}}.
	\end{eqnarray}
	
	On the other hand, taking 
	the scalar product in $H$ of (\ref{P}) with $u$ and integrating over $[t_1,t_2]$, we have
	\begin{eqnarray}\label{D}
	&&\int_{t_1}^{t_2}\left[\,\|A^{1/2}_1 u(s)\|^2+\kappa\|A^{1/2} u(s)\|^2+\left(\,f(u(s)),u(s)\right)-\left(h_\lambda,u(s)\right)\,\right]ds\nonumber\\
	&&\qquad\qquad\qquad =\int_{t_1}^{t_2}\|u_t(s)\|^2ds+\sum_{i=1}^2{F}_i,
	\end{eqnarray}
	where
	\begin{eqnarray*}
		&&	{F}_1:=\left(u_t(t_1),u(t_1)\right)-\left(u_t(t_2),u(t_2)\right) \\
		&&	{F}_2:=-\,\gamma\int_{t_1}^{t_2}||(u(s),u_t(s))||^{2q}_{\mathcal{H}^{\alpha}}\left( u_t(s),u(s)\right)ds.
	\end{eqnarray*}
	From assumption (\ref{assumption_f}), we have
	\begin{eqnarray*}
		\int_{t_1}^{t_2}\left(f(u(s)),u(s)\right)\,ds\ge \int_{t_1}^{t_2}\left(\widehat{f}(u(s)),1\right)ds-\frac{c_f}{2\sigma_1}\int_{t_1}^{t_2}\|A^{1/2}_1 u(s)\|^2ds.
	\end{eqnarray*}
	Returning to (\ref{D}) and using that $\omega=1-\frac{c_f}{\sigma_1}$, we get
	\begin{eqnarray}\label{F}
	&&\int_{t_1}^{t_2}\widetilde{E}(s)\,ds+\frac{1}{2}\int_{t_1}^{t_2}\left[\,\kappa\|A^{1/2} u(t)\|^2+\omega\|A^{1/2}_1 u(t)\|^2\,\right]ds\nonumber\\
	&&\qquad\qquad\qquad \le K_\lambda+\frac{3}{2}\int_{t_1}^{t_2}\|u_t(s)\|^2ds+\sum_{i=1}^2{F}_i.
	\end{eqnarray}
	Let us estimate the terms $F_1$ and $F_2$ as follows. Firstly, we note that through H\"older's inequality, \eqref{energy_norm}, 
	(\ref{eq37}) and Young's inequality, we obtain  
	\begin{eqnarray*}
		{F}_1&\le &\sum_{i=1}^2\|u_t(t_i)\|\|u(t_i)\|\le \frac{1}{\sigma_1^{1/2}}\sum_{i=1}^2\|u_t(t_i)\|\|A^{1/2}_1 u(t_i)\|\\
		&\le&
		\frac{8}{(\omega\sigma_1)^{1/2}\gamma^{\frac{1}{2(q+1)}}}[\,D(t)\,]^{\frac{1}{q+1}}
		\sup_{t_1\le s\le t_2}[\widetilde{E}(s)]^{1/2}\\
		&\le & 
		\frac{128}{\omega\sigma_1\gamma^{\frac{1}{q+1}}}
		[\,D(t)\,]^{\frac{2}{q+1}}
		+\frac{1}{8}\sup_{t_1\le s\le t_2}\widetilde{E}(s).
	\end{eqnarray*}
	Additionally, using again that $\mathcal{D}(A^{1/2})\hookrightarrow \mathcal{D}(A^{\alpha/2})$ with embedding constant $C_\alpha$, and (\ref{global-in-time}), we have
	$$
	\| (u(t),u_t(t))\|^{2q}_{\mathcal{H}^{\alpha}}\le C_{\alpha}^{q}\|(u(t),u_t(t))\|^{2q}_{\mathcal{H}}\le \frac{2^qC_{\alpha}^{q}}{\omega^q}\left[\,\widetilde{E}(0)\,\right]^q.
	$$
	From this, using (\ref{G}), and proceeding as in the estimate for $F_1$, we also get
	\begin{eqnarray*}
		{F}_2&\le& \frac{\gamma 2^{q}C_{\alpha}^{q}}{\omega^q\sigma_1^{1/2}}\left[\,\widetilde{E}(0)\,\right]^q\int_{t_1}^{t_2}\|u_t(s)\|\|A^{1/2}_1 u(s)\|\,ds\\
		&\le&
		\frac{2^q\gamma^{\frac{2q+1}{2(q+1)}}C_{\alpha}^{q}}{\omega^q\sigma_1^{1/2}}[\,\widetilde{E}(0)\,]^q[\,D(t)\,]^{\frac{1}{q+1}}
		\sup_{t_1\le s\le t_2}\|A^{1/2}_1 u(s)\|_2\\
		&\le &
		\frac{2^{q+1}\gamma^{\frac{2q+1}{2(q+1)}}C_{\alpha}^{q}}{\omega^{q+1/2}\sigma_1^{1/2}}[\,\widetilde{E}(0)\,]^q[\,D(t)\,]^{\frac{1}{q+1}}
		\sup_{t_1\le s\le t_2}\widetilde{E}^{1/2}(s)\nonumber\\
		&\le &
		\frac{2^{2q+3}\gamma^{\frac{2q+1}{q+1}}C^{2q}_{\alpha}}{\omega^{2q+1}\sigma_1}
		[\,\widetilde{E}(0)\,]^{2q}[\,D(t)\,]^{\frac{2}{q+1}}
		+\frac{1}{8}\sup_{t_1\le s\le t_2}\widetilde{E}(s).
	\end{eqnarray*}
	
	Regarding again (\ref{G}) and replacing the estimates for $F_i, \, i=1,2,$ in  (\ref{F}), we obtain
	\begin{eqnarray}\label{A3}
	\int_{t_1}^{t_2}\widetilde{E}(s)\,ds\le \overline{C}_{\widetilde{E}(0)}[\,D(t)\,]^{\frac{2}{q+1}}+\frac{1}{4}\sup_{t_1\le s\le t_2}\widetilde{E}(s)+K_\lambda,
	\end{eqnarray}
	where we set $$\overline{C}_{\widetilde{E}(0)}:=\left[\,\frac{3}{2\gamma^{\frac{1}{q+1}}}+\frac{128}{\omega\sigma_1\gamma^{\frac{1}{q+1}}}+\frac{2^{2q+3}\gamma^{\frac{2q+1}{q+1}}C^{2q}_{\alpha}}{\omega^{2q+1}\sigma_1}
	[\,\widetilde{E}(0)\,]^{2q}\,\right]>0.$$
	
	Using once more the Mean Value Theorem for integrals and the fact that $\widetilde{E}(t)$ is non-increasing, there exists $\zeta\in[t_1,t_2]$
	such that
	$$
	\int_{t_1}^{t_2}\widetilde{E}(s)\,ds=\widetilde{E}(\zeta)(t_2-t_1)\ge \frac{1}{2}\widetilde{E}(t+1),
	$$
	and then 
	$$
	\sup_{t\le s\le t+1}\widetilde{E}(s)=\widetilde{E}(t)=\widetilde{E}(t+1)+[\,D(t)\,]^2 \leq 2\int_{t_1}^{t_2}\widetilde{E}(s)\,ds+[\,D(t)\,]^2.
	$$
	Thus, from this and (\ref{A3}), we arrive at
	\begin{eqnarray*}
		\sup_{t\le s\le t+1}\widetilde{E}(s)&\le& [\,D(t)\,]^2+2\int_{t_1}^{t_2}\widetilde{E}(s)ds\\
		&\le&
		[\,D(t)\,]^2+2\overline{C}_{\widetilde{E}(0)}[\,D(t)\,]^{\frac{2}{q+1}}+\frac{1}{2} \sup_{t\le s\le t+1}\widetilde{E}(s)+2K_\lambda,
	\end{eqnarray*}
	and since $0<\frac{2}{q+1}\le 2$, we obtain
	\begin{eqnarray}\label{xxxx}
	\sup_{t\le s\le t+1}\widetilde{E}(s)\le
	[\,D(t)\,]^{\frac{2}{q+1}}\left[\,2[\,D(t)\,]^{\frac{2q}{q+1}}+4\overline{C}_{\widetilde{E}(0)}\right]+4K_\lambda.
	\end{eqnarray}
	Observing that 
	\begin{eqnarray*}
		2[\,D(t)\,]^{\frac{2q}{q+1}}\le 2\left[\widetilde{E}(t)+\widetilde{E}(t+1)\right]^{\frac{q}{q+1}}\le2^{\frac{2q+1}{q+1}} \left[ \widetilde{E}(0)\,\right]^{\frac{q}{q+1}},
	\end{eqnarray*}
	and denoting by 
	$$
	{C}_{\widetilde{E}(0)}:= 2^{q+1}\left[2^{\frac{2q+1}{q+1}} \left[ \widetilde{E}(0)\,\right]^{\frac{q}{q+1}}+4\overline{C}_{\widetilde{E}(0)}\right]^{q+1}>0,
	$$
	and also recalling the definition of $[D(t)]^2$ in (\ref{C''}), we obtain from (\ref{xxxx}) that 
	\begin{eqnarray*}
		\sup_{t\le s\le t+1}[\widetilde{E}(s)]^{q+1}\le
		C_{\widetilde{E}(0)}[\,\widetilde{E}(t)-\widetilde{E}(t+1)\,]+[8K_\lambda]^{q+1}.\nonumber
	\end{eqnarray*}
	
	Hence, applying  \cite[Lemma 2.1]{Nakao} with $\widetilde{E}=\phi$, $C_{\widetilde{E}(0)}=C_0,$ and $K=[8K_\lambda]^{q+1},$ we conclude
	\begin{equation*}
	\widetilde{E}(t)\le\left[\frac{q}{C_{\widetilde{E}(0)}}(t-1)^{+}+\frac{1}{\big[\widetilde{E}(0)\big]^{q}}\right]^{-1/q}+8K_\lambda,
	\end{equation*}
	which ends the proof of the second inequality in (\ref{inequality_pricipal}).
	
	The proof of  Proposition \ref{P2} is therefore complete.
\end{proof} 

\begin{remark}
	It is worth point out that we always have
	\begin{equation}\label{ineq}
	\left[\frac{q}{C_{\alpha,q,\gamma}}t+\frac{1}{\big[\widetilde{E}(0)\big]^{q}}\,\right]^{-1/q}  \le \  \left[\frac{q}{C_{\widetilde{E}(0)}}(t-1)^{+}+\frac{1}{\big[\widetilde{E}(0)\big]^{q}}\right]^{-1/q},
	\end{equation}	
	so that 	it makes sense to express $\widetilde{E}(t)$ between the inequalities in \eqref{inequality_pricipal}. Indeed, from the definitions of $C_{\alpha,q,\gamma}$ and $C_{\widetilde{E}(0)}$ in the proof of Proposition $ \ref{P2} ,$  one can easily see that $C_{\alpha,q,\gamma}\leq C_{\widetilde{E}(0)},$ from where one concludes directly \eqref{ineq}.
\end{remark} 

Some prompt consequences of  Proposition $ \ref{P2}  $ are given below.

\begin{Corollary} [Dissipativity]\label{dissipative}
	Under the assumptions of  Proposition $ \ref{P2} ,$ 
	the dynamical system
	$(\mathcal{H},S_\lambda(t))$  given by \eqref{sist-din}  is dissipative, that is, it has a bounded absorbing set $\mathcal{B}\subset \mathcal{H}$, which is uniformly bounded w.r.t. $\lambda\in[0,1]$. In particular, 
	there exists a   positively invariant bounded absorbing set.
\end{Corollary}
\begin{proof}
	From \eqref{inequality_pricipal} and (\ref{energy_norm}), we obtain
	$$
	||(u(t),u_t(t))||^2_{\mathcal{H}}\le \frac{4}{\omega}\left[\frac{q}{C_{\widetilde{E}(0)}}(t-1)^{+}+\frac{1}{\big[\widetilde{E}(0)\big]^{q}}\right]^{-1/q}+\frac{32K_1}{\omega},  \quad t>0,
	$$
	for all  $\lambda\in[0,1]$.
	Thus, given an  arbitrary bounded set $B\subset \mathcal{H}$ and taking  $(u_0,u_1)\in B$,  there exists a time $t_B>0$ such that 
	\begin{equation}\label{bdd}
	||(u(t),u_t(t))||^2_{\mathcal{H}} \le \frac{64K_1}{\omega}, \quad \forall \, t\geq t_B.
	\end{equation}
	Therefore,
	$
	\mathcal{{B}}=\overline{B\left(0,8 \sqrt{\frac{ K_1}{\omega}} \right)}^{\mathcal{H}}
	$
	constitutes a bounded absorbing set (uniformly  w.r.t. $\lambda\in[0,1]$) for  $(\mathcal{H},S_\lambda(t))$. The construction of a positively invariant bounded absorbing set is standard.
\end{proof}

\begin{Corollary}[Uniform Global  Boundedness] \label{cor-bdd}	Under the assumptions of  Proposition $ \ref{P2} ,$ 
	the trajectory solutions of problem  \eqref{P} are  globally bounded in time (uniformly w.r.t.  $\lambda\in[0,1]$)  for initial data lying in bounded sets. 
	More precisely, 
	given a bounded set $B\subset \mathcal{H}$ and initial data $(u_0,u_1)\in B$, then 
	there exists a constant $C_B>0$ (depending only on $B$) such that 
	$$
	||S_\lambda(t)(u_{0},u_{1})||_{\mathcal{H}} \, =\,
	||(u(t),u_t(t))||_{\mathcal{H}} \,\le\, C_B, \quad \forall \, t\geq0.
	$$			
	\begin{proof}
		It is a direct consequence of Corollary \ref{dissipative} (see \eqref{bdd}) and   Theorem \ref{theo-existence}-$(i)$. 
	\end{proof}

\end{Corollary}

\begin{Corollary}[Gradient Property]\label{cor-gradient}
	Under the assumptions of  Proposition $ \ref{P2} ,$
	the dynamical system
	$(\mathcal{H},S_\lambda(t))$  given by \eqref{sist-din}  is gradient, that is, there exists a strict Lyapunov functional $\Phi_\lambda:=\Phi$ for 	$(\mathcal{H},S_\lambda(t))$. Moreover, the Lyapunov functional $\Phi$ is bounded from above on any
	bounded subset of $\mathcal{H}$ and the set  
	$$\overline{\Phi}_R=\{U\in \mathcal{H} \ ; \ \Phi(U)\le R\}
	$$
	is bounded in $\mathcal{H}$ for every $R>0$.
\end{Corollary}
\begin{proof}
	Let us define $\Phi:=E$. From 	
	\eqref{energy_relation_AA} one sees that the mapping 
	$$
	t\mapsto E(u(t),u_t(t))=\Phi(S_\lambda(t)U_0)
	$$  
	is non-increasing for every $U_0:=(u_0,u_1)\in\mathcal {H}$. Additionally,   from \eqref{energy_relation} and \eqref{absorbing_BB} one gets
	\begin{equation}\label{inequality grad}
	\Phi(S_\lambda(t)U_0)
	+\gamma\int_0^t \|
	u_t(\tau)\|^{2(q+1)} \,d\tau = \Phi(U_0), \quad  t>0,
	\end{equation}
	for every $U_0\in\mathcal {H}$. From \eqref{inequality grad}, it   easily concludes that
	$$
	\Phi(S_\lambda(t)U_0)=\Phi(U_0) \quad  \Rightarrow \quad U_0\in \mathcal{N}_\lambda,     \quad  t>0,
	$$
	where $\mathcal{N}_\lambda$ is defined in \eqref{state-set}. Since we know that 
	$$
	\quad U_0\in \mathcal{N}_\lambda \quad  \Leftrightarrow \quad  S_\lambda(t)(U_0)=U_0,   \quad  t>0,
	$$ 
	then   $\Phi$ is a
	strict Lyapunov functional for the dynamical system  	$(\mathcal{H},S_\lambda(t))$. 
	
	Moreover, from \eqref{inequality grad} we have $\Phi(S_\lambda(t)U_0)\leq\Phi(U_0)$ and, therefore, it is trivial to conclude that  $\Phi$ is bounded from above on bounded subsets of
	$\mathcal {H}$.  Finally,  if $\Phi(S_\lambda(t)U_0)\leq R,$ then 
	in view of (\ref{energy_norm}) we obtain
	that
	$S_\lambda(t)U_0=(u(t),u_t(t))$ satisfies
	$$
	\|S_\lambda(t)U_0\|^2_{\mathcal {H}}\leq \frac{4}{\omega}(R+K_1), \quad t\geq0, \ \lambda\in[0,1].
	$$
	Hence, $\overline{\Phi}_R$ is a bounded set of $\mathcal {H}$  for every $R>0$.
\end{proof}

\begin{Corollary}[Polynomial Decay Range]\label{P3} 
	Under the assumptions of  Proposition $ \ref{P2} ,$ and    additionally assuming that $h\equiv 0$    and $C_f=0$ in \eqref{assumption_f}, then	
	energy $E(t)$ satisfies
	\begin{equation}\label{polin-decay}
	\left[\frac{q}{C_{\alpha,q,\gamma}}t+\frac{1}{\big[E(0)\big]^{q}}\right]^{-1/q}\le E(t)\le\left[\frac{q}{C_{E(0)}}(t-1)^{+}+\frac{1}{\big[E(0)\big]^{q}}\right]^{-1/q},
	\end{equation}
	for all $t>0$.
\end{Corollary}
\begin{proof}
	It is a promptly consequence of \eqref{inequality_pricipal} with $K_\lambda=0$ and   $\widetilde{E}(t)=E(t)$ in \eqref{modified-energy}.
\end{proof}

\begin{remark}\label{rem-sque-polin}
	Under the conditions of 	Corollary $ \ref{P3} $, the homogeneous problem  related \eqref{P} is polynomially stable with rate $1/q$. Moreover, \eqref{polin-decay} shows  
	that the energy $E(t)$  	is  squeezed  in a polynomial decay  range with optimal rate  $1/q$ in the sense that it cannot be improved (by using \eqref{polin-decay}). Therefore, under these circumstances, problem  \eqref{P}  is never exponential stable in the homogeneous scenario. 
\end{remark}

\subsection{The set of stationary solutions}

\begin{lemma}\label{lem4-2}
	Under the assumptions of Theorem  $\ref{theo-existence} $, the
	set 
	\begin{equation}\label{state-set-lemma}
	\mathcal{N}_\lambda=\left\{ (u,0)\in\mathcal{H}; \ \kappa A u +
	A_1 u +f(u)= h_\lambda \right\}, \ \ \lambda\in[0,1],
	\end{equation}
	is bounded in $\mathcal {H}$, uniformly w.r.t. $\lambda\in[0,1]$. In particular, if $h\equiv 0$    and $C_f=0$ in \eqref{assumption_f}, then
	$\mathcal{N}_0=\{(0,0)\}$.
\end{lemma}
\begin{proof}  
	The existence of (at least one non-trivial) solution $u^\lambda:=u$ for the equation	
	\begin{equation}\label{elipt}
	\kappa A u +
	A_1 u +f(u)= h_\lambda \quad \mbox{in} \quad H,
	\end{equation}
	will be given in the next result. Taking the inner product in $H$ of  \eqref{elipt} with $u$, we have
	\begin{equation}\label{est-sol}
	\|A^{1/2}_1 u\|^2 +
	\kappa\|A^{1/2} u\|^2=- \big(f(u),u\big)+\big(h_\lambda ,u\big).
	\end{equation}
	From assumption (\ref{assumption_f}), we have
	\begin{equation*}
	-\big(f(u),u\big)\leq 
	C_f|\Omega|+{c_f} \|u\|^2 \leq
	C_f|\Omega|+\frac{c_f}{\sigma_1} \|A^{1/2}_1 u\|^2.
	\end{equation*}
	Also, from  Young's inequality with $\omega=1-\frac{c_f}{\sigma_1}>0$, we infer
	$$
	\big(h_\lambda ,u\big) \leq \frac{\omega}{2}\|A^{1/2}_1 u\|^2+\frac{2}{\sigma_1\omega}\|h_\lambda\|^2.
	$$
	Going back to \eqref{est-sol} we obtain
	\begin{equation}\label{bdd-sol}
	\frac{\omega}{2} \|A^{1/2}_1 u\|^2 +\kappa\|A^{1/2} u\|^2 \leq C_f|\Omega|+\frac{2}{\sigma_1\omega}\|h_\lambda\|^2,
	\end{equation}
	from where one concludes that $\mathcal {N}_\lambda$ is
	bounded  in $\mathcal {H}$, uniformly w.r.t. $\lambda\in[0,1]$. 
	In particular, if $h\equiv 0$    and $C_f=0$, 
	then \eqref{bdd-sol} also implies that $u=0$ and thus $\mathcal{N}_0$ is the  trivial null set.
\end{proof}

\begin{lemma}\label{lem-nontrivial}
	Under the above notations
	and 	 assumptions of Theorem  $\ref{theo-existence} $, if $C_f>0$ in \eqref{assumption_f} and even if $h\equiv 0$, 
	then it is possible to show that the set
	$$
	\mathcal{N}_0=\left\{ (u,0)\in\mathcal{H}; \ \kappa A u +
	A_1 u +f(u)= 0 \right\}
	$$
	has a nontrivial  (weak) solution $u\neq0$ in $H$.
	Therefore,  it is possible to conclude that the set    $\mathcal{N}_0$ (and more generally  $\mathcal{N}_\lambda$) has at least two stationary solutions.
\end{lemma}
\begin{proof}
	Since $f(0)=0$, then obviously  $u=0$ is the trivial  solution of 
	\begin{equation}\label{elipt-null}
	\kappa A u +
	A_1 u +f(u)= 0 \quad \mbox{in} \quad H.
	\end{equation}
	
	In what follows, let us deal with the case of nontrivial  weak solution  for  $u$ for \eqref{elipt-null}. To fix the ideas, 
	we are going to assume, by simplicity, 
	the following concrete example for $f$
	\begin{equation}\label{f-concrete}
	f(s)=|s|^{\delta}s-\sigma|s|^{r}s, \ \  \sigma>0, \, 0 < r < \delta \leq p.
	\end{equation}
	The Euler-Lagrange functional $I_f: \mathcal{D}(A_1^{1/2}) \to \mathbb{R}$ corresponding to \eqref{elipt-null}-\eqref{f-concrete}
	is given by
	$$
	I_f(u) = \frac{1}{2}\|A^{1/2}_1 u\|^2+\frac{\kappa}{2}\|A^{1/2} u\|^2 +\frac{1}{\delta+2} \| u\|^{\delta+2}_{\delta+2}-
	\frac{\sigma}{r+2}\|u\|^{r+2}_{r+2}.
	$$

	We claim that for all $\sigma > 0$, $I_f$ is coercive and bounded from below. In fact, 
	since $0 < r < \delta$, then from the embedding $L^\delta\hookrightarrow L^r$ (with constant $C>0$) we get
	\begin{eqnarray}
	\nonumber I_f(u) & \geq & \frac{1}{2}\|A^{1/2}_1 u\|^2 + \frac{1}{\delta+2} \| u\|^{\delta+2}_{\delta+2} -  \frac{\sigma C}{r+2} \|u\|_{\delta+2}^{r+2}\\
	\label{coervice} & \geq &  \frac{1}{2}\|A^{1/2}_1 u\|^2 + D_0,
	\end{eqnarray}
	where
	$$
	D_0 = \inf_{\tau \geq 0}\left\{\frac{\tau^{\delta+2}}{\delta+2} - \sigma C\frac{\tau^{r+2}}{r+2}\right\}.
	$$
	Then, because of (\ref{coervice}), one has that  $I_f$ is clearly coercive\footnote{Given a Banach space $(X,||\cdot||)$, we recall that a functional $I:X\to\mathbb{R}$ is {\it coercive} if $\left\|u_{n}\right\| \rightarrow \infty$ implies $I\left(u_{n}\right) \rightarrow \infty$.} and bounded from below.

	On the other hand, by   fixing  $0\neq u \in \mathcal{D}(A_1^{1/2})$,  and
	regarding the embedding chain 
	$\mathcal{D}(A_1^{1/2})\hookrightarrow\mathcal{D}(A) \hookrightarrow L^\delta\hookrightarrow L^r$,
	we observe that  there exists a parameter $\sigma_0 > 0$ such that $I_f(u) < 0$. For such a number 
	$\sigma_0$, we consider a
	minimizing sequence, namely, a sequence $(u_n) \subset \mathcal{D}(A_1^{1/2})$ such that
	$$
	\lim_{n \to +\infty} I_f(u_n) \ =  \inf_{v\in\mathcal{D}(A_1^{1/2})}I_f(v):= \xi.
	$$
	From the coerciveness of $I_f$, we have that $(u_n)$ is bounded, and passing to a sub-sequence if necessary,  
	we infer  $u_n \to u$ weakly in $\mathcal{D}(A_1^{1/2})$. Due to the compactness of the embeddings $\mathcal{D}(A_1^{1/2})
	\hookrightarrow L^{\delta+2}$ and $\mathcal{D}(A_1^{1/2})\hookrightarrow\mathcal{D}(A)\hookrightarrow L^{r+2}$ (once $r+2 < \delta +2\leq p^*$), then
	\begin{equation*}
	\xi \, \leq \, I_f(u)
	\, \leq \, \liminf_{n \to +\infty}I_f(u_n)
	\, = \, \xi.
	\end{equation*}
	This implies that 
	$u$ is a global minimizer and, therefore, a nontrivial critical point  of $I_f$, which in turn corresponds to 
	a weak solution of (\ref{elipt-null}) as desired.
\end{proof}

\subsection{Key estimates w.r.t. critical parameters $\alpha$ and $p^*$}

The next result is a generalized version the one presented in \cite[Lemma 2.1]{Nakao}. 
It will be useful in the next key stability inherent to our dynamical system.

\begin{proposition}[Nakao's Generalized Lemma]\label{Nak-lemma}
	Let $\phi(t)$ be a non-negative continuous and $K(t)$ be a non-negative non-decreasing   functions on $[0, T), T>1,$ possibly $T=\infty,$  such that
	\begin{equation}\label{nak-hyp}
	\sup _{t \leqslant s \leqslant t+1} [\phi(s)]^{1+\rho} \leqslant C_{0}(\phi(t)-\phi(t+1))+K(t), \quad 0 \leqslant t \leqslant T-1,
	\end{equation}
	for some $C_{0}>0$ and $\rho\geq0 .$ Then, the following estimates hold:
	\begin{enumerate}[label={\rm(\textbf{N.\arabic*})}]
		
		\item\label{N1-nak} For $\rho>0$, we have:
		\begin{equation}\label{nak-thesis1}
		\phi(t) \leqslant\left(C_{0}^{-1} \rho(t-1)^{+}+\left(\sup _{0 \leqslant s \leqslant 1} \phi(s)\right)^{-\rho}\right)^{-1 / \rho}+\big[K(t)\big]^{1 /(\rho+1)}, \quad 0 \leqslant t<T,
		\end{equation}
		where we consider	usual notation $s^+:=(s+|s|)/2$.
		
		\item\label{N2-nak} For $\rho=0$, we have:
		\begin{equation}\label{nak-thesis2}
		\phi(t) \leqslant \sup _{0 \leqslant s \leqslant 1} \phi(s)\left(\frac{C_{0}}{1+C_{0}}\right)^{[t]}+K(t), \quad 0 \leqslant t<T,
		\end{equation}
		where   $[s]$ stands for the largest integer less than or equal to  $s\geq0$. 
	\end{enumerate}
	
\end{proposition}
\begin{proof}
	It follows similar patterns as done in \cite[Lemma 2.1]{Nakao}  with proper modifications in what concerns the function $K(t)$, which in \cite{Nakao} is only assumed as a positive constant $K(t):=K>0$. Here, we address a more general class of functions by  assuming that  $K(t)$ can be a non-decreasing function, instead of a positive constant only.

	For the sake of completeness, we present a detailed proof in  Appendix \ref{sec-proof-nak-lemma}.
	
\end{proof}

\begin{remark}
	Given  a measurable  function $g:[0,t+1]\to \mathbb{R}$,  $ t\geq 0$, we observe that for any $a,b \in [t, t+1]$ with $a\leq b$, it holds the following estimate
	\begin{equation}\label{dom-lema}
	\left|\int_a^b g(s)\, ds\right|\leq \sup_{s \in [0,t+1]}\sup_{r \in [0,s]}\left|\int_r^s g(\tau)\, d\tau\right|.
	\end{equation}	
	Indeed,  to reach it we simply note that it holds the chain of inequalities
	\begin{equation*}
	\left|\int_a^b g(\tau)\, d\tau\right|  \leq \sup_{s \in [a,t+1]}\left|\int_a^s g(\tau)\, d\tau\right|\\
	\leq  \sup_{s \in [a,t+1]}\sup_{r \in [a,s]}\left|\int_r^s g(\tau)\, d\tau\right|\\
	\leq \sup_{s \in [0,t+1]}\sup_{r \in [0,s]}\left|\int_r^s g(\tau)\, d\tau\right|.
	\end{equation*}	
\end{remark}

The next result provides a key inequality in the present article. It gives a suitable estimate encompassing the difference of two
trajectory solutions of \eqref{P}. To its proof, we combine refined and new arguments  along with the  above Nakao  generalized result and 
\eqref{dom-lema}.

Before stating such a key result, let us remark that by virtue of Corollary \ref{cor-bdd} any 
trajectory solution is globally bounded in time on bounded sets (uniformly w.r.t.  $\lambda\in[0,1]$).  
This fact will be highly used in the next result. Also, in what follows we will denote by $C_B>0$ several different constants  depending on a general bounded set $B\subset\mathcal{H}$.
\begin{proposition}[Key Inequality]\label{theo_inequality} 
	Under the assumptions of Theorem $\ref{theo-existence}$,  let us also consider a bounded set $B \subset {\mathcal{H}}$.
	Given  $U^i=(u^i_{0},u^i_{1})\in B, \, i=1,2, $ we denote by $S_\lambda(t)U^i=(u^i(t),u^i_{t}(t)), i=1,2,$ the respective trajectory solution
	corresponding to the dynamical system \eqref{sist-din}. Then, there exists a constant $C_B>0$ such that
	\begin{align}\label{key-inequality}
	||S_\lambda(t)U^1-S_\lambda(t)U^2||_{\mathcal{H}}^2\le &\;\left[C_B^{-1}q(t-1)^{+}+\left(\sup_{0\le s\le
		1}||(w(s),w_t(s))||_{\mathcal{H}}^2\right)^{-q}\right]^{-\frac{1}{q}}\nonumber\\
	&\;+C_{B}\sup_{0\le s\le t+1}\left[\|A^{\alpha}w(s)\|^{\frac{2(q+1)}{2q+1}}+\|w(s)\|^{2}_{p+2}\right]^{\frac{1}{q+1}}\\
	&\;+\left[\,J_f(w(t),w_t(t))\,\right]^{\frac{1}{q+1}}.\nonumber
	\end{align}
	for all $t>0$ and $\lambda\in[0,1]$, where we  set  $w=u^1-u^2$,  $F(w)=f(u^1)-f(u^2)$, and	
	\begin{align}\label{critical}
	J_f(w(t),w_t(t)) = & \ 2 
	\sup_{s \in [0,t+1]}\sup_{r \in [0,s]} \left|\int_{r}^{s}(F(w(\tau)),w_t(\tau))d\tau\right|
	\nonumber \\
	& +4^{q+1}\sup_{s \in [0,t+1]}\sup_{r \in [0,s]} \left|\int_{r}^{s}(F(w(\tau)),w_t(\tau))d\tau\right|^{q+1} .
	\end{align}
\end{proposition}
\begin{proof}
	Since we are dealing with the difference $S_\lambda(t)U^1-S_\lambda(t)U^2:=(w(t),w_t(t))$ of two trajectory solutions of problem $(\ref{P})$, then we are going to borrow some notations used 	
	in Theorem \ref{theo-dependence}. Indeed, we first note that  $(w,w_t)$ is a solution (in the weak and strong sense) of 	  problem (\ref{dependenceI}). Also, the functional $\mathcal{E}_w$ defined in \eqref{def_V} satisfies 
	\begin{align}\label{equiv}
	||(w(t),w_t(t))||_{\mathcal{H}}^2 \, \leq \, \mathcal{E}_w(t) \, \leq \,
	\mu_\kappa||(w(t),w_t(t))||_{\mathcal{H}}^2,
	\end{align} 
	form some constant $\mu_\kappa:=1+\kappa
	\mu_0>0$, see \eqref{dependenceIII}, and
	the identity 
	\begin{eqnarray}\label{dependencexxx}
	\frac{d}{dt}\mathcal{E}_w(t)+\gamma\Pi_1(t)\|
	w_t(t)\|^2=-\,\gamma\Pi_2(t)\left(u^1_t(t)+
	u^2_t(t),w_t(t)\right)+2\left(F(w(t)),w_t(t)\right),
	\end{eqnarray}
	which is essentially \eqref{dependenceII}. In what follows, our first goal is to estimate the right-hand side of \eqref{dependencexxx}.
	
	Firstly, using that $[a^{2q}+b^{2q}]\ge [a-b]^{2q}$, we get
	\begin{equation}\label{key-1}
	\gamma\Pi_1(t)\|w_t(t)\|^2\ge\gamma\left[\,\|u^1_t(t)\|^{2q}+\|u^2_t(t)\|^{2q}\,\right]\|w_t(t)\|^2\ge\frac{\gamma}{2^{2q}}\|w_t(t)\|^{2(q+1)}.
	\end{equation}

 Now, for $q\geq1/2$, we claim that
\begin{equation}\label{estimate-Pi2}
-\gamma\Pi_2(t)\left(u^1_t(t)+
	u^2_t(t),w_t(t)\right) \leq C_{B}
	\|A^{\alpha}w(t)\|^{\frac{2(q+1)}{2q+1}}+\frac{\gamma}{2^{2(q+1)}}\|w_t(t)\|^{2q+1}, 
	\end{equation}
for some constant $C_B>0$. Indeed, if $||S_\lambda(t)U^1||_{\mathcal{H}^{\alpha}}=||S_\lambda(t)U^2||_{\mathcal{H}^{\alpha}}=0$, there is nothing to do. Let us suppose then 
$||S_\lambda(t)U^1||_{\mathcal{H}^{\alpha}}+||S_\lambda(t)U^2||_{\mathcal{H}^{\alpha}}>0$.	From the Mean Value Theorem (MVT for short), there
	exists a number
	$$\xi_{\vartheta}=\vartheta
	||S_\lambda(t)U^1||_{\mathcal{H}^{\alpha}}+(1-\vartheta)||S_\lambda(t)U^2||_{\mathcal{H}^{\alpha}},\quad
	\vartheta\in (0,1),$$ such that
	\begin{eqnarray*}
		\Pi_2(t)&=& 2q\int_0^1\left|\xi_{\vartheta}\right|^{2q-1}d\vartheta\left[\,||S_\lambda(t)U^1||_{\mathcal{H}^{\alpha}}-||S_\lambda(t)U^2||_{\mathcal{H}^{\alpha}}\,\right]\\
		&=&2q\int_0^1\left|\xi_{\vartheta}\right|^{2q-1}d\vartheta\left\{\frac{\left[\,\|u^1_t(t)\|^2-\|u^2_t(t)\|^2\,\right]+\left[\,\|A^{\alpha}u^1(t)\|^2-\|A^{\alpha}u^2(t)\|^2\,\right]}{||S_\lambda(t)U^1||_{\mathcal{H}^{\alpha}}+||S_\lambda(t)U||_{\mathcal{H}^{\alpha}}}\right\}.
	\end{eqnarray*}
	Thus, we can write
	\begin{equation}\label{key-2}
	-\gamma\Pi_2(t)\left(u^1_t(t)+u^2_t(t),w_t(t)\right)=\eta(t)+\chi(t),
	\end{equation}
	where
\begin{eqnarray*}
\eta(t) & = & -2q\gamma\int_0^1\left|\xi_{\vartheta}\right|^{2q-1}d\vartheta\frac{\left[\,\|u^1_t(t)\|^2-\|u^2_t(t)\|^2\,\right]^2}{||S_\lambda(t)U^1||_{\mathcal{H}^{\alpha}}+||S_\lambda(t)U||_{\mathcal{H}^{\alpha}}}\\
\chi(t) &=& -2q\gamma\int_0^1\left|\xi_{\vartheta}\right|^{2q-1}d\vartheta
	\frac{\left[\,\|A^{\alpha}u^1(t)\|^2-\|A^{\alpha}u^2(t)\|^2\,\right]}{||S_\lambda(t)U^1||_{\mathcal{H}^{\alpha}}+||S_\lambda(t)U||_{\mathcal{H}^{\alpha}}}\left(u^1_t(t)+	u^2_t(t),w_t(t)\right).	
\end{eqnarray*}	
	Using Young inequality with $\frac{2q+1}{2(q+1)}+\frac{1}{2(q+1)}=1,$ the term $\chi(t)$ can be estimated as follows
	\begin{eqnarray*}
	|\chi(t)|&\leq & C_B\frac{\left[\,||S_\lambda(t)U^1||^2_{\mathcal{H}^{\alpha}}+||S_\lambda(t)U^2||^2_{\mathcal{H}^{\alpha}}\,\right]}{||S_\lambda(t)U^1||_{\mathcal{H}^{\alpha}}+||S_\lambda(t)U||_{\mathcal{H}^{\alpha}}}\|A^{\alpha}w(t)\|\|w_t(t)\|\nonumber\\
	&\leq &C_{B}
	\|A^{\alpha}w(t)\|^{\frac{2(q+1)}{2q+1}}+\frac{\gamma}{2^{2(q+1)}}\|w_t(t)\|^{2(q+1)}, 
	\end{eqnarray*}
for some constant $C_B>0$, once $2q-1\geq0$. Noting that $\eta(t)\leq 0$ and connecting the last estimate in \eqref{key-2}, we arrive at \eqref{estimate-Pi2}.
	
	Thus, collecting \eqref{key-1} and \eqref{estimate-Pi2} with \eqref{dependencexxx} , we obtain  
	\begin{equation}\label{dependencexxxx}
	\frac{d}{dt}\mathcal{E}_w(t)+\frac{\gamma}{2^{2q+1}}\|
	w_t(t)\|^{2(q+1)}\le C_{B}\left[\,\|A^{\alpha}w(t)\|^{\frac{2(q+1)}{2q+1}}+2\left(F(w(t)),w_t(t)\right)\,\right].
	\end{equation}
	Integrating (\ref{dependencexxxx}) from $t$ to $t+1$, we have
	\begin{align}\label{asymptotic_z}
	&\frac{\gamma}{2^{2q+1}}\int_t^{t+1}\|
	w_t(s)\|^{2(q+1)}ds\\
	&\quad\le\; \mathcal{E}_w(t)-\mathcal{E}_w(t+1)+ C_{B}\int_{t}^{t+1}\|A^{\alpha}w(s)\|^{\frac{2(q+1)}{2q+1}}ds+2\left|\int_{t}^{t+1}(F(w(s)),w_t(s))ds\right|\nonumber\\
	&\quad:=\;\left[\,G(t)\,\right]^2.\nonumber
	\end{align}
	Now, from H\"older's inequality with $\frac{q}{q+1}+\frac{1}{q+1}=1$ and (\ref{asymptotic_z}), we have
	\begin{align}\label{asymptotic_zz}
	\int_t^{t+1}\|
	w_t(s)\|^{2}ds\le \left[\int_t^{t+1}\|w_t(s)\|^{2(q+1)}ds\right]^{\frac{1}{q+1}}\le \frac{2^{\frac{2q+1}{q+1}}}{\gamma^{\frac{1}{q+1}}}\left[\,G(t)\,\right]^{\frac{2}{q+1}},
	\end{align}
	which implies (from the MVT for Integrals) that there exists $t_1\in [t,t+\frac{1}{4}]$ , $t_2\in [t+\frac{3}{4},t+1]$ such that
	\begin{align}\label{asymptotic_zzz}
	\|
	w_t(t_i)\|^{2}\le 4\left[\int_t^{t+1}\|w_t(s)\|^{2}ds\right]\le \frac{2^{\frac{4q+3}{q+1}}}{\gamma^{\frac{1}{q+1}}}\left[\,G(t)\,\right]^{\frac{2}{q+1}}.
	\end{align}

	Let us keep the above estimates in mind to apply them properly in the next computations.
	
	Working on the other hand, we take now the multiplier $w$ in (\ref{dependenceI}) and 
	integrating from $t_1$ to $t_2$, we have
	\begin{align}\label{asymptotic_zzzz}
	\int_{t_1}^{t_2} \mathcal{E}_w(t)\, ds=2\int_{t_1}^{t_2}\|w_t(s)\|^2ds+\sum_{i=1}^4\mathcal{L}_i,
	\end{align}
	where
	\begin{eqnarray*}
		\mathcal{L}_1&=&-\,\int_{t_1}^{t_2}\left(F(w(s)),w(s)\right)ds,\\
		\mathcal{L}_2&=&-\,\left[\,\left(w_t(t_2),w(t_2)\right)-\left(w_t(t_1),w(t_1)\right)\,\right],\\
		\mathcal{L}_3&=&-\,\frac{\gamma}{2}\int_{t_1}^{t_2}\Pi_1(s)\left( w_t(s),w(s)\right)ds,\\
		\mathcal{L}_4&=&-\,\frac{\gamma}{2}\int_{t_1}^{t_2}\Pi_2(s)\left(\left[\,
		u^1_t(s)+ u^2_t(s)\,\right],w(s)\right)ds.
	\end{eqnarray*}
	The terms $\mathcal{L}_i$, $i=1,\cdots,4$ can be estimated as
	follows. Firstly, from MVT, the hypothesis 
	(\ref{assumption_f'}), H\"older's inequality and the embedding $\mathcal{D}(A_1^{1/2})\hookrightarrow L^{p+2}$, we have
	\begin{align*}
	\mathcal{L}_1&\le\int_{t_1}^{t_2}\left(\,\left[f(u^1)-f(u^2)\right],\,w\right)ds\\
	&\le
	C_{f'}\int_{t_1}^{t_2}\int_{\Omega}\left[\,(1+|u_1|+|u^2|)^{p}\,\right]|w|^2dxds\\
	&\le C_{f'}\int_{t_1}^{t_2}\left[\int_{\Omega}(1+|u^1|^{p}+|u^2|^{p})^{\frac{p+2}{p}}\right]^{\frac{p}{p+2}}\|w(s)\|^2_{p+2}ds\\
	&\le C_{B}\int_{t_1}^{t_2}\|w(s)\|^2_{p+2}ds,
	\end{align*}
	for some constant $C_B>0$. Also, using now the embedding  $\mathcal{D}(A_1^{1/2})\hookrightarrow H$, Young's inequality, (\ref{asymptotic_zz})-(\ref{asymptotic_zzz}), we get
	\begin{align*}
	\mathcal{L}_2&\le \frac{1}{\sigma_1^{1/2}}\sum_{i=1}^2\|w_t(t_i)\|\|A^{1/2}_1w(s)\| \\
	&\le \frac{2}{\sigma_1^{1/2}}\left\{\frac{2^{\frac{4q+3}{2(q+1)}}}{\gamma^{\frac{1}{2(q+1)}}}\left[\,G(t)\,\right]^{\frac{1}{q+1}}\right\}\sup_{t_1\le s\le t_2}\|A^{1/2}_1w(s)\|\nonumber\\
	&\le \frac{2^{\frac{4q+3}{q+1}}}{\delta \sigma_1\gamma^{\frac{1}{q+1}}}\left[\,G(t)\,\right]^{\frac{2}{q+1}}+\delta\sup_{t_1\le s\le t_2}\left[ \mathcal{E}_w(t)\right],
	\end{align*}
	and
	\begin{equation*}
	\mathcal{L}_3\le C_B\left(\int_{t_1}^{t_2}\|
	w_t(s)\|ds\right)\sup_{t_1\le s\le t_2}\|A^{1/2}_1w(s)\|\le
	C_B\left[\,G(t)\,\right]^{\frac{2}{q+1}}+\delta\sup_{t_1\le s\le t_2}\left[\mathcal{E}_w(t)\right],
	\end{equation*}
	for some constant $C_B>0$.
	Third,   using Lemma \ref{Lemma-Haraux} and noting that 
	$$\mathcal{D}(A^{1/2}_1)\hookrightarrow \mathcal{D}(A^{\alpha}),\quad L^{p+2}(\Omega)\hookrightarrow H,$$
	we obtain for $q\geq1/2$ that
	\begin{eqnarray*}
		\mathcal{L}_4&\le& \gamma q \int_{t_1}^{t_2}\max\{||S_\lambda(s)U^1||^2_{\mathcal{H}^{\alpha}},||S_\lambda(s)U^1||^2_{\mathcal{H}^{\alpha}}\}^{2q-1}[\mathcal{E}_w(t)]^{1/2}\|w(s)\|ds\\
		& \leq & C_B\int_{t_1}^{t_2} [\mathcal{E}_w(t)]^{1/2}\|w(s)\|ds\\
		&\le &\delta \sup_{t_1\le s \le t_2}[\mathcal{E}_w(s)]+C_B\int_{t_1}^{t_2}\|w(s)\|^2_{p+2}ds,
	\end{eqnarray*}
	for some constant $C_B>0$ and $\delta>0$. 
	
	Replacing the latter estimates for  $\mathcal{L}_1,\cdots,\mathcal{L}_4$ in
	(\ref{asymptotic_zzzz}), we arrive at
	$$
	\int_{t_1}^{t_2}\left[\mathcal{E}_w(s)\right]ds \le C_{B}\left[\,G(t)\,\right]^{\frac{2}{q+1}}+3\delta\sup_{t_1\le s\le t_2}\left[\mathcal{E}_w(s)\right]+\,C_B\int_{t}^{t+1}\|w(s)\|^{2}_{p+2}\,ds,
	$$
	for some constant $C_B>0$. Again from the MVT, there exists $\tau_1\in [t_1,t_2]\subset [t,t+1]$ such that
	\begin{align}\label{z}
	\mathcal{E}_w(\tau_1)&\le C_{B}\left[\,G(t)\,\right]^{\frac{2}{q+1}}+3\delta\sup_{t_1\le s\le t_2}\left[\mathcal{E}_w(s)\right]+C_B\int_{t}^{t+1}\|w(s)\|^{2}_{p+2}\,ds.
	\end{align}
	Let us also consider $\tau_2\in [t,t+1]$ such that
	$$\mathcal{E}_w(\tau_2):=\sup_{t\le s\le t+1}\left[\mathcal{E}_w(s)\right].$$
	Now, integrating \eqref{dependencexxxx} over $[t,\tau_2]$ and over $[\tau_1,t+1]$, using \eqref{z} and noting that
	$$\mathcal{E}_w(t)\le \mathcal{E}_w(t+1)+\left[\,G(t)\,\right]^{2},$$
	we obtain
	\begin{align*}
	\mathcal{E}_w(\tau_2)&\le \mathcal{E}_w(t)+C_{B}\int_{t}^{t+1}\|A^{\alpha}w(s)\|^{\frac{2(q+1)}{2q+1}}ds+2\left|\int_{t}^{\tau_2}(F(w(s)),w_t(s))ds\right|\\
	&\le \mathcal{E}_w(t+1)+\left[\,G(t)\,\right]^{2}+C_{B}\int_{t}^{t+1}\|A^{\alpha}w(s)\|^{\frac{2(q+1)}{2q+1}}ds+2\left|\int_{t}^{\tau_2}(F(w(s)),w_t(s))ds\right|\\
	&\le \mathcal{E}_w(\tau_1)+\left[\,G(t)\,\right]^{2}+C_{B}\int_{t}^{t+1}\|A^{\alpha}w(s)\|^{\frac{2(q+1)}{2q+1}}ds\\
	&\quad+2\left|\int_{\tau_1}^{t+1}(F(w(s)),w_t(s))ds\right|+2\left|\int_{t}^{\tau_2}(F(w(s)),w_t(s))ds\right|\\
	&\le \left[\,G(t)\,\right]^{2}+C_{B}\left[\,G(t)\,\right]^{\frac{2}{q+1}}+3\delta\mathcal{E}_w(\tau_2)\\
	&\quad+C_{B}\int_{t}^{t+1}\|A^{\alpha}w(s)\|^{\frac{2(q+1)}{2q+1}}ds+C_B\int_{t_1}^{t_2}\|w(s)\|^{2}_{p+2}\,ds\\
	& \quad +2\left|\int_{\tau_1}^{t+1}(F(w(s)),w_t(s))ds\right|+2\left|\int_{t}^{\tau_2}(F(w(s)),w_t(s))ds\right|,
	\end{align*}
	for some constant $C_B>0.$
	Choosing $\delta>0$ small enough, noting that $\left[\,G(t)\,\right]^{q+1} +C_B\le C_B,$ for some constant $C_B>0$, and using \eqref{dom-lema}, 
	we infer
	\begin{align}\label{123}
	\mathcal{E}_w(\tau_2)\le &\; C_B\left[\,G(t)\,\right]^{\frac{2}{q+1}}+C_{B}\int_{t}^{t+1}\|A^{\alpha}w(s)\|^{\frac{2(q+1)}{2q+1}}ds+C_B\int_{t}^{t+1}\|w(s)\|^{2}_{p+2}\,ds\nonumber\\
	& \; +2\left|\int_{\tau_1}^{t+1}(F(w(s)),w_t(s))ds\right|+2\left|\int_{t}^{\tau_2}(F(w(s)),w_t(s))ds\right| \\ 
	\leq & \  C_B\left[\,G(t)\,\right]^{\frac{2}{q+1}}+C_{B}\sup_{t\le s\le t+1}\big[\|A^{\alpha}w(s)\|^{\frac{2(q+1)}{2q+1}}+\|w(s)\|^{2}_{p+2}\big]\nonumber\\
	& \; +4\sup_{s \in [0,t+1]}\sup_{r \in [0,s]} \left|\int_{r}^{s}(F(w(\tau)),w_t(\tau))d\tau\right|.\nonumber
	\end{align}
	In this way, from \eqref{123}, the definition of $G(t)$ in (\ref{asymptotic_z}), and using again \eqref{dom-lema}, we arrive at
	\begin{equation}\label{abcd}
	\sup_{t\le s\le t+1}\left[\mathcal{E}_w(s)\right]^{q+1}\,\le \, C_B\left[\mathcal{E}_w(t)- \mathcal{E}_w(t+1)\,\right] + K_f(t),
	\end{equation}
	where $K_f(t) := K_f(w(t),w_t(t))$ is set by 
	\begin{align*}
	K_f(t) =& \ 2 
	\sup_{s \in [0,t+1]}\sup_{r \in [0,s]} \left|\int_{r}^{s}(F(w(s)),w_t(s))ds\right|
	\nonumber\\
	& +4^{q+1}\sup_{s \in [0,t+1]}\sup_{r \in [0,s]} \left|\int_{r}^{s}(F(w(s)),w_t(s))ds\right|^{q+1} \nonumber \\ 
	&+C_{B}\sup_{0\le s\le t+1}\big[\|A^{\alpha}w(s)\|^{\frac{2(q+1)}{2q+1}}+\|w(s)\|^{2}_{p+2}\big]
	\nonumber \\ 
	= &\ J_f(t) +C_{B}\sup_{0\le s\le t+1}\big[\|A^{\alpha}w(s)\|^{\frac{2(q+1)}{2q+1}}+\|w(s)\|^{2}_{p+2}\big],
	\end{align*}
	and $J_f(t):=J_t(w(t),w_t(t))$ is defined in \eqref{critical}.
	Hence, applying Proposition \ref{Nak-lemma} with $\phi:=\mathcal{E}_w$ and $K=K_f$ in \eqref{abcd}, by noting that $K_f(t)$ is a non-decreasing function, and regarding \eqref{equiv}, we conclude
	$$
	\mathcal{E}_w(t) \leqslant\left[C_{B}^{-1} q(t-1)^{+}+\left(\sup _{0 \leqslant s \leqslant 1} ||(w(s),w_t(s))||_{\mathcal{H}}^2\right)^{-q}\right]^{-1 / q}+\big[K_f(t)\big]^{1 /(q+1)}.
	$$
	Finally, since $0<\frac{1}{q+1}<1$, we have $|a+b|^{\frac{1}{q+1}}\leq |a|^{\frac{1}{q+1}}+|b|^{\frac{1}{q+1}}$, and then \eqref{key-inequality} holds true with $J_f(t)$ given in \eqref{critical}. 
\end{proof}

\begin{Corollary}[Stabilizability Estimate]\label{cor-key-subcritic}
	Under the same assumptions and statements of Proposition 
	$ \ref{theo_inequality} $, there exists 
	a constant $C_B>0$ such that
	\begin{align}\label{key-inequality-cor}
	||S_\lambda(t)U^1-S_\lambda(t)U^2||_{\mathcal{H}}^2\le &\;\left[C_B^{-1}q(t-1)^{+}+\left(\sup_{0\le s\le
		1}||(w(s),w_t(s))||_{\mathcal{H}}^2\right)^{-q}\right]^{-\frac{1}{q}}\nonumber\\
	&\;+C_{B}\sup_{0\le s\le t+1}\left[\|A^{\alpha}w(s)\|^{\frac{2(q+1)}{2q+1}}+||w(s)||_{p^*}^{\frac{2(q+1)}{2q+1}}\right]^{\frac{1}{q+1}}.
	\end{align}
	for all $t>0$ and $\lambda\in[0,1]$.
\end{Corollary}
\begin{proof}
	By means of the condition (\ref{assumption_f'}), MVT, 
	H\"older's inequality with  $\frac{p}{p^*}+\frac{1}{p^*}+\frac{1}{2}=1,$ Young's inequality with $\frac{2 q+1}{2(q+1)}+\frac{1}{2(q+1)}=1$, and the embedding $\mathcal{D}(A_{1}^{1 / 2}) \hookrightarrow L^{p^*},$ one can estimate
	\begin{equation*}
	|\left(F(w(t)),w_t(t)\right)| \leq C_B ||w(t)||_{p^*}^{\frac{2(q+1)}{2q+1}}+\frac{\gamma}{2^{2(q+1)}}\left\|w_{t}(t)\right\|^{2(q+1)}, \quad t \geq0, 
	\end{equation*}	
	for come constant $C_B>0$. Therefore,
	we go back to \eqref{dependencexxxx} and proceed verbatim the proof of Proposition 
	$ \ref{theo_inequality} $ to conclude
	\eqref{key-inequality-cor}, where we also note that $L^{p^*}\hookrightarrow L^{p+2}$. 
\end{proof}

\begin{remark}\label{rem-finite}
	Although Corollary $\ref{cor-key-subcritic}$ provides a 
	``milder'' stabilizability estimate
	with respect to $L^{p^*}$-norm, 
	it will be useful 
to reach 
	an estimate for Kolmogorov's
	$\varepsilon$-entropy of the 
	attractor $\mathfrak{A}_{\lambda}$ corresponding to 
	the dynamical system $(\mathcal{H},S_\lambda(t))$ in sub-critical aspects w.r.t. parameters $\alpha$ and $p^*$.
\end{remark}

The next result is an extended version of some results presented in \cite{Kham-Jde} for second-order wave problems. We have raised it to our concerns 
in higher-order Sobolev spaces on bounded domains, namely, in $\mathcal{H}=\mathcal{D}(A^{1/2}_1)\times H$, but we notify that its proof follows the same lines as in \cite[Section 2]{Kham-Jde}. It will be helpful in our next result on asymptotic smoothness.
\begin{proposition}\label{Khan-lemma}
	Let us consider $f$ satisfying 
	{Assumption $\ref{assumption}$}.	If 
	$\left\{\left(u^{n}, u_{t}^{n}\right)\right\}$ is  a weakly-star convergent  sequence in $L^{\infty}\left(s, T ; \mathcal{H}\right), \, 0\leq s<T,$ then
	\begin{equation}\label{Khan-limit}
	\lim _{n \rightarrow \infty} \lim _{m \rightarrow \infty} \int_{s}^{T}\big( f\left(u^{n}(t)\right)-f\left(u^{m}(t)\right), u_{t}^{n}(t)-u_{t}^{m}(t)\big) d t=0.
	\end{equation}
\end{proposition}
\begin{proof}
	It is similar to the statements 
	provided in \cite[Lemmas 2.1 and 2.2]{Kham-Jde}, but with proper modification on the functional spaces and additional minor adjustments.	
	
	For the sake of the reader,  and in order to guarantee we can do state such a result in our case, we present the detailed proof in  Appendix \ref{sec-proof-khan-lemma}.
\end{proof}

\subsection{Lipschitz property w.r.t. index $\lambda$}\label{subsec-lip-cont}

Below, we are going to prove the $\lambda$-Lipschitz property, and therefore continuity, of the following mapping
$$
[0,1]\ni\lambda \ \mapsto \ S_\lambda(t)U_0 \in \mathcal{H},
$$
for all given $t\geq0$ and $U_0\in B$, where $B\subset \mathcal{H}$ is bounded 
set.	
To analyze such a Lipschitz continuity in terms of the parameter $\lambda$, we turn ourselves back to the notation $u^\lambda$ to the solution of \eqref{P}, that is, to the dynamical system $S_\lambda(t)U_0=(u^\lambda(t),u^\lambda_{t}(t))$ given in \eqref{sist-din}.

\begin{proposition}[$\lambda$-Lipschitz Property]\label{prop-esp-conti}
	Under the assumptions of  Theorem 
	$\ref{theo-existence} ,$ 
	let us consider an arbitrary bounded set $B \subset {\mathcal{H}}$  and denote by $S_\lambda(t)U_0=(u^{\lambda}(t),u^{\lambda}_{t}(t))$ the trajectory solution corresponding to initial data  $U_0=(u_{0},u_{1})\in B$.
	Then, for any given index $\lambda_0\in[0,1]$, there exists a positive non-decreasing function 
	$\mathcal{\hat{Q}}(t)=\mathcal{\hat{Q}}\big(B,||h||,t\big)$ such that
	\begin{equation}\label{dependenceVI-eps}
	\|S_\lambda(t)U_0-S_{\lambda_0}(t)U_0\|_{\mathcal{H}} \leq \mathcal{\hat{Q}}(t) |\lambda-\lambda_0| ,\quad t\geq0.
	\end{equation}
\end{proposition}
\begin{proof}
	Let us fix
	$\lambda_0\in[0,1]$ and set $w^{\lambda}:=u^{\lambda}-u^{\lambda_0}$. Then, $w^\lambda$ is a solution (in the weak and strong sense) of the following problem
	\begin{equation}\label{dependence-epsilon}
	\left\{\begin{array}{l}
	\displaystyle{w_{tt}^{\lambda}+\kappa A w^{\lambda}+A_1 w^{\lambda}+\frac{\gamma}{2}\Pi_1^{\lambda}\, w^{\lambda}_t+\frac{\gamma}{2}\Pi_2^{\lambda}\,\left[\,u^{\lambda}_t+u^{\lambda_0}_t\,\right]+F(w^{\lambda})=(\lambda-\lambda_0)h,} \smallskip \\
	\displaystyle{w^{\lambda}(0)=0,\quad w_t^{\lambda}(0)=0,}
	\end{array}\right.
	\end{equation}
	where hereafter we take the advantage of notations and estimates introduced in the proof of Theorem \ref{theo-dependence}, namely, we firstly set  $F(w^{\lambda})=f(u^{\lambda})-f(u^{\lambda_0})$ and 
	$$\Pi_i^{\lambda}(t)=\|S_\lambda(t)U_0\|^{2q}_{\mathcal{H}^{\alpha}}+(-1)^{1-i}\|S_{\lambda_0}(t)U_0\|^{2q}_{\mathcal{H}^{\alpha}},\quad i=1,2.$$
	
	Taking the multiplier $w_t^{\lambda}$ in 
	\eqref{dependence-epsilon}, we have
	\begin{equation}\label{dependenceII-eps}
	\frac{1}{2}\frac{d}{dt}\mathcal{E}^{\lambda}_w(t)+\frac{\gamma}{2}\Pi_1^{\lambda}(t)\|w_t^{\lambda}(t)\|^2=\mathcal{J}_1^{\lambda}(t)+\mathcal{J}_2^{\lambda}(t)+ (\lambda-\lambda_0)(h,w_t^{\lambda}(t)),
	\end{equation}
	where
	\begin{eqnarray*}
		\mathcal{E}^{\lambda}_w(t)& = & \|w_t^{\lambda}(t)\|^2+\kappa\|A^{1/2} w^{\lambda}(t)\|^2+\|A^{1/2}_1 w^{\lambda}(t)\|^2,\\
		\mathcal{J}_1^{\lambda}(t)&=&-\,\left(F(w^{\lambda}(t)),w_t^{\lambda}(t)\right),\\
		\mathcal{J}_2^{\lambda}(t)&=&-\,\frac{\gamma}{2}\Pi_2^{\lambda}(t)\left(u_t^{\lambda}(t)+
		u_t^{\lambda_0}(t),w_t^{\lambda}(t)\right).
	\end{eqnarray*}
	Repeating the same arguments as in \eqref{dependenceIII}, \eqref{J1} and \eqref{J2}, we infer
	\begin{equation}\label{equiv-eps}
	||S_\lambda(t)U_0-S_{\lambda_0}(t)U_0||_{\mathcal{H}}^2\le \mathcal{E}_w(t)\le \mu_\kappa||S_\lambda(t)U_0-S_{\lambda_0}(t)U_0||_{\mathcal{H}}^2, \quad t\geq 0,
	\end{equation}
	and
	$$|\mathcal{J}_1^{\lambda}(t)|,\, | \mathcal{J}_2^{\lambda}(t)|  \leq C_{{B}} \mathcal{E}_w^{\lambda}(t), \quad t\geq 0,$$
	for some constant $C_{{B}}>0$.
	Additionally, 
	using H\"older and Young's inequalities and since $\mathcal{D}(A^{1/2}_1)\hookrightarrow H$, we get
	\begin{equation*}
	|(\lambda-\lambda_0)(h,w^{\lambda}_t(t))|\leq \frac{1}{2}|\lambda-\lambda_0|^2\|h\|^2+\frac{1}{2}\mathcal{E}_w^{\lambda}(t).
	\end{equation*}
	Replacing the latter two estimates in 
	\eqref{dependenceII-eps} we have
	\begin{equation*}
	\frac{d}{dt}\mathcal{E}_w^{\lambda}(t)\le
	C_{B}\mathcal{E}_w^{\lambda}(t)+|\lambda-\lambda_0|^2\|h\|^2,
	\end{equation*}
	for some constant $C_B>0$, and integrating it on $(0,t)$, we arrive at 
	\begin{equation*}
	\mathcal{E}_w^{\lambda}(t)\le 
	C_{{B}}\int_0^t \mathcal{E}_w^{\lambda}(s)\, ds +t|\lambda-\lambda_0|^2\|h\|^2, \quad t>0.
	\end{equation*}
	Therefore, from Gronwall's inequality and \eqref{equiv-eps}, we finally conclude that  \eqref{dependenceVI-eps} holds true with 
	$\mathcal{\hat{Q}}(t):=t^{1/2}e^{C_{{B}}t/2}\|h\|$.
\end{proof}

\section{Long-time dynamics: case \ref{k1}}\label{sec-dynamics}

Before proceeding with the main results of this article, 
we notify for the reader's guidance that all abstract concepts and results  on dynamical systems are reminded in Appendix \ref{sec-appendixB},  by following e.g. the references \cite{babin,carvalho-sonner1,chueshov-solo,chueshov-lasiecka-2005, chueshov-white, chueshov-yellow,eden-etal,  hale,Robinson-Olson-Hoang,Lady,Robinson-book, temam1}.

\subsection{Attractors and continuity: cases 
$\alpha\in[0,1)$ and $p^*\leq\frac{2n}{n-4}$}\label{subsec4.1-k1}

We initially remark that all results up to now hold for any critical parameters $\alpha\in[0,1]$ and $p^*\leq\frac{2n}{n-4}$. However,  due to compactness issues inside this case \ref{k1}, the next theorem dealing with the existence of a family of attractors and its continuity will require the subcritical assumption w.r.t. to fractional powers $\alpha$, namely, $\alpha\in[0,1)$, but we still work in the critical scenario w.r.t. the source growth exponent  $p^*\leq\frac{2n}{n-4}.$

 \begin{theorem}\label{theo_main} Let us take on the same assumptions of Theorem  $\ref{theo-existence}, $ with the additional condition $\alpha\in
 	[0,1)$, and let $(\mathcal{H},S_\lambda(t))$ be the dynamical system given by \eqref{sist-din}. Then. we have:
 
\begin{enumerate}[label={\rm(\textbf{I.\arabic*})}]
	
	\item\label{I1}  {\bf Asymptotic Smoothness/Compactness.} 	
	For every $\lambda \in [0,1]$,
	the dynamical system
	$(\mathcal{H},S_\lambda(t))$  is	asymptotically smooth/compact.
	
		\item\label{I2} {\bf Family of Attractors.}  For every $\lambda \in [0,1]$, the dynamical system $(\mathcal{H},S_\lambda(t))$   possesses a   global attractor $\mathfrak{A}_\lambda\subset \mathcal{H}$,  which is compact and connected.

			\item\label{I3}  {\bf Geometrical Structure.} The family of global attractors $\{\mathfrak{A}_{\lambda}\}_{\lambda \in [0,1]}$  is characterized by the unstable manifold emanating from the set of stationary solutions, namely, we have 
		$
		\mathfrak{A}_\lambda =\mathbf{M}^u(\mathcal{N}_\lambda)   
		$
with
		\begin{equation}\label{state-set}
\mathcal{N}_\lambda=\left\{ (u,0)\in\mathcal{H}; \ \kappa A u +
		A_1 u +f(u)= h_\lambda \right\}, \ \ \lambda\in[0,1].
		\end{equation}

	\item\label{I4} {\bf Equilibria Set.} 	
	Every trajectory stabilizes to the set
$\mathcal{N}_\lambda$ in the sense that 
$$
\lim_{t\to +\infty} \mbox{\rm dist}(S_\lambda(t)U_0,\mathcal{N}_\lambda) = 0, \quad \forall \, U_0\in \mathcal{H}.
$$
In particular, the set $\mathcal{N}_\lambda:=\mathfrak{A}_\lambda^{\min}$ is a global minimal attractor  
for every $\lambda\in[0,1].$ Moreover, 
any trajectory from $\mathfrak{A}_\lambda$ has an upper bound   in terms of  $\mathfrak{A}_\lambda^{\min}$, namely, 
\begin{equation}\label{upper-bound}
\sup \big\{ || (u,u_t)||_{\mathcal{H}}; \ (u,u_t) \in \mathfrak{A}_\lambda\big\} \leq \sup \big\{||(u, 0)||_{\mathcal{H}}; \ (u, 0) \in \mathfrak{A}_\lambda^{\min}\big\}.
\end{equation}

	\item\label{I5} {\bf Non-triviality.} 	 The family of  minimal attractors
	$\{\mathfrak{A}_{\lambda}^{\min}\}_{\lambda \in [0,1]}$
	 is nontrivial. In other words, even if $h\equiv0$, the minimal attractor    $\mathfrak{A}_\lambda^{\min}$ has at least two stationary solutions   for every $\lambda\in[0,1].$

   \item\label{I6} {\bf Triviality.}  
If we additionally suppose that $h\equiv 0$    and $C_f=0$ in \eqref{assumption_f},  the attractor  $\mathfrak{A}_0$ is trivial.  More precisely, $\mathfrak{A}_0=\{(0,0)\}$ with polynomial  attraction depending  on the exponent $ q\geq \frac{1}{2}$ as follows
\begin{equation}\label{exp-atract}
\mbox{\rm dist}_{\mathcal{H}}(S_\lambda(t)B,\mathfrak{A}_0) = \sup_{U_0\in B} 
\|S_\lambda(t)U_0\|_{\mathcal{H}} \leq \frac{C_B}{\Big[ c_{B}+q\hat{c}_{B} \,  t\Big]^{1/2q}}, \quad t\to\infty,
\end{equation}
for any initial data $U_0$ lying in bounded sets $B\subset \mathcal{H}$, where  $C_B, \,c_{B},\, \hat{c}_B>0$ are constants depending on $B$.

	\item\label{I7} {\bf Upper Semicontinuity.}   
	The family of global attractors $\{\mathfrak{A}_{\lambda}\}_{\lambda \in [0,1]}$   is upper semicontinuous at any fixed
	$\lambda_0 \in [0,1]$, that is, 
	$$\lim_{\lambda \to \lambda_0}{\rm dist}_{\mathcal{H}}\left(\mathfrak{A}_{\lambda},\mathfrak{A}_{\lambda_0}\right)=0.
	$$			
	
	\item\label{I8} {\bf Residual Continuity.} The family of global attractors $\{\mathfrak{A}_{\lambda}\}_{\lambda \in [0,1]}$ is  continuous in a residual\footnote{Let $X$  be a complete metric space and $Y\subset X$. We recall  that  $Y$ is {\it residual} in $X$ if  $X \backslash Y$ is a countable union of nowhere dense sets.}  set $I \subset[0,1]$, that is, for any $\lambda_0 \in I$,
	$$\lim_{\lambda \to \lambda_0}\left[{\rm dist}_{\mathcal{H}}\left(\mathfrak{A}_{\lambda_0},\mathfrak{A}_{\lambda}\right)+{\rm dist}_{\mathcal{H}}\left(\mathfrak{A}_{\lambda},\mathfrak{A}_{\lambda_0}\right)\right]=0.
	$$
	In particular, the set of
	continuity points of $\mathfrak{A}_{\lambda}$ is dense in  $[0,1]$.		
	
\end{enumerate}
\end{theorem}

\begin{proof}
To the proof, we gather the 	
 ingredients	coming from
Section \ref{sec-tech-results} 
along with the abstract results reminded in Appendix \ref{sec-appendixB}. 
 	
\smallskip 	
\noindent \ref{I1}  
Let us initially consider a  bounded positively invariant set $B\subset\mathcal{H}$,  take two trajectory solutions $S_\lambda(t)U^i=(u^i(t),u^i_{t}(t)), \, i=1,2,$ corresponding to initial data $U^i=(u^i_{0},u^i_{1})\in B, \, i=1,2, $ and consider any  $\varepsilon>0$. 

From the key inequality \eqref{key-inequality}, there
exists a time large enough $T:=T_{\lambda,B}>0$ such that
\begin{eqnarray}\label{tt3}
\|S_\lambda(T) U^1-S_\lambda(T) U^2\|_{\mathcal{H}}  & \le &  \varepsilon+\psi_{\varepsilon,B,T}(U^1,U^2),
\end{eqnarray}
where we set $\psi_{\varepsilon,B,T}:\mathcal{H}\times
\mathcal{H} \to \mathbb{R}$ by
\begin{align}\label{tt2}
\psi_{\varepsilon,B,T}(U^1,U^2)  :=  \, & C_{B}\sup_{0\le s\le T+1}\left[\|A^{\alpha}u^1(s)-A^{\alpha}u^2(s)\|^{\frac{2(q+1)}{2q+1}}+\|u^1(s)-u^2(s)\|^{2}_{p+2}\right]^{\frac{1}{2(q+1)}}\nonumber\\
&+\left[J_f(T)\right]^{\frac{1}{2(q+1)}}
\end{align}
for some constant $C_B>0$ and  $J_f(t)$ given in \eqref{critical}.

Now, given a sequence of initial data $U^n=(u_0^{n},u_1^{n})\in B$, as before, we write $S_\lambda(t)U^n = (u^{n}(t),u_{t}^{n}(t))$.
Since $B$ is invariant by $S_\lambda(t)$, $t\ge 0$, it follows that
$(u^{n}(t),u_{t}^{n}(t))$ are uniformly bounded in
$\mathcal{H}=\mathcal{D}(A^{1/2}_1) \times H$ from Corollary \ref{cor-bdd}. Thus,
$$
(u^{n},u^{n}_{t}) \; \mbox{is bounded in} \; C([0,T+1], \mathcal{H}).
$$

Below is the precise moment we invoke the assumption $\alpha\in[0,1)$.
Indeed, for such fractional powers and any $p^*\leq\frac{2n}{n-4}$,
we use the fact that $\mathcal{D}(A^{1/2}_1)\hookrightarrow\mathcal{D}(A^{\alpha})$ and $\mathcal{D}(A^{1/2}_1)\hookrightarrow L^{p+2}(\Omega)$ are compact  embeddings. Then, by virtue of
 \cite[Corollary 4]{Simon}  there exists a subsequence, still denoted by
$(u^{n})$, such that
\begin{equation}\label{strong}
(u^{n}) \ \ \mbox{converges strongly in} \ \ C([0,T+1], \mathcal{D}(A^{\alpha})\cap L^{p+2}(\Omega)).
\end{equation}
Therefore,
\begin{equation}\label{last2}
\lim_{n\rightarrow \infty}\lim_{m\rightarrow\infty}\sup_{0\le s\le T+1}\left[\|A^{\alpha}u^n(s)-A^{\alpha}u^m(s)\|^{\frac{2(q+1)}{2q+1}}+\|u^n(s)-u^m(s)\|^{2}_{p+2}\right]^{\frac{1}{2(q+1)}} =0.
\end{equation}
Additionally, from \eqref{strong} and \eqref{Khan-limit}, and the expression for  $J_f(t)$ in \eqref{critical},  we also have
\begin{equation}\label{cong-J}
 \lim_{n\rightarrow \infty}\lim_{m\rightarrow\infty} J_f(s)=\lim_{n\rightarrow \infty}\lim_{m\rightarrow\infty} J_f(u^n(s)-u^m(s),u_t^n(s)-u_t^m(s))=0.
\end{equation}
Then, from \eqref{last2} and \eqref{cong-J},  we conclude 
\begin{equation}\label{contrative}
\lim_{n\rightarrow \infty}\lim_{m\rightarrow\infty}\psi_{\varepsilon,B,T}(U^n,U^m) =0,
\end{equation}
for every sequence $U^n\in B$, which proves that $\psi_{\lambda,B,T}$ is a contractive function on $B\times B$. 

Therefore,  from \eqref{tt3} and \eqref{contrative},  we can apply
Theorem \ref{theo_khanma}  to conclude that $(\mathcal{H},S_\lambda(t))$ is  asymptotically smooth. It is also  asymptotically compact in view of  Proposition \ref{prop-equiv}.

\smallskip 
\noindent \ref{I2} It follows from Corollary \ref{dissipative} and step \ref{I1}, in combination  with 
Theorem \ref{theo-attract}.

\smallskip 
\noindent \ref{I3} It follows from  Corollary \ref{cor-gradient} along with Theorem \ref{theo_geometric_struc}.

\smallskip 
\noindent \ref{I4}  It follows from \ref{I3}, Corollary \ref{cor-gradient}, and Theorem \ref{theo_minimal_attra}. Moreover, the upper bound \eqref{upper-bound} follows from the first limit in 
Theorem \ref{theo_geometric_struc},
Corollary  \ref{cor-gradient}, and Lemma \ref{lem4-2},
once we have that the  Lyapunov function $\Phi:=E$ is topologically equivalent to the norm of the phase
space $\mathcal{H}$ (see also \cite[Remark 7.5.8]{chueshov-yellow}).

\smallskip 
\noindent \ref{I5} It follows directly from 
Lemma \ref{lem-nontrivial}.

\smallskip 
\noindent \ref{I6} It follows directly from  the
second part of Lemma \ref{lem4-2}, by applying Corollary $ \ref{P3} $ for initial data lying in bounded sets $B\subset \mathcal{H}$. 

\smallskip 
\noindent \ref{I7}
From  Corollary \ref{dissipative},  
 the dynamical system
$(\mathcal{H},S_\lambda(t))$  has a bounded absorbing set $\mathcal{B}\subset \mathcal{H}$  uniformly bounded w.r.t. $\lambda\in[0,1]$. Thus, the
 attractors   $\mathfrak{A}_{\lambda} \subset \mathcal{B}$ are uniformly bounded for all $\lambda \in [0,1]$. 
Additionally, from Proposition
\ref{prop-esp-conti}, we have
\begin{equation}\label{cont}
\lim_{\lambda \to \lambda_0}\sup_{U_0 \in  \mathcal{B}}\|S_\lambda(t)U_0-S_{\lambda_0}(t)U_0\|_{\mathcal{H}}=0, \quad t\geq t_0,
\end{equation}
for every given $t_0\geq0$. Therefore, employing Theorem \ref{conti-Robinson}, we conclude that  the family $\{\mathfrak{A}_{\lambda}\}_{\lambda \in [0,1]}$ is upper semicontinuous  
at any fixed
$\lambda_0 \in [0,1].$

\smallskip 
\noindent \ref{I8}   
By noting that \eqref{cont} is valid for all $t_0>0$ and any bounded set $B\subset\mathcal{H}$ (see again Proposition
\ref{prop-esp-conti}), then we are under the assumptions of Theorem \ref{theo-ROH}. Therefore, the desired conclusion on residual continuity follows.
\end{proof}

\subsection{Kolmogorov
	$\varepsilon$-Entropy: cases 
	$\alpha\in[0,1)$ and $p^*<\frac{2n}{n-4}$}

As observed in Remark \ref{rem-finite}, we are going to appeal to Corollary $\ref{cor-key-subcritic}$ in order to reach an estimate for the 
 Kolmogorov's
$\varepsilon$-entropy of the 
attractor $\mathfrak{A}_{\lambda}$ corresponding to 
the dynamical system $(\mathcal{H},S_\lambda(t))$. For this reason (see the ``weaker'' stabilizability inequality \eqref{key-inequality-cor}), we must employ the subcritical case 
 w.r.t. parameters $\alpha$ and $p^*$.
Our second main result reads as follows.

\begin{theorem}[{\bf Kolmogorov
		$\varepsilon$-Entropy}]\label{theo_mainIII}  Let us consider the assumptions of Theorem  $\ref{theo_main}$, with the additional condition $p<\frac{4}{n-4}$  in Assumption $\ref{assumption}$. Then, there  exists $0<\varepsilon_{0}<1$ such that for all $\varepsilon \leq \varepsilon_{0}<1$, the 
	Kolmogorov $\varepsilon-$entropy $H_{\varepsilon}(\mathfrak{A}_{\lambda})$ of the existing global attractor $\mathfrak{A}_{\lambda}, \, \lambda \in [0,1],$ satisfies the following estimate for  arbitrary $\delta \in(0,1)$
	\begin{equation}\label{eps-entrop1}
	H_{\varepsilon}(\mathfrak{A}_{\lambda}) \leq \frac{2}{1-\delta}\int_{\varepsilon}^{\varepsilon_{0}} \frac{\ln m\left(g_{\delta}^{-1}(s), z(s)\right)}{s} d s+H_{g_{\delta}\left(\varepsilon_{0}\right)}(\mathfrak{A}_{\lambda}),
	\end{equation}
	where  $g_{\delta}(s)=\frac{1+\delta}{2} s$ and $z(s)=\frac{1}{2} (\delta s)^{2(q+1)}, \,  0<s<\varepsilon_{0},$ and    
	\begin{equation*}
	m(r,a)=\sup \{m(B, a); \  B \subseteq \mathfrak{A}_{\lambda}, \, \operatorname{diam} B \leq 2r\},
	\end{equation*}
	with $m(B,a)$ being the maximal number  of elements $u_{j}^{B} \in B$
	such that, for any $a>0$, we have
	$
	\varrho\left(S_\lambda(T^{\star})u_{j}^{B},S_\lambda(T^{\star})u_{i}^{B}\right)>a, i \neq j, i, j=1, \ldots, m(B,a),
	$ for $T^{\star}>0$ large enough, and $\varrho$ is a  
	pseudometric on $\mathcal{H}$.	
\end{theorem} 
\begin{proof}
Let us  consider two trajectory solutions $S_\lambda(t)U^i=(u^i(t),u^i_{t}(t)), \, i=1,2,$ corresponding to initial data $U^i=(u^i_{0},u^i_{1})\in \mathfrak{A}_{\lambda}, \, i=1,2, \, \lambda\in[0,1],$ and still denote $w:=u^1-u^2.$
Since $\mathfrak{A}_{\lambda}$ is compact (and then bounded) and invariant $S_\lambda(t)\mathfrak{A}_{\lambda}=\mathfrak{A}_{\lambda}$, then $S_\lambda(t)U^i\in\mathfrak{A}_{\lambda}$ for all $t>0$. Additionally,
from Corollary $\ref{cor-key-subcritic}$, there
exists a time   $T^{\star}:=T^\star({\mathfrak{A}_{\lambda}})>0$ such that
 \begin{eqnarray}\label{kolmog-1}
\|S_\lambda(T^\star) U^1-S_\lambda(T^\star) U^2\|_{\mathcal{H}}  & \le &  \frac{1}{2}\|U^1-U^2\|_{\mathcal{H}} \\ && +\left[ C_{\mathfrak{A}_{\lambda}}\sup_{0\le s\le T^\star+1}\left(\|A^{\alpha}w(s)\|+\|w(s)\|_{p^*}\right)\right]^{\frac{1}{2(q+1)}},
\nonumber
\end{eqnarray}
for some constant $C_{\mathfrak{A}_{\lambda}}>0$ and all  $U^1,U^2\in \mathfrak{A}_{\lambda}$, where (by virtue of  Corollary \ref{cor-bdd} and the identity $\frac{2(q+1)}{2q+1}=1+\frac{1}{2q+1}$) we have used
\begin{equation}\label{estimate}
\|A^{\alpha}w(s)\|^{\frac{2(q+1)}{2q+1}}+||w(s)||_{p^*}^{\frac{2(q+1)}{2q+1}} \leq C_{\mathfrak{A}_{\lambda}}  \left(\|A^{\alpha}w(s)\|+||w(s)||_{p^*}\right), \ \ s\geq0.
\end{equation}
Moreover,  from  Theorem \ref{theo-dependence}, there exists a constant $L^{\star}>0$ such that
\begin{equation}\label{kolmog-2}
||S_\lambda(T^\star) U^1-S_\lambda(T^\star) U^2||_{\mathcal{H}}  \, \le  \,  L^{\star}
\|U^1-U^2\|_{\mathcal{H}},\quad  \forall \, U^1,U^2\in \mathfrak{A}_{\lambda}.
\end{equation}

Therefore, one can see from \eqref{kolmog-1}-\eqref{kolmog-2} that  assumptions 1 to 4 of Theorem \ref{theo-kolmog-entrop} are fulfilled with
\begin{equation}\label{notations}
\begin{array}{l}
 \mathscr{M}:=\mathfrak{A}_{\lambda}, \ \ V:=S_\lambda(T^\star), \ \ g(s)=\frac{1}{2}s, \ \  h(s):=s^{\frac{1}{2(q+1)}}, \ \ \varrho_1\equiv0, \medskip \\
 \displaystyle \varrho_2(S_\lambda(T^\star) U^1,S_\lambda(T^\star) U^2):=C_{\mathfrak{A}_{\lambda}}\sup_{0\le s\le T^\star+1}\left(\|A^{\alpha}w(s)\|+\|w(s)\|_{p^*}\right).
\end{array}\end{equation}
It is worth mentioning that due to the compactness embedding $\mathcal{D}(A^{1/2}_1)\hookrightarrow\mathcal{D}(A^{\alpha}) \cap L^{p^*},$ then $\varrho:=\varrho_2$ is a compact seminorm\footnote{We recall that a seminorm $n_X(\cdot)$ defined on a Banach space
	$X$ is {\it compact} if whenever a sequence $x_j \to 0$
	weakly in $X$ one has $n_X(x_{j}) \to 0$.} on $\mathcal{D}(A^{1/2}_1),$ and item 4 of Theorem \ref{theo-kolmog-entrop} follows by a standard argument, see e.g.   \cite[p. 55]{chueshov-lasiecka-2005} for a generic proof of this fact.

Hence, the estimate for the Kolmogorov $\varepsilon$-entropy $H_{\varepsilon}(\mathfrak{A}_{\lambda})$ provided in \eqref{eps-entrop1} follows from the conclusion of Theorem \ref{theo-kolmog-entrop}. 
\end{proof}

\begin{remark}
	As final information in this subsection, we note that both  Theorems $ \ref{theo_main} $ and $ \ref{theo_mainIII} $  remain unchanged when we neglect the potential energy in the 
	damping coefficient $\mathcal{E}_{\alpha}(u,u_t)$ given in \eqref{Ealpha}. In this case  the equation in \eqref{P} reduces to the particular  one with nonlinear averaged damping
\begin{equation*}
u_{tt}+\kappa A u + A_1 u+\gamma \|u_t\|^{2q} u_t+f(u)= h_\lambda ,\quad t> 0, \ \lambda\in[0,1]. 
	\end{equation*} 
Therefore, all results in the present section provide a generalization of the ones in \cite{zhao1,zhao2} when one considers a constant coefficient of extensibility $\kappa$. Additionally, the same happens if one takes a slightly more general situation with respect to $\kappa(\cdot)$ as a nonlocal function under suitable properties like in \cite{jorge-narciso-DCDS,jorge-narciso-EECT,zhao1,zhao2}.	
\end{remark}

\subsection{Attempts for finite dimensionality}
\label{sec-final-comments}

Below we try to clarify  how hard (if not impossible) is to achieve the finiteness of the fractal dimension ($
\mbox{dim}^{f}_{\mathcal{H}}\mathfrak{A}_{\lambda}$) of the attractors  $\mathfrak{A}_\lambda,\, \lambda\in[0,1],$  corresponding to the   dynamical system $(\mathcal{H},S_\lambda(t))$ defined in  \eqref{sist-din} in case \ref{k1}. 
A first attempt is trying to invoke Theorem \ref{theo-dim-finite}. In this way, in view of notations in \eqref{notations}, we can rewrite \eqref{kolmog-1} as 
 \begin{equation}\label{finite-1}
\|V U^1-V U^2\|_{\mathcal{H}}    \le   \frac{1}{2}\|U^1-U^2\|_{\mathcal{H}}   +  \left[\varrho(V U^1,VU^2) \right]^{\frac{1}{2(q+1)}}.
\end{equation}
Therefore, all hypotheses of   Theorem \ref{theo-dim-finite} are satisfied except for item $(ii)$ that requires the linearity of function $h(s)=s_0 s$. Indeed, this fact is impossible in our case since $h(s)=s^{\frac{1}{2(q+1)}}$ for $ q\geq \frac{1}{2} $. The only chance  to achieve \eqref{finite-1} with linear function $h(s)$ is to consider the particular (and already known) case $q=0$. But  this latter especial case  reduces the damping term in \eqref{P} to the linear one $\gamma  u_t, \, \gamma>0$, which in turn  
is a particular case of \ref{k2} in what concerns function $k(\cdot)$. For such a case, we show  later that $
\mbox{dim}^{f}_{\mathcal{H}}\mathfrak{A}_{\lambda}<\infty$, the regularity of any trajectory from the attractors and (generalized) fractal exponential attractors.

Going back to this ``worse'' case \ref{k1}, another attempt in trying to achieve  
\eqref{finite-1} (or else \eqref{kolmog-1} which comes from Corollary $\ref{cor-key-subcritic}$) with a proper power concerning its last term is to regard  perturbed energy computations instead of Nakao's method as in the  proof of Proposition 
$ \ref{theo_inequality} $ and Corollary \ref{cor-key-subcritic}. In such a way, the following result can be proved. 
\begin{proposition} \label{prop-atemp}
	Under the same assumptions and statements of Proposition 
	$ \ref{theo_inequality} $, and given any $\epsilon>0$,   there exist 
 constants $c_B, C_B>0$ depending on $B$  such that
	\begin{align}\label{key-inequality-atemp}
	||S_\lambda(t)U^1-S_\lambda(t)U^2||_{\mathcal{H}}^2\le &\;
	 C_B \left\|U^{1}-U^{2}\right\|_{\mathcal{H}}^2 e^{- c_B t} + \epsilon^2
\\
	&\;+C_{B} \int_{0}^{t}e^{- c_B(t-s)}\Big(\|A^{\alpha}w(s)\|^{\frac{2(q+1)}{2q+1}}+||w(s)||_{p^*}^{\frac{2(q+1)}{2q+1}}\Big)ds,	\nonumber
	\end{align}
for every $t>0,$ where we still denote $w=u^{1}-u^{2}$.
\end{proposition}
\begin{proof}
	The proof relies on  energy perturbation and similar technical estimates  as used in the proof of Proposition 
 $ \ref{theo_inequality} $ and Corollary \ref{cor-key-subcritic}, along with proper Young's inequality. Thus, it will be omitted.
\end{proof}

Hence, as a consequence of \eqref{key-inequality-atemp} on $B:=\mathfrak{A}_\lambda$,  and using again \eqref{estimate}, we arrive at 
 \begin{equation}\label{finite-2}
\|S_\lambda(T^\star) U^1-S_\lambda(T^\star) U^2\|_{\mathcal{H}}    \le   \frac{1}{2}\|U^1-U^2\|_{\mathcal{H}}  + \epsilon +  \left[\varrho(S_\lambda(T^\star) U^1,S_\lambda(T^\star)U^2)\right]^{\frac{1}{2}},
\end{equation}
for some  time  $T^{\star}:=T^\star({\mathfrak{A}_{\lambda}})>0$ large enough, some compact seminorm $\varrho$, and any $\epsilon>0$. Nonetheless, \eqref{finite-2} is not enough
to achieve the finiteness of  the
fractal dimension $
\mbox{dim}^{f}_{\mathcal{H}}\mathfrak{A}_{\lambda}$ by means of 
Theorem \ref{theo-dim-finite}.

The above approaches \eqref{finite-1} and \eqref{finite-2} for 
 studying the dimensionality of  the attractors $\mathfrak{A}_{\lambda}$  raise similar issues 
 as presented in \cite{chueshov-lasiecka-2005} in terms of (more) general stabilizability estimates. 
  Indeed, this is the exact moment where we explore the difficulty imposed by function $k(s)= \gamma s^q$ for any $q\geq \frac{1}{2}$ in case \ref{k1} because under this structure    the character of the functions $h(s)=s^{\frac{1}{2(q+1)}}$ and $h(s)=s^{\frac{1}{2}}$ present in the  lower order terms (LOT) 
 $$
 \mbox{LOT}\left(U^{1}, U^{2}\right):=\varrho(S_\lambda(T^\star) U^1,S_\lambda(T^\star)U^2)
 $$
are determined from the behavior of the nonlinear damping term
\begin{equation}\label{damping-term}
k\big(\mathcal{E}_{\alpha}(u,u_t)\big)u_t=\gamma\big[\mathcal{E}_{\alpha}(u,u_t)\big]^q u_t, \ \ q\geq \frac{1}{2},
\end{equation}
with $\mathcal{E}_{\alpha}(u,u_t)$ being set in \eqref{Ealpha}.
Moreover, due to the computations  in the proof of Proposition 
$ \ref{theo_inequality} $  (see \eqref{key-1}-\eqref{key-2})
it seems that the nonlocal nonlinear damping term \eqref{damping-term}
does not even provide a suitable coercivity property as usual for nonlinear damping like $D(u_t)$, where $D:\mathbb{R}\to\mathbb{R}$  is a real function with growth exponent $q$,  namely, 
\begin{equation}\label{ass-impl}
(D(s)-D(r))(s-r) \geq {c_{q}}|s-r|^{q+2}, \quad \forall \ s,r\in\mathbb{R}, \ q\geq0,
\end{equation}
or 
\begin{equation}\label{ass-vando}
(D(s)-D(r))(s-r) \geq {c_{q}}\left(|s|^q+|r|^q\right)|s-r|^{2}, \quad \forall \ s,r\in\mathbb{R}, \ q\geq0,
\end{equation}
or else,
for any given $\epsilon>0,$ there exists a constant
$C_{\epsilon}>0$ such that
\begin{equation}\label{ass-dom}
C_{\epsilon}(D(s)-D(r))(s-r) \geq|s-r|^{2}-\epsilon, \quad \forall \, s, r \in \mathbb{R}.
\end{equation}

It is worth mentioning that the above 
assumptions \eqref{ass-impl}-\eqref{ass-dom} have been extensively regarded in the literature in what concerns long-time dynamics of  hyperbolic-type
second-order evolution problems  with nonlinear damping, see for instance the works
\cite{Aloiui-Ben-Haraux,chueshov-strong,chueshov-eller-lasiecka1,chueshov-eller-lasiecka2,chueshov-lasiecka-2004-JDDE,chueshov-white,chueshov-yellow,chueshov-lasiecka-tound1,chueshov-lasiecka-tound2,feireisl,jorge-narciso-EECT,Kham-2010,lasiecka-ruzmainkina,Nakao,prazak,sun} where (at least) some of them do play an important role in finding the asymptotic smoothness  of  the corresponding nonlinear infinite-dimensional dynamical system.

Although conditions \eqref{ass-impl}-\eqref{ass-dom} have shown to be very effective in building the  existence of compact global attractors, the picture is much more delicate when  one deals with respect to regularity and especially 
finite-dimensionality of such global attractors. Indeed, among the above-mentioned works with nonlinear damping, those that proved the finiteness of fractal (or Hausdorff) dimension  
strongly used more hypotheses on $D$ and its derivatives $D', \, D''$, see for instance \cite[Theorem 1.4]{chueshov-eller-lasiecka2},  \cite[Theorem 1.5]{chueshov-lasiecka-2004-JDDE}, \cite[Theorem 5.8]{chueshov-white}, \cite[Theorem 9.2.6]{chueshov-yellow}, \cite[Assumption 3 and Theorem 3.5]{chueshov-lasiecka-tound1}, 
\cite[Proposition 4 and Theorem 3.2]{jorge-narciso-EECT}, \cite[Theorems 3.1 and 3.2]{Kham-Jde}, \cite[Theorem 8]{lasiecka-ruzmainkina}, and \cite[Theorem 1.1]{prazak}, just to quote a few.
Nonetheless, such damping  controlling by means of its derivative does not seem to be applicable to the nonlinear damping \eqref{damping-term} due to its nonlocal structure. This is similar to the nonlocal damping addressed in \cite{Aloiui-Ben-Haraux}, where the author involves differential operators  in space and  covers a wide class of average nonlocal damping. Therefore, we conclude that \eqref{damping-term} represents a generalization of the linear damping in a different way of  the existing literature w.r.t. the nonlinear damping terms $D(u_t)$.

In conclusion, to our best knowledge, there is no keen theory to conclude the finiteness of the fractal (or at least Hausdorff) dimension for problems with non-local damping like \eqref{damping-term} where the behavior of lower order terms are determined, in general, by the possibly degenerate function $k(s)=\gamma s^q$ over the linear energy coefficient $\mathcal{E}_{\alpha}(u,u_t)$. 
A way of circumvent this situation is when $k$ is bounded from below but such a case is covered by the case 
\ref{k2} to be analyzed next.

\section{Long-time dynamics: case \ref{k2}}\label{sec-dynamics-k2}

As we are going to show below, this case is much more touchable in the sense that it provides the finiteness and regularity of the global attractor $\mathfrak{A}_\lambda\subset \mathcal{H}, \lambda\in[0,1],$ related to 
the dynamical system $(\mathcal{H},S_\lambda(t))$ given by \eqref{sist-din} in case \ref{k2}. Moreover, in this case, we also prove the existence of (generalized) fractal exponential attractors $\mathfrak{A}^{\exp}_\lambda\subset \mathcal{H}$. Such statements are  due to the positiveness of function $k(\cdot)$ in case \ref{k2}, which allows us to control the behavior of the nonlocal damping (especially at the origin) and this makes this case ``smoother'' than the previous  one.

We also note that due to the structure of damping in case  \ref{k2}, we will be able to work with critical parameters  $\alpha\in[0,1]$ and $p^*\leq\frac{2n}{n-4}$ for obtaining technical estimates  and the computations rely on similar techniques as previously used in the literature, see e.g. \cite{cavalcanti-mato-luci,chueshov-white,gatti-miranv-pata-zelik,jorge-narciso-DCDS,jorge-narciso-EECT}.

In this case, the essential technical results proved in Section  \ref{sec-tech-results} can be refined as follows.

\begin{proposition}[Dissipativity]\label{Prop-exp-decay-k2} 
	Under the assumptions of Theorem  $\ref{theo-w-p-general} $ with 
 $k(\cdot)$ given in case  {\rm \ref{k2}}, there exist positive constants $c=c_{\widetilde{E}(0)},$  $C=C_{\widetilde{E}(0)}$ (which may depend on initial data)
 such that  $\widetilde{E}(t)$ given in \eqref{modified-energy} satisfies
 \begin{equation}\label{ineq-exp-k2}
\widetilde{E}(t)\le C\widetilde{E}(0) e^{-c \,t}+8K_\lambda, \quad t>0.
\end{equation} 
In particular, 
 	the dynamical system
 $(\mathcal{H},S_\lambda(t))$  given by \eqref{sist-din}  is dissipative, say with (positively invariant) 
bounded absorbing set $\mathcal{B}\subset \mathcal{H}$, which is uniformly bounded w.r.t. $\lambda\in[0,1]$.
\end{proposition}
\begin{proof}
	Since $k(\cdot)$ is a  $C^1$-function on $[0,\infty)$ such that $k(s)>0, \, s\geq0$, then the proof follows exactly the same lines as in \cite[Proposition 1]{jorge-narciso-EECT}. 
\end{proof}

Additionally, Corollaries \ref{cor-bdd} to  \ref{P3} can  be adapted to this case as well. Also, and much more important, the following stabiliability inequality can be reached.

\begin{proposition}[Stabilizability Estimate]\label{prop-stabi-ineq-k2}
Under the assumptions of Theorem  $\ref{theo-w-p-general} $ with 
$k(\cdot)$ given in case  {\rm \ref{k2}}, let us  consider a bounded set $B \subset {\mathcal{H}}$ with initial data $U^i=(u^i_{0},u^i_{1})\in B, \, i=1,2.$
Still denoting by $S_\lambda(t)U^i=(u^i(t),u^i_{t}(t)), i=1,2,$ the 
corresponding  dynamical system \eqref{sist-din}, then there exist constants $c_B, C_B>0$ (depending on $B$) such that
	\begin{align}\label{estabil-ineq-k2}
	||S_\lambda(t)U^1-S_\lambda(t)U^2||_{\mathcal{H}}^2\le &\; C_B e^{-c_B t } ||U^1-U^2||_{\mathcal{H}}^2 \nonumber\\
	&\;+C_{B}\int_{0}^{t}e^{-c_B (t-s)} \Big(\|A^{\alpha}w(s)\|^{2}+||w(s)||_{p+2}^{2}\big)ds,
	\end{align}
	for all $t>0$ and $\lambda\in[0,1]$, where  $w=u^1-u^2$.
\end{proposition}
\begin{proof}
	The proof is analogous to \cite[Proposition 1]{jorge-narciso-DCDS} with help of similar arguments as in \cite[ Lemma 4.9]{cavalcanti-mato-luci} to handle critical exponent. 
\end{proof}

The motivation for getting \eqref{estabil-ineq-k2} with respect to critical growth exponent $p$ came from the results in \cite[Proposition 4.13]{chueshov-white} and 
\cite[Lemma 7.1]{gatti-miranv-pata-zelik}. Finally, we also note that Proposition \ref{prop-esp-conti} can  be proved in this case in a very similar way. Therefore, we are able to state our main results in the present section as follows. 

\subsection{Attractors, continuity, finite dimensionality, and regularity: cases $\alpha\in[0,1)$ and $p^*\leq\frac{2n}{n-4}$}

As in Subsection \ref{subsec4.1-k1}, 
the above technical estimates in the present case 
hold for any critical parameters $\alpha\in[0,1]$ and $p^*\leq\frac{2n}{n-4}$, but due to the compactness issues involving the parameter $\alpha$  (which comes now from \eqref{estabil-ineq-k2}), we must work in the subcritical case w.r.t. it.
Nonetheless, nothing changes w.r.t.  
the source growth exponent  $p^*\leq\frac{2n}{n-4}.$

\begin{theorem}\label{theo_main-k2} Let us take on the same assumptions of Theorem  $\ref{theo-w-p-general} $ with $k(\cdot)$ given in case  {\rm \ref{k2}}. We additionally assume that $\alpha\in
	[0,1)$. Then, the dynamical system 	$(\mathcal{H},S_\lambda(t))$ set in \eqref{sist-din} has the following properties:
	
	\begin{enumerate}[label={\rm(\textbf{J.\arabic*})}]
		
		\item\label{J1k2}  {\bf 
		Quasi-Stability.} 	
		For every $\lambda \in [0,1]$,
		the dynamical system
		$(\mathcal{H},S_\lambda(t))$  is asymptotically quasi-stable on any positively invariant bounded set $B\subset\mathcal{H}$. In particular, it is 		asymptotically smooth/compact.

		\item\label{J2k2} {\bf Family of Attractors.}  For every $\lambda \in [0,1]$, the dynamical system $(\mathcal{H},S_\lambda(t))$   possesses a   global attractor $\mathfrak{A}_\lambda\subset \mathcal{H}$,  which is compact and connected.

		\item\label{J3}  {\bf Geometrical Structure.} The family of global attractors $\{\mathfrak{A}_{\lambda}\}_{\lambda \in [0,1]}$  is characterized by the unstable manifold emanating from the set of stationary solutions, namely, we have 
		$
		\mathfrak{A}_\lambda =\mathbf{M}^u(\mathcal{N}_\lambda)   
		$
		with
		\begin{equation*}
		\mathcal{N}_\lambda=\left\{ (u,0)\in\mathcal{H}; \ \kappa A u +
		A_1 u +f(u)= h_\lambda \right\}, \ \ \lambda\in[0,1].
		\end{equation*}

		\item\label{J4} {\bf Equilibria Set.} 	
		Every trajectory stabilizes to the set
		$\mathcal{N}_\lambda$ in the sense that 
		$$
		\lim_{t\to +\infty} \mbox{\rm dist}(S_\lambda(t)U_0,\mathcal{N}_\lambda) = 0, \quad \forall \, U_0\in \mathcal{H}.
		$$
		In particular, the set $\mathcal{N}_\lambda:=\mathfrak{A}_\lambda^{\min}$ is a global minimal attractor  
		for every $\lambda\in[0,1].$ Moreover, 
		any trajectory from $\mathfrak{A}_\lambda$ has an upper bound   in terms of  $\mathfrak{A}_\lambda^{\min}$, namely, 
		\begin{equation*}
		\sup \big\{ || (u,u_t)||_{\mathcal{H}}; \ (u,u_t) \in \mathfrak{A}_\lambda\big\} \leq \sup \big\{||(u, 0)||_{\mathcal{H}}; \ (u, 0) \in \mathfrak{A}_\lambda^{\min}\big\}.
		\end{equation*}

		\item\label{J5} {\bf Non-triviality.} 	 The family of  minimal attractors
		$\{\mathfrak{A}_{\lambda}^{\min}\}_{\lambda \in [0,1]}$
		is nontrivial. In other words, even if $h\equiv0$, the minimal attractor    $\mathfrak{A}_\lambda^{\min}$ has at least two stationary solutions   for every $\lambda\in[0,1].$

		\item\label{J6} {\bf Triviality.}  
		If we additionally suppose that $h\equiv 0$    and $C_f=0$ in \eqref{assumption_f},  the attractor 
		$\mathfrak{A}_0=\{(0,0)\}$
		is trivial with exponential attraction as follows
		\begin{equation*}
		\mbox{\rm dist}_{\mathcal{H}}(S_\lambda(t)B,\mathfrak{A}_0) = \sup_{U_0\in B} 
		\|S_\lambda(t)U_0\|_{\mathcal{H}} \leq {C_B} e^{-c_B t}, \quad t\to\infty,
		\end{equation*}
		for any initial data $U_0\in B\subset \mathcal{H}$, where  $C_B, \,c_B>0$ are constants depending on $B$.

		\item\label{J7} {\bf Upper Semicontinuity.}   
		The family of global attractors $\{\mathfrak{A}_{\lambda}\}_{\lambda \in [0,1]}$   is upper semicontinuous at any fixed
		$\lambda_0 \in [0,1]$, that is, 
		$$\lim_{\lambda \to \lambda_0}{\rm dist}_{\mathcal{H}}\left(\mathfrak{A}_{\lambda},\mathfrak{A}_{\lambda_0}\right)=0.
		$$			
		
		\item\label{J8} {\bf Residual Continuity.} The family of global attractors $\{\mathfrak{A}_{\lambda}\}_{\lambda \in [0,1]}$ is  continuous in a residual set $J \subset[0,1]$, that is, for any $\lambda_0 \in J$,
		$$\lim_{\lambda \to \lambda_0}\left[{\rm dist}_{\mathcal{H}}\left(\mathfrak{A}_{\lambda_0},\mathfrak{A}_{\lambda}\right)+{\rm dist}_{\mathcal{H}}\left(\mathfrak{A}_{\lambda},\mathfrak{A}_{\lambda_0}\right)\right]=0.
		$$
		In particular, the set of
		continuity points of $\mathfrak{A}_{\lambda}$ is dense in  $[0,1]$.		
		
	\item\label{J9} {\bf Finite Dimensionality.} The compact global attractor $\mathfrak{A}_{\lambda}$ has finite fractal dimension
	$$
	\mbox{dim}^{f}_{\mathcal{ H}}\big(\mathfrak{A}_\lambda\big)<\infty, \quad \lambda\in[0,1].
	$$

\item\label{J10} {\bf Regularity.} Any trajectory $\Gamma=\{(u(t); u_t(t));  \ t\in\mathbb{R} \} \subset \mathfrak{A}_{\lambda}$ has the following   regularity
		\begin{equation}\label{regularit_traj1}
		(u_t,u_{tt})\in L^{\infty}(\mathbb{R}; \mathcal{H}).
		\end{equation}
		Moreover, there exists a constant $R > 0$ such that
		\begin{equation}\label{regularit_traj2}
		\sup_{\Gamma\subset\mathfrak{A}_{\lambda}}\sup_{t\in\mathbb{R}} 
		\|(u_t(t),u_{tt}(t))\|^2_{\mathcal{H}}  \leq R^2.
		\end{equation}	
		
		\item\label{J11} {\bf Generalized Fractal Exponential Attractor.} The dynamical system $(\mathcal{H},S_\lambda(t))$  possesses a generalized fractal exponential
attractor $\mathfrak{A}^{\exp}_{\lambda}$ with finite fractal dimension $ \big(\mbox{dim}^{f}_{\mathcal{H}^{-s}}\big(\mathfrak{A}^{\exp}_{\lambda}\big)<\infty\big) $ in the
extended space
\begin{equation}\label{interp-space}
\mathcal{H}^{-s}:= \mathcal{D}\big(A_1^{(1-s)/2}\big)\times \mathcal{D}\big(A_1^{-s/2}\big), \quad 0<s\leq1.
\end{equation}
\end{enumerate}
\end{theorem}
\begin{proof}
The proof is  a consequence of Propositions \ref{Prop-exp-decay-k2} and \ref{prop-stabi-ineq-k2} in combination with
the abstract results reminded in Appendix \ref{sec-appendixB}. It follows analogously to 
Theorem \ref{theo_main} 
except for the item \ref{J1k2} and the further items \ref{J9}-\ref{J11}.

\smallskip 	
\noindent \ref{J1k2}  
Let us consider a  bounded positively invariant set $B\subset\mathcal{H}$  and two trajectory solutions $S_\lambda(t)U^i=(u^i(t),u^i_{t}(t)), \, i=1,2,$ with initial data $U^i=(u^i_{0},u^i_{1})\in B, \, i=1,2.$

Now, in view of the stabilizability estimate \eqref{estabil-ineq-k2} we have 
\begin{align}\label{quasi-stab}
||S_\lambda(t)U^1-S_\lambda(t)U^2||_{\mathcal{H}}^2\le a_1(t) ||U^1-U^2||_{\mathcal{H}}^2 +a_2(t)\sup_{0<s<t} \big[
n(u^1(s)-u^2(s)) \big]^2, 
\end{align}
where 
$$
a_1(t) := C_B e^{-c_B t } 
, \quad  a_2(t):=C_{B}\int_0^t
e^{- c_B(t-s)} \, ds, \quad t> 0,
$$
and
$$
n(u):= \|A^{\alpha}u\|^{2}+\|u\|_{p+2}^{2}=\|u\|_{\mathcal{D}(A^{\alpha}) \cap L^{p+2}}
$$
is a compact seminorm once here the embedding  $\mathcal{D}(A^{1/2}_1)\hookrightarrow\mathcal{D}(A^{\alpha}) \cap L^{p^*}$ is compact.

 From \eqref{dependence} and \eqref{quasi-stab} one can see that the 
dynamical system
$(\mathcal{H},S_\lambda(t))$ 
satisfies the required conditions \eqref{inforce1}-\eqref{stabili} 
with
\begin{equation}
\begin{array}{l}
X:=\mathcal{D}(A^{1/2}_1), \ \  Y:= H, \ \  Z:=\{0\}, \medskip  \\
a(t):=\mathcal{Q}(t), \ \  b(t):=a_1(t), \ \ c(t):=a_2(t).
\end{array}\end{equation}

Therefore, $(\mathcal{H},S_\lambda(t))$  is asymptotically quasi-stable on  $B\subset\mathcal{H}$. In particular, by Proposition \ref{prop-quasi-stabil} it is also asymptotically smooth.

\smallskip 
\noindent \ref{J2k2}-\ref{J8}  It follows verbatim the same arguments as in \ref{I2}-\ref{I8}. 

\smallskip 
\noindent \ref{J9}-\ref{J10} 
From the first items \ref{J1k2} and \ref{J2k2}, $(\mathcal{H},S_\lambda(t))$  is asymptotically quasi-stable on the global attractor $\mathfrak{A}_\lambda\subset \mathcal{H}$ for every $\lambda\in[0,1].$ 
Thus, from Theorem \ref{theo_fractal}
one gets $\mbox{dim}^{f}_{\mathcal{ H}}\big(\mathfrak{A}_\lambda\big)<\infty$ as desired.
Moreover, since 
$a_{\infty}=\sup_{t\in\mathbb{R}^+}a_2(t)<\infty,$ then the properties \eqref{regularit_traj1}-\eqref{regularit_traj2} follows from Theorem \ref{theo_trajectories}.

\smallskip 
\noindent \ref{J11} From Proposition \ref{Prop-exp-decay-k2} and item  \ref{J1k2}, the dynamical system $(\mathcal{H},S_\lambda(t))$  is asymptotically quasi-stable on the 
 positively invariant bounded
 absorbing set $\mathcal{B}$. In what follows, for any $U_0\in\mathcal{B}$, we are going to prove that mapping 
\begin{equation}\label{mapping}
 t \mapsto S_\lambda(t)U_0:=(u(t),u_t(t))
\end{equation}
 is H\"{o}lder continuous in $\mathcal{H}^{-s}$ in \eqref{interp-space} for any $0<s\leq1.$ 
 Let us start with 
 $s=1$. From the well-posedness result, one can infer
$$
(u_t,u_{tt})\in L^{\infty}_{\text{loc}}\big(\mathcal{H}^{-1}\big).
$$
Thus, by taking $U_0\in\mathcal{B}$, $T>0$, and any $t_1, \, t_2\in[0,T],$ we get
\begin{align*}
	\|S_\lambda(t_2)U_0 - S_\lambda(t_1)U_0\|_{\mathcal{H}^{-1}}
	 \leq & \, \int_{t_1}^{t_2} \left\| \frac{d}{ds} (u(s),u_t(s)) \right\|_{\mathcal{H}^{-1}} ds \\
\leq & \, \left( \int_{0}^{T} \left\| (u_t(s),u_{tt}(s)) \right\|_{\mathcal{H}^{-1}}^2 ds\right)^{1/2} |t_2-t_1|^{1/2} \\
\leq & \, C_{\mathcal{B},T} |t_2-t_1|^{1/2},
\end{align*}
 that is, $t\mapsto S(t)U_0$ is
H\"{o}lder continuous in $\mathcal{H}^{-1}$. Besides, for $0<s<1,$ one has from the above case and interpolation theorem that 
$$
\|S_\lambda(t_2)U_0 - S_\lambda(t_1)U_0\|_{\mathcal{H}^{-s}} \leq  C_{s,\mathcal{B},T} |t_2-t_1|^{s/2}, \quad  t_1, \, t_2\in [0,T], 
$$
for some constant $C_{s,\mathcal{B},T}>0$, which proves the 
H\"{o}lder continuity in $\mathcal{H}^{-s}$.

Hence, from   Theorem
\ref{theo_attrac_exponential} the dynamical system
$(\mathcal{H},S_\lambda(t))$ has a generalized fractal exponential
attractor $\mathfrak{A}^{\exp}_\lambda$ with finite fractal dimension in
$\mathcal{H}^{-s}$ for $0<s\leq1,$ that is,
$$
\mbox{dim}^{f}_{\mathcal{H}^{-s}}\big(\mathfrak{A}^{\exp}_{\lambda}\big)<\infty, \quad 0<s\leq1.
$$

This completes the proof of Theorem \ref{theo_main-k2}.
\end{proof}

Although Theorem \ref{theo_main-k2} - \ref{J11} provides the existence of a 
generalized fractal exponential
attractor whose fractal dimension is  finite in the extended space $\mathcal{H}^{-s}, \, 0<s\leq1$, one sees from its  proof that the same methodology can not be extended to the lower limit case $s=0$, say in $\mathcal{H}^{0}=\mathcal{H}$.
 However, among all possibilities for the function $k(\cdot)$ in case \ref{k2},  in the constant scenario 
 we are supposed to 
 reach exponential attractors, that is, with finite fractal dimensional in $\mathcal{H}$. But even so, the above approach  seems to be not applicable and to circumvent the difficulty in obtaining  
 the H\"{o}lder continuity of the mapping \eqref{mapping} in $\mathcal{H}$, we are going to   replace it by a Lipschitz continuous property on a suitable space. This is exactly the goal of the next section.

\subsection{A special case:  constant $k(\cdot)$ and $p^*<\frac{2n}{n-4}$}\label{subsec-special-linear} 

It is worth pointing out that case \ref{k2} covers the class of constant functions $k(s)=\gamma>0$  for all $s\geq0$, which in turn reflects to the case \ref{k1} with $q=0$ in \eqref{P}. However, in this very special case of linear damping $\gamma u_t$ and commutative  patterns $A=A_1^{1/2}$ (see \eqref{commutative-bc}) we can go further. Indeed, in the next result we prove that  the dynamical system $(\mathcal{H},S_\lambda(t))$ has a time-dependent exponential attractor
$\mathfrak{A}^{\exp}_\lambda=\left\{\mathfrak{A}^{\exp}_\lambda(t); \, t \in \mathbb{R}\right\}  \subset \mathcal{H}$ for every $\lambda\in[0,1],$
whose sections $\mathfrak{A}^{\exp}_\lambda(t)$  have finite fractal dimension in 
$\mathcal{H}$ for all $t\in\mathbb{R}.$  This provides, in particular, the existence of an exponential attractor   $\widetilde{\mathfrak{A}^{\exp}_\lambda}$ for  the dynamical system $(\mathcal{H},S_\lambda(t))$. To this purpose, we shall work with the decomposition method motivated by the works \cite{carvalho-sonner1,carvalho-sonner2, eden-etal,efend-etal,efend-etal2,miranv-zelik,temam1}.

For $k\equiv\gamma>0$ and requiring the case \eqref{commutative-bc}, then problem \eqref{P-intro} turns into
\begin{equation} \label{main}
\left\{\begin{array}{l}\displaystyle{
	u_{tt}+\kappa A_1^{1/2} u +
	A_1 u+ \gamma u_t+f(u)= h_\lambda ,\quad t> 0,} \medskip \\
(u(0),u_t(0))=(u_0,u_1):=U_0.
\end{array}\right.
\end{equation} 
As in \eqref{sist-din}, we still denote by $(\mathcal{H},S_\lambda(t))$  the dynamical system associated with \eqref{main} in case \ref{k2} under the assumption of constant function $k(\cdot)$. Moreover, for each $\lambda \in [0,1]$ we 
 split the semigroup $S_{\lambda}(t)=S^1_{\lambda}(t)+S^2_{\lambda}(t)$ with $S^1_\lambda(t), S^2_\lambda(t), \, t\geq0,$ given as follows. 
 
 Let us consider the evolution operator
 \begin{equation}\label{def-contraction}
 S^1_{\lambda}(t):\mathcal{H} \to \mathcal{H}, \quad S^1_{\lambda}(t)U_0:= (v(t),v_t(t)),
 \end{equation}
 where $v$ the solution of the linear problem
 \begin{equation} \label{problem-v}
 \left\{\begin{array}{l}\displaystyle{
 v_{tt}+\kappa A_1^{1/2} v+A_1 v+\gamma\, v_t=h_{\lambda}, \quad t> 0,} \medskip \\
 (v(0),v_t(0))=U_0.
 \end{array}\right.
 \end{equation}
Then, we set  $S^2_{\lambda}(t):\mathcal{H} \to \mathcal{H}$ as
\begin{equation}\label{def-smoothing}
 S^2_{\lambda}(t) U_0
=S_\lambda(t)U_0-S^1_\lambda (t)U_0:= (z(t),z_t(t)),
\end{equation} 
 where $z$ solves the following problem
 \begin{equation} \label{problema-z}
\left\{\begin{array}{l}\displaystyle{
z_{tt}+\kappa A_1^{1/2} z+A_1 z+\gamma\, z_t=-f(u), \quad t> 0,} \medskip \\
(z(0),z_t(0))=(0,0).
\end{array}\right.
\end{equation}

\begin{proposition}\label{prop-exp-decay}
	Under the above setting \eqref{main}-\eqref{problema-z}, 
let us  consider a bounded set $B \subset {\mathcal{H}}$ with initial data 	$U^1_0,U^2_0 \in B$.
Then, there exist constants $c_B, C_B>0$ (depending on $B$) such that
	\begin{equation}\label{inequality-eexp}
	\|S_{\lambda}^1(t)  U^1_0-S_{\lambda}^1(t)  U^2_0\|_{\mathcal{H}} \leq C_B e^{-c_B t } \| U^1_0- U^2_0\|_{\mathcal{H}},  \quad   t > 0.
	\end{equation}
\end{proposition}
\begin{proof}
	Since the difference 
	$S_{\lambda}^1(t)  U^1_0-S_{\lambda}^1(t)  U^2_0:=(\tilde{v}(t),\tilde{v}_t(t))$ satisfies the homogeneous  linear problem related to  \eqref{problem-v} 
 \begin{equation*} 
\left\{\begin{array}{l}\displaystyle{
	\tilde{v}_{tt}+\kappa A_1^{1/2} \tilde{v}+A_1 \tilde{v}+\gamma\, \tilde{v}_t=0, \quad t> 0,} \medskip \\
(\tilde{v}(0),\tilde{v}_t(0))=U^1_0-U^2_0,
\end{array}\right.
\end{equation*}	
then the proof is a particular case of Proposition \ref{prop-stabi-ineq-k2} neglecting the precompact component.
\end{proof}

\begin{proposition}\label{prop-compac-embedd}
Under the above setting \eqref{main}-\eqref{problema-z}, let us also 
consider Assumption $\ref{assumption}$
with subcritical exponent $p<\frac{4}{n-4}$. Then, there exists a Banach space $\mathcal{W}$ such that the embedding $\mathcal{H} \hookrightarrow \mathcal{W}$ is compact and
	\begin{equation}\label{inequality-S}
	\|S_{\lambda}^2(t)  U_0^1-S_{\lambda}^2(t)  U_0^2\|_{\mathcal{H}} \leq  \mathcal{\hat{Q}}(t)\| U_0^1- U_0^2\|_{\mathcal{W}},  \quad t > 0,
	\end{equation}
for all $U_0^1,U_0^2 \in \mathcal{H}$, where 	$\mathcal{\tilde{Q}}(t)=\mathcal{\tilde{Q}}(\|U_0^1\|_{\mathcal{H}},\|U_0^2\|_{\mathcal{H}},t)$ a positive non-decreasing function.
\end{proposition}
\begin{proof}
Given  $U_0^i =(u^i_0,u^i_1)\in \mathcal{H}, \, i=1,2,$ 
and denoting $S^2_{\lambda}(t)U_0^i= (z^i(t),z^i_t(t))$, then the function $z=z^1-z^2$  satisfies
	\begin{equation}\label{problema-z-att-diferenca}
	z_{tt}+\kappa A_1^{1/2} z+A_1 z+\gamma\, z_t=f(u^2)-f(u^1), \quad (z(0),z_t(0))=(0,0).
	\end{equation}
Taking the multiplier $z_t$ in \eqref{problema-z-att-diferenca},  we have
	\begin{equation}\label{z-identity}
	\frac{1}{2}\frac{d}{dt}\mathcal{E}_z(t)=-\gamma \|z_t(t)\|_2^2+(f(u^2(t))-f(u^1(t)), z_t(t)),
	\end{equation}
where 
$$
\mathcal{E}_z(t)=\|z_t(t)\|^2+\kappa\|A_1^{1/4} z(t)\|^2+\|A^{1/2}_1 z(t)\|^2, \quad t\geq0.
$$

Let us estimate the second term on the right-hand side of \eqref{z-identity} as follows.

We first claim that there exists a power $s \in (0,2)$ such that
\begin{equation}\label{new-inequality-I}
\|f(z)-f(w)\| \leq C_f(1+\|z\|_{p^*}^{p}+\|w\|_{p^*}^{p})\|A_1^{s/4}z-A_1^{s/4}w\|, \quad  z,w \in \mathcal{D}(A_1^{1/2}),
\end{equation}
for some constant $C_f>0$.
Indeed, 
Given  $z,w \in \mathcal{D}(A_1^{1/2})\hookrightarrow L^{p^*}(\Omega)$, from  Assumption $\ref{assumption}$ and  H\"{o}lder's inequality  with $\frac{p}{p+1}+\frac{1}{p+1}=1,$ we get
	\begin{eqnarray*}\label{f-fracional}
		\nonumber \| f(z)-f(w)\|^2 & = & \int_{\Omega}\left(\int_0^1 f'(\theta z+(1-\theta)w)(z-w)\,d \theta\right)^2\,dx\\ & \leq & C_{f'}\int_{\Omega}\left(1+|z|^{2p}+|w|^{2p}\right)|z-w|^2\,dx \\
		&\leq &  C_{f'}\left(|\Omega|+\|z\|_{p^*}^{2p}+\|w\|_{p^*}^{2p}\right)\|z-w\|_{p^*}^2.  
	\end{eqnarray*}
Also, since $p<\frac{4}{n-4},$ then 
$s:=\frac{n}{2}-\frac{n}{p^*}$ satisfies
\begin{equation*}
0<s<2<\frac{n}{2},
\end{equation*}
that is, the number  $s$ satisfies the conditions of \cite[Theorem 5.1.5]{Agranovich-book}, and thus
$$
\mathcal{D}(A_1^{1/2})\hookrightarrow\mathcal{D}(A_1^{s/4}) \hookrightarrow H^s(\Omega)\hookrightarrow L^{p^*}(\Omega),
$$
from where it follows 
\eqref{new-inequality-I} for some constant $C_f>0$. Now, using Cauchy-Schwarz's inequality  and
\eqref{new-inequality-I}, we obtain
\begin{equation}\label{est-f}
(f(u^2(t))-f(u^1(t)), z_t(t))\leq C_0\|A_1^{s/4}(u^1-u^2)(t)\|[\mathcal{E}_{z}(t)]^{1/2},
\end{equation}
for some constant $C_0=C_0(\|U_0^1\|_{\mathcal{H}},\|U_0^2\|_{\mathcal{H}})>0$.
Replacing \eqref{est-f} in  \eqref{z-identity} arrive at
	\begin{equation}\label{z-derivat}
\frac{1}{2}\frac{d}{dt}\mathcal{E}_z(t)\leq-\gamma \mathcal{E}_z(t) + C_0\|A_1^{s/4}(u^1-u^2)(t)\|[\mathcal{E}_{z}(t)]^{1/2} .
\end{equation}
From \eqref{z-derivat} and regarding 
\cite[Lemma 4.1]{Showalter-book} with $\phi:=\mathcal{E}_z$ set on $[0,t],$ 
and also taking into account that $\mathcal{E}_z(0)=0$, we have
	\begin{equation}\label{z-identity-2}
[\mathcal{E}_{z}(t)]^{1/2}\leq{C}_0\int_{0}^t\|A_1^{s/4}(u^1-u^2)(\tau)\|\, d\tau,
\end{equation}
for some constant ${C}_0>0$ depending on initial data.

On the other hand, by setting $u:=u^1-u^2$, then the function  $\tilde{u}:=A_1^{(s-2)/4}u$  fulfills the problem
 \begin{equation}\label{problema-u-tilda-att-diferenca}
\left\{\begin{array}{l}\displaystyle{
	\tilde{u}_{tt}+\kappa A_1^{1/2} \tilde{u}+A_1 \tilde{u}+\gamma\, \tilde{u}_t=A_1^{(s-2)/4}(f(u^2)-f(u^1)), \quad t> 0,} \medskip \\
(\tilde{u}(0),\tilde{u}_t(0))=(A_1^{(s-2)/4} u_0,A_1^{(s-2)/4}u_1).
\end{array}\right.
\end{equation}

Taking the multiplier  $\tilde{u}_t$ in  \eqref{problema-u-tilda-att-diferenca}, and integrating the resulting expression on $(0,\tau)$, $\tau \in (0,t)$, we get
	\begin{equation}\label{estiamte-utilde-difference}
	\frac{1}{2}\mathcal{G}_u(\tau) \leq \frac{1}{2}\mathcal{G}_u(0)+\int_{0}^{\tau}(A_1^{(s-2)/4}(f(u^2(r))-f(u^1(r))), \tilde{u}_t(r)) dr,
	\end{equation}
	where
	$$\mathcal{G}_u(t)=\|A_1^{(s-2)/4} u_t(t)\|^2+\kappa\|A_1^{(s-1)/4} u(t)\|^2+\|A_1^{s/4} u(t)\|^2.$$
From Cauchy-Schwarz's inequality, since $L^2(\Omega)\hookrightarrow D(A_1^{(s-2)/4})$, and using again \eqref{new-inequality-I},  we infer
	\begin{equation*}
	\left|\int_{0}^{\tau}(A_1^{(s-2)/4}(f(u^2(r))-f(u^1(r))), \tilde{u}_t(r)) dr\right|  \leq C_0 \int_{0}^{\tau}\mathcal{G}_u(r) dr,
	\end{equation*}
	for some constant $C_0=C_0(\|U_0^1\|_{\mathcal{H}},\|U_0^2\|_{\mathcal{H}})>0$. Plugging the last estimate in \eqref{estiamte-utilde-difference},
	\begin{equation*}\label{estiamte-utilde-difference-gronwall}
	\frac{1}{2}\mathcal{G}_u(\tau) \leq \frac{1}{2}\mathcal{G}_u(0)+C_0 \int_{0}^{\tau}\mathcal{G}_u(r) dr,
	\end{equation*}
and from Gronwall's inequality, we obtain 
	\begin{equation}\label{estiamte-utilde-difference-gronwall-I}
	\mathcal{G}_u(\tau) \leq e^{C_0 \tau}\mathcal{G}_u(0),
	\end{equation}
	for some $C_0=C_0(\|U_0^1\|_{\mathcal{H}},\|U_0^2\|_{\mathcal{H}})>0$. Combining the estimates \eqref{z-identity-2} and \eqref{estiamte-utilde-difference-gronwall-I}, we finally conclude
	\begin{align*}
	\|S_{\lambda}^2(t)  U_0^1-S_{\lambda}^2(t)  U_0^2\|_{\mathcal{H}} = & \ \|(z(t),z_t(t))\|_{\mathcal{H}}\leq	[\mathcal{E}_{z}(t)]^{1/2} \\  \leq & \ {C}_0\int_{0}^t[\mathcal{G}_u(\tau)]^{1/2} \, d\tau  \leq  \ te^{C_0 t}[\mathcal{G}_u(0)]^{1/2} \\ \leq  & \ C te^{C_0 t}\|U_0^1-U_0^2\|_{\mathcal{W}},
	\end{align*}
	for some constant $C>0$ and  $C_0=C_0(\|U_0^1\|_{\mathcal{H}},\|U_0^2\|_{\mathcal{H}})>0$, where we set 
	$$\mathcal{W}:=D(A_1^{s/4})\times D(A_1^{(s-2)/4}).
	$$
Therefore, 	once the embedding $\mathcal{H} \hookrightarrow \mathcal{W}$ is compact, the estimate \eqref{inequality-S} follows by taking $\mathcal{\tilde{Q}}(t)=C te^{C_0 t}.$
\end{proof}

We are now in a position to state the result dealing with exponential attractor for the dynamical system $(\mathcal{H},S_\lambda(t))$  corresponding to \eqref{main}. 

\begin{theorem}[{\bf Exponential Attractor}]\label{theo-exponential-attractor}
Under the assumptions of Proposition $\ref{prop-compac-embedd}$, the dynamical system  	
$(\mathcal{H},S_\lambda(t))$ has a time-dependent exponential attractor
$\mathfrak{A}^{\exp}_\lambda=\left\{\mathfrak{A}^{\exp}_\lambda(t)\right\}_{ t \in \mathbb{R}}$, for every $\lambda\in[0,1],$
whose sections $\mathfrak{A}^{\exp}_\lambda(t)$  have finite fractal dimension in 
$\mathcal{H}$ for all $t\in\mathbb{R}, $ that is,  

$$
\mbox{dim}^{f}_{\mathcal{H}}\big(\mathfrak{A}^{\exp}_\lambda(t)\big)<\infty, \quad   t \in \mathbb{R}.
$$

In particular, there exists a time $T^\star>0$ such that 
$$
\widetilde{\mathfrak{A}^{\exp}_\lambda}:=\bigcup_{t \in[T^\star, 2 T^\star]} S_\lambda(t) \overline{\mathfrak{A}^{\exp}_\lambda},
$$
is an exponential attractor for $(\mathcal{H},S_\lambda(t))$, where
$\overline{\mathfrak{A}^{\exp}_\lambda}$ denotes the exponential attractor for the corresponding discrete semigroup
$
\{S_\lambda(n T^\star)\}_{n \in \mathbb{N}}.
$	
\end{theorem}
\begin{proof}
In order to employ Theorem \ref{theo-carvalho-soner}, we note that the next conditions are verified.

\begin{itemize}
	\item[{\rm (S1)}] From Proposition \ref{Prop-exp-decay-k2}, the 
	dynamical system  	
	$(\mathcal{H},S_\lambda(t))$ 
possesses a 
bounded (uniformly w.r.t. $\lambda\in[0,1]$) absorbing set $\mathcal{B}\subset \mathcal{H}$.

\item[{\rm (S2)}] From Proposition \ref{prop-exp-decay}, there exist a time $T^\star>0$ and a constant 	
  $a:=a_{T^\star, \mathcal{B}}<\frac{1}{2}$  such that $S^1_\lambda (t)$ satisfies the contraction on $\mathcal{B}$
	$$
	\|S_{\lambda}^1(T^\star)  U^1_0-S_{\lambda}^1(T^\star)  U^2_0\|_{\mathcal{H}} \leq a \| U^1_0- U^2_0\|_{\mathcal{H}}, \quad 	U^1_0,U^2_0\in\mathcal{B}.
	$$

\item[{\rm (S3)}]  From Proposition \ref{prop-compac-embedd}, there exists 
a constant  given by $b=\mathcal{\hat{Q}}(T^\star)>0$
such that
$S^2_\lambda(t)$ satisfies the  smoothing condition on $\mathcal{B}$
$$
	\|S_{\lambda}^2(T^\star)  U_0^1-S_{\lambda}^2(T^\star)  U_0^2\|_{\mathcal{H}} \leq  b \| U_0^1- U_0^2\|_{\mathcal{W}}, \quad 	U^1_0,U^2_0\in\mathcal{B},
$$
where $\mathcal{H} \hookrightarrow \mathcal{W}$ is compactly embedded.

\item[{\rm (S4)}] From Theorem \ref{theo-w-p-general} (see \eqref{dependence}), the semigroup  $S_\lambda (t)$ is  Lipschitz on  $\mathcal{B}$ with constant $L_{t}=\mathcal{Q}(t)>0$, that is,
	$$
	\|S_\lambda(t)U^1_0 - S_\lambda(t)U^2_0\|_{\mathcal{ H}} \leq L_{t} \|U^1_0-U^2_0\|_{\mathcal{ H}}, \quad    U^1_0,U^2_0 \in \mathcal{B}, \ t\geq0.
	$$	
\end{itemize}
Therefore, by means of Theorem \ref{theo-carvalho-soner} and Remark \ref{remark-last}, the conclusion of Theorem \ref{theo-exponential-attractor} is ensured.
\end{proof}

\section{Long-time dynamics: case \ref{k3}}\label{sec-dynamics-k3}

We have finally arrived at the critical case w.r.t. power $\alpha=1$.
In this case, by virtue of \eqref{equal-norm}, the $k$-function argument  $\mathcal{E}_{1}(u,u_t)$ set in \eqref{Ealpha} can be written as
\begin{equation}\label{Ealpha=1}
\mathcal{E}_{1}(u,u_t)=\|Au\|^2+\|u_t\|^2=\|A_1^{1/2}u\|^2+\|u_t\|^2,
\end{equation}
from where one sees why we lose compactness of the damping coefficient in phase space $\mathcal{H}=\mathcal{D}(A^{1/2}_1)\times H$. Therefore, this case motivated us to keep the  closedness property in  the 
definition of global attractors instead of compactness, because here the system has a degenerate damping coefficient (in the $\mathcal{H}$-topology) without control  in a neighborhood of origin and, consequently, a noncompact global attractor comes into play.

To deal with problem in this case, let us consider problem \eqref{P-intro} with critical power $\alpha=1$ and vanishing functions $f=h=0$. In this way, \eqref{P-intro} with notation  \eqref{Ealpha=1}  can be expressed as follows
\begin{equation} \label{P-k3}
\left\{\begin{array}{l}\displaystyle{
	u_{tt}+\kappa A u +
	A_1 u+k\Big(\|A_1^{1/2}u\|^2+\|u_t\|^2\Big)u_t=0,} \medskip \\
\displaystyle{u(0)=u_0, \quad u_t(0)=u_1},
\end{array}\right.
\end{equation} 
where we remember that in case \ref{k3}
we assume that
$k(\cdot)$ is a bounded Lipschitz function on  $[0,\infty)$ such that $k\equiv0$ on $[0,1]$ and $k(s)$ is strictly increasing for $s>1.$

Also, the energy functional $E(t):=E\left(u(t), u_{t}(t)\right)$ set in \eqref{energy} boils down to the following
\begin{equation}\label{energy-k3}
E(t)=\frac{1}{2}\left[ \,\|u_t(t)\|^2 +  \|A^{1/2}_1 u(t)\|^2 +
\kappa\|A^{1/2} u(t)\|^2\,\right],
\end{equation}
and satisfies the relation 
\begin{equation}\label{energy-ident}
E(t)+\int_{0}^{t}k\left(||(u(\tau),u_t(\tau))||^{2}_{\mathcal{H}}\right)\left\|u_{t}(\tau)\right\|^{2} d \tau=E(0), \quad t>0.
\end{equation}

It is worth mentioning that, under the above conditions  in the setting of problem \eqref{P-k3}, the dynamical system related to both cases \ref{k1} and \ref{k2} possesses a trivial attractor, see Theorem \ref{theo_main}-\ref{I6} and Theorem \ref{theo_main-k2}-\ref{J6}. Nonetheless, in the present case,   even under this particular scenario concerning  the source $f$ and external force $h$, we are going to show below that the dynamical system corresponding to problem \eqref{P-k3} has a nontrivial noncompact global attractor. 

\subsection{A noncompact attractor: the critical case $\alpha=1$}

Since the critical problem 
\eqref{P-k3} arises independently of the parameter $\lambda\in[0,1]$, then  
we denote by $(\mathcal{H},S(t))$  the dynamical system associated with \eqref{P-k3}. We also remind from Subsection \ref{subsec-generation-syst-din} that 
$\{S(t)\}_{t\in\mathbb{R}}$ can be seen as an evolution $C_0$-group.

Here, our main result   concerning the dynamical system $(\mathcal{H},S(t))$ reads as follows.

\begin{theorem}[{\bf Noncompact Global Attractor}]\label{theo-noncompact-attract} Under the 
setting of problem \eqref{P-k3} with function $k(\cdot)$ in case {\rm \ref{k3}},  the corresponding  dynamical system $(\mathcal{H},S(t))$  has the following global attractor 
\begin{equation}\label{bola}
\mathfrak{A}_{\kappa}=\left\{\left(u_{0} , u_{1}\right) \in \mathcal{H}; \ \left\|u_{1}\right\|^{2}+ \big\|A_1^{1 / 2} u_{0}\big\|^{2}+\kappa\|A^{1/2} u_0\|^2 \leq 1\right\}, 
\end{equation}
for $\kappa>0$ (possibly zero) small enough.
\end{theorem}
\begin{proof}
The proof will be done  in two steps  by appealing to the definition of the global attractor.

\noindent
{\it Step 1. Fully Invariance.} 
We start by noting that the linear problem 
\begin{equation} \label{P-k3-linear}
\left\{\begin{array}{l}\displaystyle{
	u_{tt}+\kappa A u +
	A_1 u=0,} \medskip \\
\displaystyle{u(0)=u_0, \quad u_t(0)=u_1},
\end{array}\right.
\end{equation} 
has a unique solution $(u(t),u_t(t))$
such that the energy set in \eqref{energy-k3}
satisfies 
$$
E(u(t), u_{t}(t))=E(u_0, u_1), \quad   t \in \mathbb{R}.
$$
From this, since $k(s)=0$ for $s\in[0,1]$, and uniqueness of the solution,  
it is easy to verify that  $\mathfrak{A}_\kappa$ is a fully invariant set with respect to $S(t)$ and, consequently,  $\mathcal{H} \backslash \mathfrak{A}_\kappa$ so is it.
Moreover, from the energy identity, it is also easy to check that the set $\mathcal{B}_{R,\kappa}:=\mathcal{B}_{R}$ given by 
$$
\mathcal{B}_{R}=\left\{\left(u_{0} , u_{1}\right) \in \mathcal{H}; \ \left\|u_{1}\right\|^{2}+ \big\|A_1^{1 / 2} u_{0}\big\|^{2}+\kappa\|A^{1/2} u_0\|^2 \leq R^2 \right\}
$$
is forward invariant with respect to $S(t)$ for every $R>1$.

\noindent
{\it Step 2. Uniform Attracting.} Given any bounded set $B\subset\mathcal{H}$, there exists $R>1$ such that $B\subset\mathcal{B}_{R}$.  Hence, below  we only need to verify that  $S(t) \mathcal{B}_{R}$ goes to $\mathfrak{A}_\kappa$ uniformly w.r.t. $\mathcal{B}_{R}$, for every $R>1$.

Let us consider $U_0=(u_0,u_1)\in\mathcal{B}_{R}, \, R>1,$ and  $S(t) U_{0}=\left(u(t), u_{t}(t)\right)$ the corresponding semigroup solution.
There are only two possibilities $U_0\in\mathfrak{A}_\kappa$ or else $U_0\in\mathcal{B}_{R}\backslash\mathfrak{A}_\kappa.$
If $U_0\in\mathfrak{A}_\kappa$, then  $S(t)U_0\in S(t)\mathfrak{A}_\kappa=\mathfrak{A}_\kappa,$ and the conclusion follows trivially. Thus, in what follows we take $U_0\in\mathcal{B}_{R}\backslash\mathfrak{A}_\kappa.$ We first claim that 
\begin{equation}\label{limit}
\lim_{t\to +\infty} \left[2 E(S(t)U_0)\right]
=1.
\end{equation}
Indeed, let us suppose  that it does not hold. Thus,
due to the invariance of $\mathcal{H} \backslash \mathfrak{A}_\kappa$ and since the mapping $t \mapsto 2 E(S(t)U_0)$ is non-increasing,  
there exists $R_{0}>1$ such that
\begin{equation}\label{mig1}
2 E(S(t)U_0) \geq R_{0}, \ t>0, \ \ \text { and } \ \ \lim_{t\to +\infty} \left[2 E(S(t)U_0)\right]=R_{0}.
\end{equation}
Now we remember (analogously to \eqref{dependenceIII}) that 
\begin{equation}\label{mig2}
||S(t) U_{0}||_{\mathcal{H}}^2  \leq 2 E(S(t)U_0) \leq 
(1+\kappa
\mu_0)||S(t) U_{0}||_{\mathcal{H}}^2,
\end{equation}
and taking $\kappa>0$ small enough (or possibly zero) so that $\kappa<\frac{R_0-1}{\mu_0}$,  
then $R_\kappa:= \frac{R_0}{1+\kappa \mu_0}>1$. From \eqref{mig1}-\eqref{mig2} one sees that $\left\|S(t) U_{0}\right\|_{\mathcal{H}}^{2} \geq R_\kappa>1$ for all $t\geq0$, and 
 from the assumption on $k(\cdot)$  in the range $(1,\infty)$, we infer
\begin{equation}\label{mig3}
k\left(\left\|S(t) U_{0}\right\|_{\mathcal{H}}^{2}\right) \geq k(R_{\kappa}) := k_{0}>0 \ \text { for all } \ t \geq 0.
\end{equation}

Going back to problem 
\eqref{P-k3}, taking the multiplier $u_t$, and using \eqref{mig3}, we get 
\begin{equation}\label{mig4}
\frac{d}{d t} E(S(t)U_0) + 2 k_{0}\left\|u_{t}(t)\right\|^{2} \leq 0, \quad t>0.
\end{equation}
Using \eqref{mig4} and since  $k_{1}=\sup \{k(s); \, s \geq 1\}<\infty$, then problem \eqref{P-k3}
behaves like in case \ref{k2} and similar to Proposition \ref{Prop-exp-decay-k2} (see \cite[Remark 7]{jorge-narciso-DCDS} for more details)  one shows that there exists a constant $c>0$ (which may depend on $U_0$ and $k_0$) such that
$$
E(S(t)U_0) \leq 3 E(U_0) e^{- c \, t}, \quad  t>0,
$$
whenever $2E(S(t)U_0) \geq R_{0}>1,$ which is a contraction to \eqref{mig1}.
Therefore, \eqref{limit} holds true and is uniform with respect to $U_{0} \in \mathcal{B}_{R} \backslash \mathfrak{A}_\kappa$. Additionally,
\begin{equation*}
\lim_{t\to +\infty} \left[2 E(S(t)U_0)\right]^{1/2}
=1, \quad U_{0} \in \mathcal{B}_{R} \backslash \mathfrak{A}_\kappa.
\end{equation*}

Finally, for every $U_{0} \in \mathcal{B}_{R} \backslash \mathfrak{A}_\kappa$, we observe that
\begin{align*}
\operatorname{dist}\left(S(t) U_{0}, \mathfrak{A}_\kappa\right) \, \leq \, & \, \left\|S(t) U_{0}-\frac{S(t) U_{0}}{[2E(S(t)U_0)]^{1/2}}\right\|_{\mathcal{H} } \\
 \, = \, & \,
\frac{ \left\|S(t) U_{0}\right\|_{\mathcal{H}}
\Big[[2E(S(t)U_0)]^{1/2}-1\Big]}{[2E(S(t)U_0)]^{1/2}} 
\\
\, \leq \, & \, \Big[[2E(S(t)U_0)]^{1/2}-1\Big],
\end{align*}
where in the last inequality we use \eqref{mig2}, from where one concludes
$$
\lim_{t\to +\infty}
\operatorname{dist}_{\mathcal{H}}\left(S(t)(\mathcal{B}_{R} \backslash \mathfrak{A}_\kappa) , \mathfrak{A}_\kappa\right) =\lim_{t\to +\infty} \left[\sup_{U_0\in \mathcal{B}_{R} \backslash \mathfrak{A}_\kappa} 
\operatorname{dist}\left(S(t) U_{0}, \mathfrak{A}_\kappa\right)\right]=0.  
$$

This completes the proof that the
bounded closed set
$\mathfrak{A}_\kappa $ given in \eqref{bola} uniformly
attracts $\mathcal{B}_{R}$ for any $R>1$, and bounded subsets $B\subset\mathcal{H}$ as well.

Therefore,
 $\mathfrak{A}_\kappa$ is a global attractor for the dynamical  system $(\mathcal{H},S(t))$ generated by problem \eqref{P-k3}. 
\end{proof}

\begin{Corollary}
Under the hypotheses of Theorem $\ref{theo-noncompact-attract}$, the dynamical system  $(S(t),\mathcal{H})$ associated with problem \eqref{P-k3} does not have a compact global attractor.
\end{Corollary}
\begin{proof}
	It is a directly consequence of the uniqueness of a global attractor.
\end{proof}

\begin{Corollary}\label{cor-bola-attrac}
Under the hypotheses of Theorem $\ref{theo-noncompact-attract}$ with $\kappa=0$, then the closed ball 
\begin{equation*}
\mathfrak{B}_{0}=\left\{\left(u_{0} , u_{1}\right) \in \mathcal{H}; \ \left\|u_{1}\right\|^{2}+ \big\|A_1^{1 / 2} u_{0}\big\|^{2}\leq 1\right\}
\end{equation*}
is the global attractor for  $(S(t),\mathcal{H})$.
\end{Corollary}

\begin{remark}
The result stated in Corollary $\ref{cor-bola-attrac}$ corresponds to the one approached \cite[Proposition 5.3.9]{chueshov-solo}. Additionally, we note that in this limit case $\kappa=0$, one can relax the assumption of $k(\cdot)$  on the interval $[0,1]$, for instance with  arbitrary behavior on $[0,1]$ instead vanishing on it, see \cite[Subsect. 5.3.3]{chueshov-solo}. 
\end{remark}

\appendix

\section{Auxiliary proofs}\label{sec-appendix}

\subsection{Completion of the proof of Theorem \ref{theo-existence}}\label{subsec-A1}

By using  Assumption \ref{assumption}, we prove now that operator $\mathcal{M}:\mathcal{H}\to\mathcal{H}$  defined in (\ref{def_B}) is  a locally Lipschitz continuous operator on $\mathcal{H}$.   

In fact, let us consider $R >0$ and   $U= (u,v), \widetilde{U}=(\widetilde{u},\widetilde{v})\in \mathcal{H}$ such that $||U||_{\mathcal{H}},\,||\widetilde{U}||_{\mathcal{H}}\le R$. From definition (\ref{def_B}) and the $\mathcal{H}$-norm, we have
\begin{eqnarray}\label{Sem}
||\mathcal{M}(U)-\mathcal{M}(\widetilde{U})\,||_{\mathcal{H}}
&=&
\left\|\kappa [A u-A \widetilde{u}]+\gamma\left[||U||^{2q}_{\mathcal{H}^{\alpha}}v
-||\widetilde{U}||^{2q}_{\mathcal{H}^{\alpha}}\widetilde{v}\right]
-[f(u)-f(\widetilde{u})]\right\|\nonumber\\
&\le& \kappa \|U-\widetilde{U}\|_{\mathcal{H}}
+\gamma\underbrace{\left\|||U||^{2q}_{\mathcal{H}^{\alpha}}v
	-||\widetilde{U}||^{2q}_{\mathcal{H}^{\alpha}}\widetilde{v}\right\|}_{\mathcal{I}_1}
+ \underbrace{\|f(u)-f(\widetilde{u})\|}_{\mathcal{I}_2}.
\end{eqnarray}
In what follows, we are going to estimate the terms $\mathcal{I}_1$ and $\,\mathcal{I}_2$ above. To estimate the term $\mathcal{I}_1$, we first note that  
\begin{eqnarray}\label{Sem_I}\,||U||^{2q}_{\mathcal{H}^{\alpha}}v-||\widetilde{U}||^{2q}_{\mathcal{H}^{\alpha}}\widetilde{v}=\frac{1}{2}\left[||U||^{2q}_{\mathcal{H}^{\alpha}}+||\widetilde{U}||^{2q}_{\mathcal{H}^{\alpha}}\right][v-\widetilde{v}]+\frac{1}{2}\left[||U||^{2q}_{\mathcal{H}^{\alpha}}-||\widetilde{U}||^{2q}_{\mathcal{H}^{\alpha}}\right][v+\widetilde{v}].
\end{eqnarray}
Using that $\mathcal{D}(A^{1/2}_1)\subseteq \mathcal{D}(A) \subseteq \mathcal{D}(A^{\alpha})$ and $\|A^{1/2}_1u\|^2=\|A u\|$ for $u\in \mathcal{D}(A_1^{1/2})$, we can estimate the first term of the sum in (\ref{Sem_I}) as follows 
\begin{eqnarray*}
	\frac{1}{2}\left\| \left[||U||^{2q}_{\mathcal{H}^{\alpha}}
	+||\widetilde{U}||^{2q}_{\mathcal{H}^{\alpha}}\right][v-\widetilde{v}]
	\,\right\|\le C_R\|v-\widetilde{v}\|\le C_R
	||U-\widetilde{U}||_{\mathcal{H}}.
\end{eqnarray*}
We emphasize again that this is the precise moment where the assumption $q\geq1/2$ is crucial in our computations. Indeed, once we have
 $2q\geq 1$, then we can use Lemma \ref{Lemma-Haraux} to obtain
\begin{eqnarray*}
\nonumber \frac{1}{2}\left\| \left[\,||U||^{2q}_{\mathcal{H}^{\alpha}}
-||\widetilde{U}||^{2q}_{\mathcal{H}^{\alpha}}\,\right][v+\widetilde{v}]
\,\right\|& \le & q\max\left\{||U||_{\mathcal{H}^{\alpha}},||\tilde{U}||_{\mathcal{H}^{\alpha}}\right\}^{2q-1}\|v+\widetilde{v}\|||U-\widetilde{U}||_{\mathcal{H}^{\alpha}}\\
& \le & C_R ||U-\widetilde{U}||_{\mathcal{H}}
\end{eqnarray*}
for some $C_R>0$. Therefore, collecting the last two estimates, we deduce from \eqref{Sem_I} that 
$$
\mathcal{I}_1\le C_R||U-\widetilde{U}||_{\mathcal{H}}.
$$
Lastly, we are going to estimate the term $\mathcal{I}_2$. Using Mean Value Theorem, conditions (\ref{assumption_f'}), H\"older's inequality with $\frac{p}{p+1}+\frac{1}{p+1}=1$ and the embedding $\mathcal{D}(A^{1/2}_1)\hookrightarrow L^{p^*}(\Omega)$, we have
\begin{eqnarray*}
	\mathcal{I}_2&\le& 2^{2p}C_{f'}\left[|\Omega|+\|u\|^{p^*}_{p^*}+\|\widetilde{u}\|^{p^*}_{p^*}\right]^{\frac{p}{p^*}}\|u-\widetilde{u}\|_{p^*}\\
	&\le& 2^{2p}C_{f'}C_{|\Omega|}\left[|\Omega|+C^{p^*}_{|\Omega|}\left(\|A^{1/2}_1 u\|^{p^*}+\|A^{1/2}_1 \widetilde{u}\|^{p^*}\right)\right]^{\frac{p}{p^*}}\|A^{1/2}_1 u-A^{1/2}_1 \widetilde{u}\|\\
	&\le& C_R||U-\widetilde{U}||_{\mathcal{H}}.
\end{eqnarray*}
replacing $\mathcal{I}_1$ and $\mathcal{I}_2$ in (\ref{Sem}), we obtain
\begin{equation*}
||\mathcal{M}(U)-\mathcal{M}(\widetilde{U})\,||_{\mathcal{H}} \le   C_{R}||U-\widetilde{U}||_{\mathcal{H}},
\end{equation*}
for some $C_R>0$, which proves the desired.

\qed

\subsection{Proof of Proposition \ref{Nak-lemma}}\label{sec-proof-nak-lemma}

\noindent  {\bf Proof of \ref{N1-nak}.} In this case, under the assumptions of Proposition \ref{Nak-lemma}, we are going to prove the inequality \eqref{nak-thesis1}. Indeed, we first 
  observe that  it  holds true    for $0\leq t\leq 1$. Then, let us prove it for $t>1$. 
  
  Setting the function
$$
\beta(t):= \left(C_{0}^{-1} \rho(t-1)^{+}+\left(\sup _{0 \leqslant s \leqslant 1} \phi(s)\right)^{-\rho}\right)^{-1 / \rho}, \quad t\geq0,
$$
we initially claim that the following property holds: if for some $t\geq0$, the inequality holds
\begin{equation}\label{ineq1}
\phi(t) \leqslant \beta(t)+[K(t)]^{1 /(\rho+1)},
\end{equation}
then the next inequality holds as well
\begin{equation}\label{ineq2}
\phi(t+1) \leqslant \beta(t+1)+[K(t+1)]^{1 /(\rho+1)}.
\end{equation}
Indeed, let us suppose that \eqref{ineq2} does not hold, that is 
\begin{equation}\label{ineq3}
\phi(t+1) > \beta(t+1)+[K(t+1)]^{1 /(\rho+1)}.
\end{equation}
Now, if $\phi(t) \leqslant [K(t)]^{1 /(\rho+1)}$, then \eqref{ineq3} and \eqref{ineq1} imply
$$
\phi(t+1) > [K(t+1)]^{1 /(\rho+1)}\geq [K(t)]^{1 /(\rho+1)} > \phi(t),
$$
where we have used that $K(t)$ is non-increasing. Thus, from \eqref{nak-hyp}, we obtain
$$
\phi(t+1)<[K(t+1)]^{1 /(\rho+1)},
$$
which is a contradiction with \eqref{ineq3}. Thus, we can infer that  $\phi(t)>[K(t)]^{1 /(\rho+1)}$ and, consequently,  
$$
\varphi(t) :=\phi(t)-[K(t)]^{1 /(\rho+1)}>0 \quad \text { and } \quad \varphi(t+1) :=\phi(t+1)-[K(t+1)]^{1 /(\rho+1)}>0.
$$ 
From this, one sees
\begin{equation}\label{ineq4}
{[\varphi}(t)]^{\rho+1}+K(t) \leqslant\left({\varphi}(t)+[K(t)]^{1 /(\rho+1)}\right)^{\rho+1}=[\phi(t)]^{1+\rho},
\end{equation}
and using again \eqref{nak-hyp}, one gets
\begin{equation}\label{ineq5}
{[\varphi}(t)]^{1+\rho} \leqslant C_{0}({\varphi}(t)-{\varphi}(t+1)).
\end{equation}
We additionally set
$$
\psi(t):= \varphi^{-\rho}(t) \quad \text { and } \quad \psi(t+1):= \varphi^{-\rho}(t+1).
$$
Then, a straightforward integral computation along with \eqref{ineq5} lead to
$$
\begin{aligned}
\psi(t+1)-\psi(t) &=-\int_{0}^{1} \frac{d}{d \tau}(\tau \varphi(t)+(1-\tau) \varphi(t+1))^{-\rho} d \tau \\
&=\rho \int_{0}^{1}(\tau \varphi(t)+(1-\tau) \varphi(t+1))^{-(1+\rho)} d \tau(\varphi(t)-\varphi(t+1)) \\
& \geqslant \rho \varphi^{-(1+\rho)}(t)(\varphi(t)-\varphi(t+1)) \geqslant \rho C_{0}^{-1}.
\end{aligned}
$$
Thus,  we have
\begin{equation*}
\varphi^{-\rho}(t+1)=\psi(t+1) \geqslant \psi(t)+\rho C_{0}^{-1}=\varphi^{-\rho}(t)+\rho C_{0}^{-1},
\end{equation*}
that is,
\begin{equation}\label{ineq6}
\varphi(t+1)\leq \big(\varphi^{-\rho}(t)+\rho C_{0}^{-1}\big)^{-1/\rho},
\end{equation}

From \eqref{ineq6} and \eqref{ineq1}, we finally arrive at
$$
\begin{aligned}
\phi(t+1) &=\varphi(t+1)+[K(t+1)]^{1 /(\rho+1)}  \\
&  \leqslant\left(\varphi^{-\rho}(t)+\rho C_{0}^{-1}\right)^{-1 / \rho}+[K(t+1)]^{1 /(\rho+1)} \\
& \leqslant\left(\rho C_{0}^{-1}(t-1)^{+}+\sup _{0 \leqslant s \leqslant 1} [\phi(s)]^{-\rho}+\rho C_{0}^{-1}\right)^{-1 / \rho}+[K(t+1)]^{1 /(\rho+1)} \\
& \leqslant \beta(t+1)+[K(t+1)]^{1 /(\rho+1)},
\end{aligned}
$$
which contradicts  \eqref{ineq3}. This concludes the proof of our initially assertion.

To conclude the proof, we note that for any $t>1$ real, we can write $t=n+r$ with 
$n\in\mathbb{N}$ and $0\leq r<1.$ From the beginning, \eqref{ineq1} holds true for $0\leq r<1$ and from \eqref{ineq2},
\begin{equation*}
\phi(r+1) \leqslant \beta(r+1)+[K(r+1)]^{1 /(\rho+1)}.
\end{equation*}
Using this assertion $n$ times, we infer 
\begin{equation*}
\phi(t)=\phi(r+n) \leqslant \beta(r+n)+[K(r+n)]^{1 /(\rho+1)}=\beta(t)+[K(t)]^{1 /(\rho+1)},
\end{equation*}
as desired. This concludes the proof of \eqref{nak-thesis1}.

\smallskip

\noindent  {\bf Proof of  \ref{N2-nak}.} To the proof of  
\eqref{nak-thesis2}, we initially set 
$$
\beta(t):= \sup _{0 \leqslant s \leqslant 1} \phi(s)\left(\frac{C_{0}}{1+C_{0}}\right)^{[t]}, \quad t\geq0,
$$
and proceed verbatim as in the first case.
 \qed

\subsection{Proof of Proposition \ref{Khan-lemma}}\label{sec-proof-khan-lemma}

	From the hypothesis of Proposition \ref{Khan-lemma}, one sees
	\begin{equation}\label{khan1}
	\left\{\begin{array}{lll}
	u^{n} \rightarrow u & \text { weakly-star in }& L^{\infty}\big(s, T ; \mathcal{D}(A^{1/2}_1)\big), \\
	u_{t}^{n} \rightarrow u_{t} & \text { weakly-star in }& L^{\infty}(s, T ; H),
	\end{array}\right.
	\end{equation}
	and from the Aubin-Lions compactness theorem (see e.g. Simon \cite{Simon}), we also have
	\begin{equation}\label{khan2}
	u^{n} \rightarrow u \quad \text { strongly in } \quad C([s, T] ; H).
	\end{equation}
	Additionally, by using Lemma 8.1 in Lions and Magenes \cite{lions-magenes} (see on p. 275 therein), \eqref{khan1} also implies that $u^{n}$  is bounded in  $C_{s}\big(s, T ; \mathcal{D}(A^{1/2}_1)\big)$, and then $u^{n}(t)$   is bounded in $\mathcal{D}(A^{1/2}_1)$ for all $t \in[s, T] .$ From this and \eqref{khan2} one gets
	\begin{equation}\label{khan3}
	u^{n}(t) \rightarrow u(t) \quad \text { weakly in  } \quad \mathcal{D}(A^{1/2}_1), \ s\leq t\leq T,
	\end{equation}
	and due to  the compact embedding theorem, we infer 
	\begin{equation}\label{khan4}
	\widehat{f}\left(u^{n}(t)\right) \rightarrow \widehat{f}(u(t)) \quad \text { strongly in } \quad L^{1}(\Omega),  \ s\leq t\leq T.
	\end{equation}
	where we remember that  $\widehat{f}(u)=\int_0^{u}f(\tau)d\tau$. Also, from \eqref{khan1}, assumptions on $f$ and again \eqref{khan2}, we have
	\begin{equation}\label{khan5}
	\left( f\left(u^{n}\right), u_{t}^{n}\right) \rightarrow \left( f(u), u_{t}\right) \quad \text { strongly in } \quad  L^{1}(s, T).
	\end{equation}
	
	Now, regarding that
	$$
	\frac{\partial}{\partial t} \int_{\Omega} \widehat{f}\left(u^{n}(x,t)\right) d x=\left( f\left(u^{n}(t)\right), u_{t}^{n}(t)\right),
	$$	
	we get
	$$
	\int_{s}^{t}\big( f(u^n(\tau)), u^n_{t}(\tau)\big) d \tau=\int_{\Omega} \widehat{f}(u^n(x,t)) d x-\int_{\Omega} \widehat{f}(u^n(s, x)) d x.
	$$	
	From this identity (which also holds true for $u$) and from the limits \eqref{khan4}-\eqref{khan5}, we finally arrive at
	\begin{align*}
	\lim _{n \rightarrow \infty} \lim _{m \rightarrow \infty} \int_{s}^{T}\big( f\left(u^{n}(t)\right)-&f\left(u^{m}(t)\right), u_{t}^{n}(t)-u_{t}^{m}(t)\big) d t\\
	= & \, \lim _{n \rightarrow \infty} \int_{\Omega} \widehat{f}\left(u^{n}( x,T)\right) d x+\lim _{m \rightarrow \infty} \int_\Omega \widehat{f}\left(u^{m}( x,T)\right) d x \\
	&-\lim _{n \rightarrow \infty} \int_{\Omega} \widehat{f}\left(u^{n}( x,s)\right) d x-\lim _{m \rightarrow \infty} \int_{\Omega} \widehat{f}\left(u^{m}( x,s)\right) d x
	\\
	& -\lim _{n \rightarrow \infty} \lim _{n \rightarrow \infty} \int_s^T\int_{\Omega} f\left(u^{n}(t, x)\right) u_{t}^{m}( x,t) d x d t \\
	& -\lim _{n \rightarrow \infty} \lim _{n \rightarrow \infty} \int_{s}^{T}\int_\Omega f\left(u^{m}(x,t)\right) u_{t}^{n}(x,t) d x d t \\
	= & \ 2 \int_{\Omega} f(u(x,T)) d x 
	-2 \int_{\Omega} f(u(x,s)) d x \\
	& -2 \int_s^T\big(  f(u(t)) ,u_{t}( t)\big) d x d t \\
	= & \ 0,
	\end{align*}
which proves the desired in \eqref{Khan-limit}.	
 \qed

\section{Auxiliary facts on dynamical systems}\label{sec-appendixB}

The concepts and results on dynamical systems reminded below  can be found e.g. in
 \cite{babin,carvalho-sonner1,chueshov-solo,chueshov-lasiecka-2005, chueshov-white, chueshov-yellow,eden-etal, hale,Robinson-Olson-Hoang,Lady,Robinson-book, temam1}. 
 Hereafter, $(H,S_t)$ stands for  a dynamical system
consisting of a $C_0$-semigroup $S_t, \ t\geq0,$ defined on a
Banach space $(H,\|\cdot\|)$.
 
\begin{itemize}

\item  The dynamical system $({H},S_t)$ is said to be  {\it dissipative} if it
possesses a bounded {\it absorbing set}, that is, a bounded set $
\mathscr{B}\subset{H}$ such that for any bounded set $B \subset H$
there exists $t_B \ge 0$ satisfying
$$
S_tB \subset \mathscr{B}, \quad \forall \, t \geq t_B.
$$

\item   One says that $( {H},S_t)$ is
{\it asymptotically smooth} if for any bounded positively invariant set
$B\subset  {H}$ $(S_tB \subseteq B$ for
	all $t\ge 0)$, there exists a compact set $\mathscr{M} \subset
\overline{B}$ such that
\begin{equation}\label{aymp-smooth}
\lim_{t\to +\infty} \mbox{dist}_{ {H}}(S_tB,\mathscr{M}) = 0,
\end{equation}
where $\mbox{dist}_{H}(\cdot,\cdot)$ stands for the {\it Hausdorff semidistance}\footnote{The  Hausdorff semidistance of two non-empty subsets $A, B \subset H$   is given by
		$$
		\operatorname{dist}_{H}(A, B):=\sup _{x \in A} \mbox{dist}(x,B)=\sup _{x \in A} \inf _{y \in B}\|x-y\|.
		$$} in ${H}$.

\item $\left(H, S_{t}\right)$ is said to be {\it asymptotically compact} iff there exists an {\it attracting compact set} $\mathscr{M},$ that is, for any bounded set $B$ one has that 
\eqref{aymp-smooth} holds.

\item  A {\it global attractor} for $({H},S_t)$ is a bounded closed set
$\mathscr{A} \subset  {H}$ which is fully invariant and uniformly
attracting, that is, $S_t\mathscr{A} = \mathscr{A}$ for all $t\ge
0$ and for every bounded subset $B \subset  {H}$,
$$
\lim_{t\to +\infty} \mbox{dist}_{{H}}(S_tB,\mathscr{A}) = 0.
$$

\item A {\it global minimal attractor} for $(H,S_t)$ is a bounded closed
set $\mathscr{A}_{\min} \subset  H$ which is positively
invariant $(S_t\mathscr{A}_{\min}\subseteq \mathscr{A}_{\min})$
and attracts uniformly every point, that is,
$$
\lim_{t\to +\infty}
\mbox{dist}(S_tU_0,\mathscr{A}_{\min}) = 0, \ \
\mbox{for any} \ \ U_0\in\mathcal{ H},
$$
and $\mathscr{A}_{\min}$ has no proper subsets possessing these
two properties.

\item The {\it unstable manifold}  emanating
from a set $\mathscr{N}$, denoted by $\mathbf{M}^u(\mathscr{N})$,
is a set of $H$ such that for each
$U_0\in\mathbf{M}^u(\mathscr{N})$ there exists a full trajectory
$\Gamma=\{\mathbf{U}(t); \ t\in\mathbb{R} \}$ satisfying
$$
\mathbf{U}(0)=U_0 \quad \mbox{and} \quad \lim_{t\to-\infty}
\mbox{dist}(\mathbf{U}(t),\mathscr{N})=0.
$$

 \item  The dynamical system $(H,S_t)$ is said to be
{\it gradient} iff there exists a {\it strict Lyapunov functional}
on $H$, that is, there exists a continuous functional $\Phi(z)$ such
that the function $t \mapsto \Phi(S_tz)$ is nonincreasing for any
$z\in H$, and the equation $\Phi(S_tz)=\Phi(z)$ for all $t>0$ and some
$z\in H$ implies that $S_tz=z$ for all $t>0$.

\item  The  {\it Kolmogorov  $\varepsilon-$entropy} ${H}_{\varepsilon}(\mathscr{M})$ of a compact set $\mathscr{M}\subset H$   is given by
\begin{equation}\label{def-kolmog}
{H}_{\varepsilon}(\mathscr{M})=\ln N(\mathscr{M}, \varepsilon), \quad \varepsilon>0,
\end{equation}	
where $N(\mathscr{M}, \varepsilon)$ is the minimal number of closed sets of the diameter not greater than $2 \varepsilon $ which cover the compact $\mathscr{M}.$  The {\it fractal dimension} $\operatorname{dim}_{f} \mathscr{M}$ of $\mathscr{M}$ is defined by the formula
$$
\mbox{dim}^{f}_{H}\mathscr{M}=\limsup _{\varepsilon \rightarrow 0} \frac{{H}_{\varepsilon}(\mathscr{M})}{\ln (1 / \varepsilon)}.
$$

\item A compact set $\mathscr{A}_{\exp}\subset{H}$ is said to be a
{\it fractal exponential attractor} of the dynamical system
$({H},S_t)$ if $\mathscr{A}_{\exp}$ is a positively
invariant set of finite fractal dimension  and for
every bounded set $B\subset {H}$ there exist positive
constants $t_B$, $C_B$ and $\sigma_B$ such that
$$
\mbox{dist}_{{H}}(S_tB,\mathscr{A}_{\exp}) \leq C_B \, e^{-\sigma_B(t-t_B)}, \quad t\geq t_B.
$$
If the exponential attractor has  finite fractal
dimension in some extended space $\widetilde{{H}}\supseteq
{H}$, one calls this exponentially attracting set as a
{\it generalized fractal exponential attractor}.

\end{itemize}

The first results below deal with the existence and characterization of global attractors. To their statements, we follow more closely the works \cite{chueshov-white,chueshov-yellow}.
 
 \begin{theorem}[{\cite[Proposition 2.10]{chueshov-white}}] \label{theo_khanma} Assume that for any bounded
 	positively invariant set $B\subset {H}$ and for any $\varepsilon
 	>0$, there exists $T=T(\varepsilon , B)$ such that
 	$$
 	\Vert S_tz_1-S_tz_2 \Vert \leq \varepsilon +
 	\phi_T(z_1,z_2), \quad \forall \,  z_1,z_2 \in B,
 	$$
 	where $\phi_T : B  \times B \to \mathbb{R}$ satisfies
 	\begin{equation} \label{lim-PhiT}
 	\liminf_{n \rightarrow\infty}\liminf_{m
 		\rightarrow\infty}\phi_T(z_{n},z_{m})=0,
 	\end{equation}
 	for any sequence $(z_{n})$ in $B$. Then $(H,S_t)$ is an
 	asymptotically smooth dynamical system.
 \end{theorem}


\begin{proposition}[{\cite[Proposition 7.1.4]{chueshov-yellow}}]\label{prop-equiv}
 Let $\left(H, S_{t}\right)$ be a dissipative dynamical system. Then,  $ \left(H, S_{t}\right)$ is asymptotically compact  if and only if	$\left(H, S_{t}\right)$ is asymptotically smooth.
\end{proposition}

\begin{theorem}[{\cite[Theorem 2.3]{chueshov-white}}]\label{theo-attract}
 Let $\left(H, S_{t}\right)$ be a dissipative dynamical system. Then $\left(H, S_{t}\right)$ possesses a compact global attractor $\mathscr{A}$ if and only if $\left(H, S_{t}\right)$ is asymptotically smooth.
\end{theorem}

 \begin{theorem}[{\cite[Theorem 2.28]{chueshov-white}}] \label{theo_geometric_struc}
 	Let $\mathscr{N}$ be the set of stationary points\footnote{$\mathscr{N}=\left\{v \in H: S_{t} v=v \text { for all } t \geq 0\right\}$.} of $(H,S_t)$
 	and assume that $(H,S_t)$ possesses a compact global attractor
 	$\mathscr{A}$. If there exists a strict Lyapunov functional on
 	$\mathscr{A}$, then $\mathscr{A}=\mathbf{M}^u(\mathscr{N})$. Moreover the global attractor $\mathscr{A}$ consists of full trajectories $\Gamma=\{U(t): t \in \mathbb{R}\}$ such that
 	$$
 	\lim _{t \rightarrow-\infty} \operatorname{dist}(U(t), \mathscr{N})=0 \quad \text { and } \quad \lim _{t \rightarrow+\infty} \operatorname{dist}(U(t), \mathscr{N})=0.
 	$$
 \end{theorem}

 \begin{theorem}[{\cite[Theorem 7.5.10]{chueshov-yellow}}] \label{theo_minimal_attra}
 	Assume that a gradient dynamical system $({H},S_t)$ possesses a
 	compact global attractor $\mathscr{A}$. Then for any $z\in H$ we
 	have
 	$$
 	\lim_{t\to +\infty} \mbox{\rm dist}(S_tz,\mathscr{N}) = 0,
 	$$
 	that is, any trajectory stabilizes to the set $\mathscr{N}$ of
 	stationary points. In particular,
 	$\mathscr{A}_{\min}=\mathscr{N}.$
 \end{theorem}

  In the next results we deal with a family of attractors $ \{\mathscr{A}_{\lambda}\}_{\lambda \in \Lambda}\subset  H .$\footnote{Here, $\Lambda$ stands for complete metric space.} To their statements, we follow the references \cite{Robinson-Olson-Hoang,Robinson-book}.  We first remind the concepts of   upper semicontinuity and (residual) continuity as follows.
  
 \begin{itemize}
 	\item The family of attractors  $\{\mathscr{A}_{\lambda}\}_{\lambda \in \Lambda}$ is {\it upper semicontinuous} at the point $\lambda_0$ iff 	  
 	$$\lim_{\lambda \to \lambda_0}{\rm dist}_{H}\left(\mathscr{A}_{\lambda},\mathscr{A}_{\lambda_0}\right)=0,$$
where, as above,  $\mbox{dist}_{H}(\cdot,\cdot)$ stands for the {Hausdorff semidistance} in ${H}$.

\item  The family of attractors  $\{\mathscr{A}_{\lambda}\}_{\lambda \in \Lambda}$ is {\it continuous} at the point $\lambda_0$ when 
$$\lim_{\lambda \to \lambda_0}\left[{\rm dist}_{H}\left(\mathscr{A}_{\lambda_0},\mathscr{A}_{\lambda}\right)+{\rm dist}_{H}\left(\mathscr{A}_{\lambda},\mathscr{A}_{\lambda_0}\right)\right]=0.$$
  \end{itemize} 

 \begin{theorem}[{\cite[Theorem 10.16]{Robinson-book}}]\label{conti-Robinson}
 Let $\{S_t^{\lambda}\}_{\lambda \in \Lambda}$ be a family of semigroups on $H$ possessing global attractors $\mathscr{A}_{\lambda}$ for $\lambda \in \Lambda$. Let us additionally assume:
 	\begin{itemize}
 		\item[{\rm (a)}] the attractors $\mathscr{A}_{\lambda}$ are uniformly bounded, i.e., there exists a bounded set $\mathscr{B}_0 \subset H$ such that $\mathscr{A}_{\lambda}\subset \mathscr{B}_0$, for all $\lambda \in \Lambda$;
 		
 		\item[{\rm(b)}] there exists $t_0\geq 0$ such that 
 		$$\lim_{\lambda \to \lambda_0}\sup_{z \in \mathscr{B}_0}||S_t^{\lambda}z-S_t^{\lambda_0}||=0, \quad \forall \,  t\geq t_0.$$
 	\end{itemize}
 	Then, the family of attractors  $\{\mathscr{A}_{\lambda}\}_{\lambda \in \Lambda}$ is upper semicontinuous at the point $\lambda_0$.
 \end{theorem}

 \begin{theorem}[{\cite[Theorem 5.2]{Robinson-Olson-Hoang}}]\label{theo-ROH}
 Let $\{S_t^{\lambda}\}_{\lambda \in \Lambda}$ be a family of semigroups on $H$. Let us additionally suppose:
 
 \begin{itemize}
 	\item[{\rm (a)}]  $S^{\lambda}_t$ has a global attractor $\mathscr{A}_{\lambda}$ for every $\lambda \in \Lambda$;
 	
 	\item[{\rm (b)}]  there is a bounded subset $\mathscr{B}\subset H$ such that $\mathscr{A}_{\lambda} \subset \mathscr{B}$ for every $\lambda \in \Lambda ;$
 	
 	\item[{\rm (c)}] for $t>0, S^{\lambda}_t x$ is continuous in $\lambda,$ uniformly for $x$ in bounded subsets of $X$.
 \end{itemize}

Then, the family $\mathscr{A}_{\lambda}$ is continuous in $\lambda$ for all $\lambda_{0}$ in a residual subset of $\Lambda .$ In particular, the set of continuity points of $\mathscr{A}_{\lambda}$ is dense in $\Lambda$.
 	
 \end{theorem}

 The next result deals with  an estimate for  Kolmogorov's $\varepsilon$-entropy ${H}_{\varepsilon}(\mathscr{M})$ of a compact set $\mathscr{M}\subset H,$ where now $(H,\|\cdot\|)$ means a Hilbert space. 
 For its proof, we refer to \cite{chueshov-lasiecka-2005}.

\begin{theorem}[{\cite[Theorem 4.2]{chueshov-lasiecka-2005}}]\label{theo-kolmog-entrop}
Let $H$ be a separable Hilbert space and $\mathscr{M}$ be a bounded closed set in $H .$ Assume that there exists a mapping $V: \mathscr{M} \mapsto H$ such that:

\begin{itemize}
	\item[$ 1. $] $\mathscr{M} \subseteq V \mathscr{M};$
	
	\item[$ 2. $]  $V$ is Lipschitz on $\mathscr{M},$ that is, there exists $L>0$ such that
	$$
	\left\|V z_{1}-V z_{2}\right\| \leq L\left\|z_{1}-z_{2}\right\|, \quad z_{1}, z_{2} \in \mathscr{M};
	$$
	
	\item[$ 3. $] There exist pseudometrics $\varrho_{1}$ and $\varrho_{2}$ on $H$ such that
	$$
	\left\|V z_{1}-V z_{2}\right\| \leq g\left(\left\|z_{1}-z_{2}\right\|\right)+h\left(\left[\varrho_{1}\left(z_{1}, z_{2}\right)^{2}+\varrho_{2}\left(V z_{1}, V z_{2}\right)^{2}\right]^{1 / 2}\right)
	$$
	for all $z_{1}, z_{2} \in \mathscr{M},$ where $g, h: \mathbb{R}^{+} \to \mathbb{R}^{+}$ are continuous non-decreasing functions such that
	$$
	g(0)=0, \ g(s)<s, \, s>0, \ s-g(s) \text { is nondecreasing,}
	$$
	and the function $h(s)$ is strictly increasing in the interval $\left[0, s_{0}\right]$ for some $s_{0}>0$ with  $h(0)=0.$

	\item[$ 4. $]For any $q>0$ and for any closed bounded set $B \subset \mathscr{M}$ the maximal number $m(B, q)$ of elements $x_{j}^{B} \in B$ such
	$$
	\varrho_{1}\left(x_{j}^{B}, x_{i}^{B}\right)^{2}+\varrho_{2}\left(V x_{j}^{B}, V x_{i}^{B}\right)^{2}>q^{2},  \ \ i \neq j,  \, i, j=1, \ldots, m(B, q),
	$$
	is finite.
\end{itemize}

Then $\mathscr{M}$ is a compact set and there exists $0<\varepsilon_{0}<1$ such that for all $\varepsilon \leq \varepsilon_{0}<1,$ Kolmogorov's $\varepsilon$-entropy ${H}_{\varepsilon}(\mathscr{M})$ admits the following estimate
$$
{H}_{\varepsilon}(\mathscr{M}) \leq \int_{\varepsilon}^{\varepsilon_{0}} \frac{\ln m\left(g_{\delta}^{-1}(s), q(s)\right)}{s-g_{\delta}(s)} d s+H_{g_{0}\left(\varepsilon_{0}\right)}(\mathscr{M}),
$$
where $g_{\delta}(s)=\frac{1-\delta}{2} g(2 s)+\delta s$ with arbitrary $\delta \in(0,1),$ the function $q(s)$ is defined by the formula
$$
q(s)=\frac{1}{2} h^{-1}\{\delta[2 s-g(2 s)]\}, 0<s<\varepsilon_{0},
$$
and
$$
m(r, q)=\sup \{m(B, q): B \subseteq M, \operatorname{diam} B \leq 2 r\}.
$$
	
\end{theorem}

\begin{theorem}[{\cite[Theorem 4.5]{chueshov-lasiecka-2005}}]\label{theo-dim-finite}
 Under the hypotheses of Theorem $\ref{theo-kolmog-entrop}$ with
 \begin{itemize}
 	\item[$(i)$] 
 	$\displaystyle\lim_{s \rightarrow 0} \frac{g(s)}{s}=g_0<1; 
 	$
 	
 	\item[$(ii)$] $h(s)$ being a linear function  ($h(s)=h_{0} \cdot s$);
 	
 	\item[$(iii)$]  $\varrho_{i}:=n_{i}, i=1,2$, being a  precompact seminorm  on $H$ (item $ 4 $ can be neglected).
 \end{itemize}
Then  $\mathscr{M}$ is a compact set in $H$ with finite fractal dimension $
\left(\mbox{dim}^{f}_{H}\mathscr{M}<\infty\right)$.
\end{theorem}
\begin{proof}
  See also \cite[Theorem 2.15]{chueshov-white}, \cite[Theorem 7.3.3]{chueshov-yellow}, or  \cite[Theorem 3.1.15]{chueshov-solo} for slightly new versions of this result.
 \end{proof}

Below we recall the notion on quasi-stable dynamical systems by following \cite{chueshov-solo,chueshov-yellow} and then some results on regularity, finite dimensionality, and  exponential attrators.

\begin{itemize}
	\item   The dynamical system $({H},S_t)$ is said to be
 {\it quasi-stable} on a set $\mathscr{B} \subset H$ (at time $t^{*}$) if there exist 
  time $t^{*}>0,$ a Banach space $Z,$ a globally Lipschitz mapping $V: \mathscr{B} \mapsto Z,$ and
a compact seminorm $n_{Z}(\cdot)$ on the space $Z,$ such that
\begin{equation}\label{quasi-stabil}
	\left\|S_{t^{*}} z_{1}-S_{t^{*}} z_{2}\right\| \leq q \cdot\left\|z_{1}-z_{2}\right\|+n_{H}\left(V z_{1}-V z_{2}\right), 
\end{equation}
	for every $z_{1}, z_{2} \in \mathscr{B}$ with $0 \leq q<1 .$ Here,  the space $Z,$ the operator
	$V,$ the seminorm $n_{Z},$ and the time moment $t^{*}$ may depend on $\mathscr{B}$.

\item Let $H$ be decomposed as  
${H}=X\times Y \times Z$ where   $X,Y,Z$ are reflexive Banach spaces with $X$ compactly embedded in $Y,$ and endowed with the usual norm. Additionally, let 
 the dynamical system $({H},S_t)$
be given by an evolution operator like
\begin{equation}\label{inforce1}
S_t z = (u(t),u_{t}(t),\zeta(t)), \quad z=(u_0, u_1, \zeta_0)
\in{H},
\end{equation}
where the functions $u$ and $\zeta$ possess the properties 
\begin{equation}\label{inforce2}
u \in C(\mathbb{R}^{+}; X)\cap C^{1}(\mathbb{R}^{+}; Y), \quad
\zeta \in C(\mathbb{R}^{+}; Z).
\end{equation}
Under this structure, one says that $({H},S_t)$ is {\it asymptotically quasi-stable} on a set $\mathscr{B}
\subset {H}$ if there exist a compact seminorm $n_X(\cdot)$ on $X$,
non-negative scalar functions $a(t)$ and $c(t)$ locally bounded in
$[0,\infty)$, and a function $b\in L^1(\mathbb{R}^{+})$ with
$\displaystyle\lim_{t\to \infty}b(t)=0$, such that
\begin{equation} \label{local-lips}
\Vert S_tz_1 - S_tz_2 \Vert^2 \le a(t)\Vert z_1 -
z_2 \Vert^2,
\end{equation}
and
\begin{equation} \label{stabili}
\Vert S_tz_1 - S_tz_2 \Vert^2 \le b(t)\Vert z_1 -
z_2 \Vert^2 + c(t) \sup_{0<s<t}\left[
n_X(u_1(s)-u_2(s)) \right]^2 ,
\end{equation}
for any $z_1,z_2 \in \mathscr{B}$,
where we denote $S_{t} z_{i}=\left(u_{i}(t), u_{t,i}(t) , \zeta_{i}(t)\right), i=1,2.$	
\end{itemize}

The next result shows that the quasi-stability notion generalizes the concept of asymptotically quasi-stability for structural systems like \eqref{inforce1}-\eqref{inforce2}. To avoid repetitions, whenever we (only) state below that $({H},S_t)$ is asymptotically quasi-stable, it is implicit that  $({H},S_t)$ is given by \eqref{inforce1}-\eqref{inforce2}.

\begin{proposition}[{\cite[Proposition 3.4.17]{chueshov-solo}}]
If $(H, S_{t})$  is asymptotically quasi-stable on some set
	$\mathscr{B}\subset H,$ then it is quasi-stable on $\mathscr{B}$ at every time $T>0$ such that $b(T)<1$.
\end{proposition}
  
   The following result is a direct consequence of  
   \cite[Proposition 3.4.3 and Corollary 3.4.4]{chueshov-solo} and
   \cite[Proposition 7.9.4 and  Corollary 7.9.5]{chueshov-yellow}. 
  \begin{proposition}\label{prop-quasi-stabil}
  	Let us assume that the dynamical system $({H},S_t)$ is dissipative and (asymptotically)
  	quasi-stable	 on every bounded forward invariant set $\mathscr{B}\subset H .$ 
  	Then, $\left(H, S_{t}\right)$ is 	asymptotically smooth and, consequently, it  possesses a compact global
  	attractor $\mathscr{A}\subset H$.
  \end{proposition}

\begin{theorem}\label{theo_fractal}
	Let us assume that the dynamical system $({H},S_t)$ 	
	possesses a compact global attractor $\mathscr{A}$ and is (asymptotically) quasi-stable on $\mathscr{A}$ at some point $t^{*}>0$. Then,  $\mathscr{A}$ has  finite fractal dimension 
	$\mbox{dim}^{f}_{H}\mathscr{A}<\infty$.
\end{theorem}  
\begin{proof}
See \cite[Theorem 7.9.6]{chueshov-yellow} for asymptotically quasi-stable	systems and  \cite[Theorem 3.4.5]{chueshov-solo} for more general quasi-stable systems.
\end{proof}

 For asymptotically  quasi-stable systems one can reach the following regularity of trajectories from the attractor.
 
\begin{theorem}[\cite{chueshov-yellow}, Theorem. 7.9.8] \label{theo_trajectories}
Let us assume that  $({H},S_t)$ satisfies the structure \eqref{inforce1}-\eqref{inforce2}, 
possesses a compact global attractor $\mathscr{A}$ and is asymptotically quasi-stable on $\mathscr{A}$. 	 Additionally if \eqref{stabili} holds with  $c(t)$ satisfying $c_{\infty}=\sup_{t\in\mathbb{R}^+}c(t)<\infty,$ then any
	full trajectory $\Gamma=\{(u(t), u_t(t), \zeta(t));  \
	t\in\mathbb{R} \} \subset\mathscr{A}$ enjoys the following regularity properties
	$$
	u_t\in  L^{\infty}(\mathbb{R}; X)\cap C(\mathbb{R},Y), \quad
	u_{tt} \in  L^{\infty}(\mathbb{R}; Y), \quad \zeta_t\in
	L^{\infty}(\mathbb{R}; Z).
	$$
Besides, there exists a constant $R > 0$ such that
	$$
	\sup_{\Gamma\subset\mathscr{A}}\sup_{t\in\mathbb{R}}\big(
	\|u_t(t)\|^2_X + \|u_{tt}(t)\|^2_Y + \|\zeta_t(t)\|_Z^2 \big)\leq
	R^2,
	$$
	where $R$ depends on the constant $c_{\infty}$, on the seminorm
	$n_X(\cdot)$, and  on the embedding  $X\hookrightarrow Y$.
\end{theorem}

Generalized fractal exponential attractors can be also reached for quasi-stable and asymptotic quasi-stable systems as follows. The next version can be found in  \cite[Theorem 7.9.9]{chueshov-yellow}. See also 
\cite[Theorem 3.4.7]{chueshov-solo}.

\begin{theorem}\label{theo_attrac_exponential}
Let us assume that the dynamical system $({H},S_t)$ is dissipative and 
asymptotically quasi-stable	 on some bounded absorbing set $\mathscr{B}$. In
	addition, let us suppose that there exists a space
	$\widetilde{H}\supseteq H$ such that mapping $t\mapsto S_tz$ is
	H\"{o}lder continuous in $\widetilde{H}$ for each
	$z\in\mathscr{B},$ that is, there exist $0<\sigma\leq1$ and
	$C_{\mathscr{B},T}>0$ ($T>0$ given) such that
	\begin{equation}\label{holder_prop}
	\| S_{t_2}z-S_{t_1}z\|_{\widetilde{H}} \, \leq \,
	C_{\mathscr{B},T}|t_2-t_1|^{\sigma}, \quad t_1,t_2\in[0,T], \ z
	\in\mathscr{B}.
	\end{equation}
	Then, $({H},S_t)$ possesses a generalized
	fractal exponential attractor $\mathscr{A}_{\exp}$ whose dimension is finite in the
	space $\widetilde{H}$ $(\mbox{dim}^{f}_{\widetilde{H}}\mathscr{A}_{\exp}<\infty)$.
\end{theorem}

In Theorem \ref{theo_attrac_exponential}, unless $\widetilde{H}=H$ we can only guarantee the finite fractal dimension in an extended phase space $\widetilde{H}$. Thus, to achieve such finiteness of fractal dimensional in $H$, one must prove \eqref{holder_prop} in $H$, which sometimes seems to be a hard task. Therefore, in order  to present a tangible result with exponential attractor whose fractal dimensional is finite in $H$ 
($\mbox{dim}^{f}_{H}\mathscr{A}_{\exp}<\infty$), we finally remind an useful result by following \cite{carvalho-sonner1,carvalho-sonner2}, which relies on the construction of time-dependent exponential attractors $\mathscr{A}_{\exp}=\left\{\mathscr{A}_{\exp}(t); \, t \in \mathbb{R}\right\}$ for (continuous) dynamical systems $({H},S_t)$ under suitable 
 decomposition and Lipschitz properties. For previous results on the subject, we also refer to \cite{eden-etal,efend-etal,efend-etal2,miranv-zelik,temam1}. The next concept and result are based on the more recent construction developed in \cite[Section 4]{carvalho-sonner1}.

\begin{itemize}
	\item  The familiy $\mathscr{A}_{\exp}=\left\{\mathscr{A}_{\exp}(t); \, t \in \mathbb{R}\right\}$ is called a {\it time-dependent exponential attractor} for  $({H},S_t)$  if there exists
	$0<a<\infty$ such that  $\mathscr{A}_{\exp}(t)=\mathscr{A}_{\exp}(a+t)$ for all $t \in \mathbb{R}$, and 
	
	(i)  the subsets $\mathscr{A}_{\exp}(t) \subset H$ are non-empty and compact in $H$ for all $t \in \mathbb{R}$;
	
	(ii)  the family is positively semi-invariant, that is
	$$
	S_t \mathscr{A}_{\exp}(s) \subset \mathscr{A}_{\exp}(t+s), \quad \forall  \,t \geq 0, s \,\in \mathbb{R};
	$$
	
  (iii) the fractal dimension of the sets $\mathscr{A}_{\exp}(t), \, t \in \mathbb{R},$ is uniformly bounded;
	
	(iv)  the family  attracts exponentially all bounded subsets of $H$ uniformly, that is, there exists a positive constant $\omega>0$ such that for any bounded subset $\mathscr{B} \subset H$
	$$
	\lim _{\tau \rightarrow \infty} \sup _{t \in[0, a]} e^{\omega \tau} \operatorname{dist}_{\mathrm{H}}(S_\tau\mathscr{B}, \mathscr{A}_{\exp}(t))=0.
	$$
\end{itemize}

\begin{theorem}[{\cite[Theorem 4.4]{carvalho-sonner1}}]\label{theo-carvalho-soner}
Let us assume that the dynamical system $({H},S_t)$ can be splitted into $S_t=S^1_t+S^2_t:H\to H$ and 
let 
 $W$ be another normed space with compact embedding $(H,\|\cdot\|_H )\hookrightarrow \left(W,\|\cdot\|_{W}\right)$.
  Let us  additionally  suppose the following conditions:
 \begin{itemize}
 	\item[{\rm (S1)}]  $({H},S_t)$ is dissipative, that is, it has a bounded absorbing set $\mathscr{B}\subset H;$
 	
 	\item[{\rm (S2)}] there exist a constant  $0 \leq c_1<\frac{1}{2}$ and a time $T>0$ such that $S^1_t$ satisfies the contraction property on $\mathscr{B}$
 	$$
 	\|S^1_{T}z_1 -S^1_{T} z_2 \|_{H} \leq c_1 \|z_1-z_2\|_{H}, \quad    z_1 ,z_2 \in \mathscr{B};
 	$$

 	\item[{\rm (S3)}]  there exists 
 	a constant  $c_2>0$
 	such that
 	$S^2_t$ satisfies  smoothing property within $\mathscr{B}$ at time $T>0$
 	$$
 	\|S^2_{T}z_1 -S^2_{T} z_2\|_{H} \leq  c_2\|z_1-z_2\|_{W}, \quad    z_1,z_2 \in \mathscr{B};
 	$$

 	\item[{\rm (S4)}] there exists a time $T_0\geq0$ such that $S_t$ is  Lipschitz on  $\mathscr{B}$ for $t\geq T_0$, that is, for some constant $L_{t}>0$  it holds
 $$
 \|S_t z_1 - S_t z_2\|_{H} \leq L_{t} \|z_1-z_2\|_{H}, \quad    z_1,z_2 \in \mathscr{B}, \ t\geq T_0.
 $$	
 \end{itemize}
Then, $({H},S_t)$ possesses a
time-dependent exponential attractor $\mathscr{A}_{\exp}=\left\{\mathscr{A}_{\exp}(t); \, t \in \mathbb{R}\right\}$, whose 
sections are compact subsets of $H$  with finite fractal dimension in $H$, that is, 
$$
\mbox{dim}^{f}_{H}\big(\mathscr{A}_{\exp}(t)\big)<\infty, \quad \forall \, t \in \mathbb{R}.
$$
\end{theorem}

\begin{remark}\label{remark-last}
Finally, according to \cite[Remark 6]{carvalho-sonner1}, under the assumptions of Theorem $ \ref{theo-carvalho-soner} $ one can construct an  exponential attractor for the dynamical system $({H},S_t)$. Indeed, to this end it is enough to consider (in general) the union
$$
\widetilde{\mathscr{A}_{\exp}}:=\bigcup_{t \in[T, 2 T]} S_t \mathscr{A}_{\exp}^{d},
$$
where $\mathscr{A}_{\exp}^{d}$ denotes the exponential attractor for the corresponding discrete semigroup
$
\{S_{n T}\}_{n \in \mathbb{N}}.
$
\end{remark}

\paragraph*{Acknowledgments.}

The authors would like to express their gratitude to Professors Alain Miranville, Irena Lasiecka,  Vittorino Pata, To Fu Ma, and Igor Chueshov (in memory) for several helpful advices on quasi-stability, global and exponential attractors.

 \addcontentsline{toc}{section}{\numberline{}References}

\end{document}